\documentclass[12pt]{article}
\usepackage{amsmath,graphicx,amssymb,amsthm, latexsym}
\allowdisplaybreaks
\def\cc{\mathbb{C}}

\def\nn{{\bf N}}

\def\GG{\mathfrak G}
\def\FF{\mathfrak F}

\def\NN{\mathfrak N}
\def\A{\mathcal A}
\def\B{\mathcal B}

\def\F{\mathcal F}

\def\H{\mathcal H}

\def\L{\mathcal L}
\def\M{\mathcal M}
\def\N{\mathcal N}
\def\P{\mathcal P}

\def\S{\mathcal S}
\def\U{\mathcal U}

\def\amslatex{$\mathcal{A}\kern-.1667em\lower.5ex\hbox{$\mathcal{M}$}\kern-.125em\mathcal{S}$-\LaTeX}

\DeclareMathOperator{\esssup}{ess\, sup}
\DeclareMathOperator{\essinf}{ess\, inf}
\DeclareMathOperator{\Aut}{Aut}

\def \sp {\hspace{0.3cm}}

\newtheorem{set}{set}[section]
\newtheorem{Corollary}[set]{Corollary}
\newtheorem{Example}[set]{Example}
\newtheorem{Definition}[set]{Definition}

\newtheorem{Lemma}[set]{Lemma}

\newtheorem{Proposition}[set]{Proposition}
\newtheorem{Remark}[set]{Remark}
\newtheorem{Theorem}[set]{Theorem}
\newcommand{\define}{\mathrel{\hbox{$\equiv$\hskip -.90em \lower .47ex \hbox{$\leftharpoondown$}}}}
\newcommand{\enifed}{\mathrel{\hbox{$\equiv$\hskip -.90em \lower .47ex \hbox{$\rightharpoondown$}}}}

\numberwithin{equation}{section}
\begin{document}
\date{}
\title {Tracial gauge norms on finite von Neumann algebras satisfying the
weak Dixmier property}
\author{Junsheng Fang\thanks{jfang\@@cisunix. unh.edu}\, Don Hadwin\thanks{don\@@math.unh.edu}
 \, Eric Nordgren\thanks{ean@christa.unh.edu} \, Junhao Shen\thanks{junhao.shen@unh.edu}\\
Department of Mathematics\\
 University of New Hampshire\\
 Durham, NH  03824}
\maketitle
\begin{abstract}  In this paper we set up a representation theorem
for tracial gauge norms on finite von Neumann algebras satisfying
the weak Dixmier property  in terms of Ky Fan  norms. Examples of
tracial gauge norms on finite von Neumann algebras satisfying the
weak Dixmier property include unitarily invariant norms on finite
factors (type ${\rm II}\sb 1$ factors and $M_n(\cc)$) and symmetric
gauge norms on $L^\infty[0,1]$ and $\cc^n$. As the first
application, we  obtain that the class of unitarily invariant norms
on a type ${\rm II}\sb 1$ factor coincides with the class of
symmetric gauge norms on $L^\infty[0,1]$ and von Neumann's classical
result~\cite{vN} on  unitarily invariant norms on $M_n(\cc)$.  As
the second application, Ky Fan's dominance theorem~\cite{Fan} is
obtained for finite von Neumann algebras satisfying the weak Dixmier
property. As the third application, some classical results in
non-commutative $L^p$-theory (e.g., non-commutative
H$\ddot{\text{o}}$lder's inequality, duality and reflexivity of
non-commutative $L^p$-spaces) are obtained for general unitarily
invariant norms on finite factors.  We also investigate the extreme
points of $\NN(\M)$, the convex compact set (in the pointwise weak
topology) of normalized unitarily invariant norms (the norm of the
identity operator is 1) on a finite factor $\M$. We obtain all
extreme points of $\NN(M_2(\cc))$ and some extreme points of
$\NN(M_n(\cc))$ ($n\geq 3$). For a type ${\rm II}\sb 1$ factor $\M$,
we prove that if $t$ ($0\leq t\leq 1$) is a rational number then the
Ky Fan $t$-th norm is an extreme point of $\NN(\M)$.
\end{abstract}

{\bf Keywords:} finite von Neumann algebras, the weak Dixmier property,
tracial gauge norms, $s$-numbers, Ky Fan norms.\\

 \vskip1.0cm

 \centerline{\Large{$\mathbf{Contents}$}}
\vskip 0.5cm

\noindent 1.\,\, Introduction\\
\noindent 2.\,\, Preliminaries\\
\indent 2.1.\,\, Nonincreasing rearrangements of functions\\
\indent 2.2.\,\, Invertible measure-preserving transformations on $[0, 1]$\\
\indent 2.3.\,\, $s-$numbers of operators in type ${\rm II}\sb 1$ factors\\
\indent 2.4.\,\, $s-$numbers of operators in  finite von Neumann
algebras\\
\noindent 3.\,\, Tracial gauge semi-norms on finite von
Neumann algebras satisfying the weak Dixmier\\\indent property\\
\indent 3.1.\,\, Gauge semi-norms\\
\indent 3.2.\,\, Tracial gauge semi-norms \\
\indent 3.3.\,\, Symmetric gauge semi-norms \\
\indent 3.4.\,\, Unitarily invariant semi-norms \\
\indent 3.5.\,\, Weak Dixmier property\\
\indent 3.6.\,\, A comparison theorem\\
\noindent 4.\,\, Proof of  Theorem~{\bf B}\\
\noindent 5.\,\, Ky Fan norms\\
\noindent 6.\,\, Dual norms of tracial gauge norms on finite von
Neumann algebras satisfying the\\\indent weak Dixmier property\\
\indent 6.1.\,\, Dual norms\\
\indent 6.2.\,\, Dual norms of Ky Fan norms\\
\indent 6.3.\,\, Proof of Theorem~{\bf C}\\
\noindent 7.\,\, Proof of Theorem~{\bf A}\\
\noindent 8.\,\, Proof of Theorem~{\bf D} and Theorem~{\bf E}\\
\noindent 9.\,\, Proof of Theorem~{\bf F}\\
\noindent 10.\,\,Extreme points of normalized unitarily invariant
norms on finite factors\\
\indent 10.1\,\, $\NN(\M)$\\
\indent 10.2\,\, $\NN_e(M_n(\cc))$\\
\indent 10.3\,\, Proof of Theorem~{\bf J}\\
\indent 10.4\,\, Proof of Theorem~{\bf K}\\
\noindent 11.\,\,Proof of Theorem~{\bf G}\\
\noindent 12.\,\,Completion of type ${\rm II}\sb 1$ factors with
respect to
unitarily invariant norms\\
\indent 12.1\,\, Embedding of $\overline{\M_{|\!|\!|\cdot|\!|\!|}}$
into $\widetilde{\M}$\\
\indent 12.2\,\, $\widetilde{\M}$ and $L^1(\M,\tau)$\\
\indent 12.3\,\, Elements in $\overline{\M_{|\!|\!|\cdot|\!|\!|}}$\\
\indent 12.4\,\, A generalization of H$\ddot{\text{o}}$lder's inequality\\
\noindent 13.\,\,Proof of Theorem~{\bf H} and Theorem~{\bf I}\\
\noindent References

\vskip 1.00cm

\section{Introduction}

The unitarily invariant norms were introduced by von
Neumann~\cite{vN} for the purpose of metrizing matrix spaces. Von
Neumann, together with his associates, established that the class of
unitarily invariant norms of $n\times n$ complex matrices  coincides
with the class of symmetric gauge functions of their $s$-numbers.
These norms have now been variously generalized and utilized in
several contexts. For example, Schatten~\cite{Sc1,Sc2} defined norms
on two-sided ideals of completely continuous  operators on an
arbitrary Hilbert space; Ky Fan~\cite{Fan} studied Ky Fan norms and
obtained his dominance theorem. The unitarily invariant norms play a
crucial role in the study of function spaces and group
representations (see e.g.~\cite{Ku}) and in obtaining certain bounds
of importance in quantum field theory (see~\cite{Si}). For
historical perspectives and surveys of unitarily invariant norms,
see Schatten~\cite{Sc1,Sc2}, Hewitt and Ross~\cite{H-R}, Gohberg and
Krein~\cite{G-K} and Simon~\cite{Si}.\\

The theory of non-commutative $L^p$-spaces has been developed under
the name ``non-commutative integration" beginning with pioneer work
of Segal, Dixmier, and Kunze. Since then the theory has been
extensively studied, extended and applied  by Nelson, Haagerup,
Fack, Kosaki, Junge, Xu, and many others. The recent survey by
Pisier and Xu~\cite{P-X} presents a rather complete picture on
noncommutative integration and contains a lot of references.  This
theory is still a very active subject of investigation. Some tools
in the study of the usual commutative $L^p$-spaces still work in the
noncommutative setting. However, most of the time, new techniques
must be invented. To illustrate the difficulties one may encounter
in studying the noncommutative $L^p$-spaces, we mention here one
basic well-known fact. Let $\H$ be a complex Hilbert space, and let
$\B(\H)$ denote the algebra of all bounded linear operators on $\H$.
The basic fact states that the usual triangle inequality for the
absolute values of complex numbers is no longer valid for the
absolute values of operators, namely, in general, we do not have
$|S+T| \leq  |S|+|T|$ for $S, T\in \B(\H)$, where $|S| =
(S^*S)^{1/2}$ is the absolute value of $S$. Despite such
difficulties, by now the strong parallelism between noncommutative
and classical Lebesgue integration is well-known. \\

Motivated by von Neumann's theorem and the analogies between
noncommutative and classical $L^p$-spaces,
 in this paper, we will  systematically study tracial gauge norms
on finite von Neumann algebras that satisfy the weak Dixmier
property. Before stating the main theorem and its consequences, we
explain some of the notation and terminology that will be used throughout the paper.\\

In this paper, a finite von Neumann algebra $(\M,\tau)$ means a  von
Neumann algebra $\M$ with a faithful normal tracial state $\tau$. A
finite von Neumann algebra $(\M,\tau)$ is said to satisfy the
\emph{weak Dixmier property} if for every positive operator $T\in
\M$, $\tau(T)$ is in the \emph{operator norm} closure of the convex
hull of $\{S\in \M:\,\, \text{$S$ and $T$ are equi-measurable,}\,
\text{i.e.}, \tau(S^n)=\tau(T^n)\, \text{for all
$n=0,1,2,\cdots$}\}$. Recall that finite factors satisfy the
\emph{Dixmier property}: if $T\in \M$, then $\tau(T)$ is in the
\emph{operator norm} closure of the convex hull of $\{UTU^*:\ U\in\M
\,\,\text{is a unitary operator}\}$ and hence satisfy the weak
Dixmier property. In section 3.5, we prove that a finite von Neumann
algebra $(\M,\tau)$
 satisfies the weak Dixmier property if
and only if either  $(\M,\tau)$ can be identified as a von Neumann
subalgebra of $(M_n(\cc),\tau_n)$ that contains all diagonal
matrices, where $\tau_n$ is the normalized trace on $M_n(\cc)$, or
$\M$ is diffuse. Throughout the paper, we will reserve the notation
$\|\cdot\|$ for the \emph{operator norm} on von Neumann algebras.\\

A \emph{tracial gauge} norm $|\!|\!|\cdot|\!|\!|$ on a finite von
Neumann algebra $(\M,\tau)$ is a norm on $\M$ satisfying
 $|\!|\!|T|\!|\!|=|\!|\!|\,|T|\,|\!|\!|$ for
all $T\in \M$ (gauge invariant) and
$|\!|\!|S|\!|\!|=|\!|\!|T|\!|\!|$ if $S$ and $T$ are two
equi-measurable positive operators in $\M$ (tracial). For a finite
von Neumann algebra $(\M,\tau)$, let $\Aut(\M,\tau)$ be the set of
$\ast$-automorphisms on $\M$ that preserve the trace. A
\emph{symmetric gauge}
 norm $|\!|\!|\cdot|\!|\!|$ on a
finite von Neumann algebra $(\M,\tau)$ is a gauge norm on $\M$
satisfying $|\!|\!|\theta(T)|\!|\!|=|\!|\!|T|\!|\!|$ for all
positive operators $T\in \M$ and $\theta\in \Aut(\M,\tau)$.  A
\emph{unitarily invariant} norm $|\!|\!|\cdot|\!|\!|$ on a finite
von Neumann algebra $(\M,\tau)$ is a norm on $\M$ satisfying
$|\!|\!|UTW|\!|\!|=|\!|\!|T|\!|\!|$ for all $T\in \M$ and unitary
operators $U, W$ in $\M$. On $(L^\infty[0,1],\int_0^1\,dx)$ and
$(\cc^n,\tau)$, where
$\tau((x_1,\cdots,x_n))=\frac{x_1+\cdots+x_n}{n}$, a norm  is a
tracial gauge norm if and only if it is a symmetric gauge norm. A
norm on a finite factor is a tracial gauge norm if and only if it is
a unitarily invariant norm.
 A \emph{normalized} norm is one that assigns the value 1 to the identity
 operator (which is also denoted by 1).\\

In~\cite{F-K}, Fack and Kosaki defined $\mu_s(T)$, \emph{the
generalized $s$-numbers of an operator $T$} in a finite von
Neumann algebra $(\M,\tau)$, by
\[\mu_s(T)=\inf\{\|TE\|: E\,\, \text{is a projection in $\M$ with $\tau(1- E) \leq
s$}\},\quad 0\leq s \leq 1.
\]
For $0<t\leq 1$, \emph{the Ky Fan $t$-th  norm},
$|\!|\!|T|\!|\!|_{(t)}$, on a finite von Neumann algebra $(\M,\tau)$
is defined  by
\[|\!|\!|T|\!|\!|_{(t)}=\frac{1}{t}\int_0^t\mu_s(T)ds.
\]
Then $|\!|\!|\cdot|\!|\!|_{(t)}$ is a tracial gauge norm on
$(\M,\tau)$.  Note that $|\!|\!|T|\!|\!|_{(1)}=\tau(|T|)=\|T\|_1$ is
the trace norm. \\

Let $n\in \nn$, $a_1\geq a_2\geq \cdots\geq a_n\geq a_{n+1}=0$ and
$f(x)=a_1\chi_{[0,\frac{1}{n})}(x)+a_2\chi_{[\frac{1}{n},\frac{2}{n})}(x)+\cdots+
 a_n\chi_{[\frac{n-1}{n},1]}(x)$. For $T\in \M$, define
$|\!|\!|T|\!|\!|_{f}=\int_0^1 f(s)\mu_s(T)ds.$
  Then \[
|\!|\!|T|\!|\!|_{f}=\sum_{k=1}^n \frac{k(a_k - a_{k+1})}{n}
|\!|\!|T|\!|\!|_{\left(\frac{k}{n}\right)}.\] Therefore,
$|\!|\!|T|\!|\!|_{f}$ is a tracial gauge norm on $(\M,\tau)$.  Note
that if $f(x)$ is  the constant $1$ function on $[0,1]$,
then $|\!|\!|T|\!|\!|_{f}=|\!|\!|T|\!|\!|_{(1)}=\|T\|_1=\tau(|T|)$. \\

Let
$\F=\{f(x)=a_1\chi_{[0,\frac{1}{n})}(x)+a_2\chi_{[\frac{1}{n},\frac{2}{n})}(x)+\cdots+
 a_n\chi_{[\frac{n-1}{n},1]}(x):\, a_1\geq a_2\geq \cdots\geq
 a_n\geq 0, \frac{a_1+\cdots+a_n}{n}\leq 1, n=1,2,\cdots\}$.  In  section 7,
we prove the following representation theorem, which is the main
result of this paper.\\

\noindent {\bf Theorem A.}\,\,\emph{Let $(\M,\tau)$ be a finite von
Neumann algebra satisfying the weak Dixmier property. If
 $|\!|\!|\cdot|\!|\!|$ is a normalized tracial gauge norm on
$\M$, then there is a subset $\F'$ of $\F$ containing the constant
$1$ function on $[0,1]$ such that for every $T\in \M$,
\[|\!|\!|T|\!|\!|=\sup\{|\!|\!|T|\!|\!|_{f}:\, f\in
\F'\},\] where $|\!|\!|T|\!|\!|_{f}$ is defined as above.}
\bigskip

To prove Theorem~{\bf A}, we firstly prove the following technical
theorem in section 4. \\

\noindent {\bf Theorem B.}\,\,\emph{Let $(\M,\tau)$ be a finite
von Neumann algebra satisfying the weak Dixmier property and let
$|\!|\!|\cdot|\!|\!|$ be a tracial gauge norm on $\M$. Then
$\M_{1,\,|\!|\!|\cdot|\!|\!|}=\{T\in\M:\, |\!|\!|T|\!|\!|\leq 1\}$
is closed in the weak operator topology.}
\bigskip

The Russo-Dye Theorem~\cite{R-D} and the Kadison-Peterson
Theorem~\cite{K-P} on convex hulls of unitary operators in von
Neumann algebras  and the idea of Dixmier's averaging
Theorem~\cite{Di1} play  fundamental roles in the proof of
Theorem~{\bf B}. An important consequence of Theorem~{\bf B} is the
following corollary which enables us to apply the powerful
techniques of normal conditional expectations from finite von
Neumann algebras to abelian von Neumann
subalgebras.\\

\noindent {\bf Corollary 1.}\,\,\emph{Let $(\M,\tau)$ be a finite
von Neumann algebra satisfying  the weak Dixmier property and let
$|\!|\!|\cdot|\!|\!|$ be a tracial gauge  norm on $\M$. If $\A$ is
a separable  abelian von Neumann subalgebra of $\M$ and
$\mathbf{E}_\A$ is the normal conditional expectation from $\M$
onto $\A$ preserving $\tau$, then
$|\!|\!|\mathbf{E}_{\A}(T)|\!|\!|\leq |\!|\!|T|\!|\!|$ for all
$T\in \M$. }
\bigskip

The notion of \emph{dual norms} plays a key role in the proof of
Theorem~{\bf A}. Let $|\!|\!|\cdot|\!|\!|$ be a norm on a finite von
Neumann algebra $(\M,\tau)$. Then the dual norm
$|\!|\!|\cdot|\!|\!|^\#$ is defined by
\[|\!|\!|T|\!|\!|^\#=\sup\{|\tau(TX)|:\, X\in\M,\, |\!|\!|X|\!|\!|\leq
1\},\quad T\in \M.
\] In section 5,  we study the dual norms systematically. By
applying  Corollary 1  and careful analysis, we prove
the following theorem. \\

\noindent {\bf Theorem C.}\,\,\emph{Let $(\M,\tau)$ be a finite von
Neumann algebra satisfying the weak Dixmier property and
$|\!|\!|\cdot|\!|\!|$ be a tracial gauge norm on $\M$. Then
$|\!|\!|\cdot|\!|\!|^\#$ is also a tracial gauge norm on $\M$ and
$|\!|\!|\cdot|\!|\!|^{\#\#}=|\!|\!|\cdot|\!|\!|$.}
\bigskip

Combining  Corollary~{\bf 1}, Theorem~{\bf C} and the following
theorem on non-increasing rearrangements of functions (see 10.13
of~\cite{HLP}
for instance), we prove Theorem~{\bf A} in section 7.\\

\noindent {\bf
Hardy-Littlewood-P$\acute{\text{o}}$lya.}\,\,\emph{Let $f(x),
g(x)$ be non-negative Lebesgue measurable functions on $[0,1]$ and
let $f^*(x), g^*(x)$ be the non-increasing rearrangements of
$f(x),g(x)$, respectively, then $\int_0^1 f(x)g(x)dx\leq \int_0^1
f^*(x)g^*(x)dx$. }
\bigskip

Now we state some important consequences of Theorem~{\bf A}.  Since
there is a natural  one-to-one correspondence between Ky Fan $t$-th
norms on  finite von Neumann algebras (satisfying the weak Dixmier
property) and Ky Fan $t$-th norms on $(L^\infty[0,1],\int_0^1 dx)$
or $(\cc^n,\tau)$, the first
application of  Theorem~{\bf A} is the following\\

\noindent{\bf Theorem D.}\,\,\emph{Let $(\M,\tau)$ be a diffuse
finite von Neumann algebra $($or a von Neumann subalgebra of
$\M_n(\cc)$, $\tau=\tau_n|_\M$, such that $\M$ contains all diagonal
matrices$)$. Then there is a one-to-one correspondence between
tracial gauge norms on $(\M, \tau)$ and symmetric gauge norms on
$(L^\infty[0,1], \int_0^1dx)$ $($or $(\cc^n, \tau')$,
$\tau'((x_1,\cdots,x_n))=\frac{x_1+\cdots+x_n}{n}$, respectively
$)$. Namely:
\begin{enumerate}
\item if $|\!|\!|\cdot|\!|\!|$ is a  tracial gauge norm on
$(\M,\tau)$ and $\theta$ is an embedding from $(L^\infty[0,1],
\int_0^1 dx)$ into $(\M,\tau)$ $($or $x_1\oplus\cdots\oplus x_n$ is
the diagonal matrix with diagonal elements $x_1,\cdots,x_n$,
respectively$)$, then
$|\!|\!|f(x)|\!|\!|'=|\!|\!|\theta(f(x))|\!|\!|$ defines a symmetric
gauge norm on $(L^\infty[0,1],\\ \int_0^1dx)$ $($or
$|\!|\!|(x_1,\cdots,x_n)|\!|\!|'=\left|\!\left|\!\left|x_1\oplus\cdots\oplus
x_n\right|\!\right|\!\right|$ defines a symmetric gauge norm on
$(\cc^n,\tau')$, respectively$)$;
\item if $|\!|\!|\cdot|\!|\!|'$ is a  symmetric gauge
norm on $(L^\infty[0,1], \int_0^1dx)$ $($or $(\cc^n,\tau')$
respectively$)$, then $|\!|\!|T|\!|\!|=|\!|\!|\mu_s(T)|\!|\!|'$
$($or $|\!|\!|T|\!|\!|=|\!|\!|(s_1(T),\cdots,s_n(T))|\!|\!|'$
respectively$)$ defines a tracial gauge norm on $(\M,\tau)$.
\end{enumerate}}
\bigskip

As consequences of Theorem~{\bf D}, we have the following corollary
and
von Neumann's Theorem.\\

\noindent {\bf Corollary 2.}\,\,\emph{There is a one-to-one
correspondence between  unitarily invariant norms on a type ${\rm
II}\sb 1$ factor $(\M, \tau)$ and  symmetric gauge norms on
$(L^\infty[0,1], \int_0^1dx)$. Namely:
\begin{enumerate}
\item if $|\!|\!|\cdot|\!|\!|$ is a  unitarily invariant norm on
$\M$ and  $\theta$ is an embedding from $(L^\infty[0,1]$,
$\int_0^1dx)$ into $(\M,\tau)$, then
$|\!|\!|f(x)|\!|\!|'=|\!|\!|\theta(f(x))|\!|\!|$ defines a symmetric
gauge norm on $(L^\infty[0,1], \int_0^1dx)$;
\item if $|\!|\!|\cdot|\!|\!|'$ is a symmetric gauge
norm on $(L^\infty[0,1], \int_0^1dx)$, then
$|\!|\!|T|\!|\!|=|\!|\!|\mu_s(T)|\!|\!|'$ defines a  unitarily
invariant norm on $\M$.
\end{enumerate}
}
\bigskip

\noindent {\bf Von Neumann.}\,\,\emph{There is a one-to-one
correspondence between unitarily invariant norms on $M_n(\cc)$ and
symmetric gauge norms on $(\cc^n,\tau)$,
$\tau((x_1,\cdots,x_n))=\frac{x_1+\cdots+x_n}{n}$. Namely:
\begin{enumerate}
\item if $|\!|\!|\cdot|\!|\!|$ is a unitarily invariant norm on
$M_n(\cc)$, then
$|\!|\!|(x_1,\cdots,x_n)|\!|\!|'=|\!|\!|x_1\oplus\cdots\oplus
x_n|\!|\!|$ defines a symmetric gauge norm on $(\cc^n,\tau)$;
\item if $|\!|\!|\cdot|\!|\!|'$ is a  symmetric gauge
norm on $(\cc^n,\tau)$, then
$|\!|\!|T|\!|\!|=|\!|\!|(s_1(T),\cdots,s_n(T))|\!|\!|'$ defines a
unitarily invariant norm on $M_n(\cc)$.
\end{enumerate}
}
\bigskip

Theorem~{\bf D}  establishes  the one to one correspondence between
tracial gauge norms on finite von Neumann algebras satisfying the
weak Dixmier property and symmetric gauge norms on abelian von
Neumann algebras. The following theorem further establishes the one
to one correspondence between the dual norms on finite von Neumann
algebras satisfying the weak Dixmier property and the dual norms on
abelian von Neumann algebras, which plays a key role in the studying
of duality and reflexivity of the completion of type ${\rm II}\sb 1$
factors
with respect to  unitarily invariant norms. \\

\noindent {\bf Theorem E.}\,\,\emph{Let $(\M,\tau)$ be a diffuse
finite von Neumann algebra $($or a von Neumann subalgebra of
$\M_n(\cc)$, $\tau=\tau_n|_\M$, such that $\M$ contains all diagonal
matrices$)$. If $|\!|\!|\cdot|\!|\!|$ is a tracial gauge norm on
$(\M,\tau)$ corresponding to the symmetric gauge norm
$|\!|\!|\cdot|\!|\!|_1$ on $(L^\infty[0,1],\int_0^1dx)$ $($or
$(\cc^n, \tau')$ respectively$)$ as in Theorem~{\bf D}, then
$|\!|\!|\cdot|\!|\!|^\#$ on $\M$ is the tracial gauge norm
corresponding to the symmetric gauge norm $|\!|\!|\cdot|\!|\!|_1^\#$
on $(L^\infty[0,1],\int_0^1dx)$ $($or $(\cc^n, \tau')$
respectively$)$ as in Theorem~{\bf D}.}
\bigskip

The second consequence of Theorem~{\bf A} is the following theorem.\\

 \noindent {\bf Theorem F.}\,\,\emph{Let
$(\M,\tau)$ be a finite von Neumann algebra satisfying the weak
Dixmier property
 and $S, T\in \M$. If $|\!|\!|S|\!|\!|_{(t)}\leq |\!|\!|T|\!|\!|_{(t)}$
 for all Ky Fan
$t$-th norms, $0\leq t\leq 1$, then $|\!|\!|S|\!|\!|\leq
|\!|\!|T|\!|\!|$ for all tracial gauge norms $|\!|\!|\cdot|\!|\!|$
on $\M$.}
\bigskip

As a  corollary, we obtain the following\\

\noindent{\bf Ky Fan's Dominance Theorem~\cite{Fan}.}\,\, \emph{If
$S,T\in M_n(\cc)$ and $|\!|\!|S|\!|\!|_{(k/n)}\leq
|\!|\!|T|\!|\!|_{(k/n)}$, i.e., $\sum_{i=1}^k s_i(S)\leq
\sum_{i=1}^k s_i(T)$ for $1\leq k\leq n$, then $|\!|\!|S|\!|\!|\leq
|\!|\!|T|\!|\!|$ for all unitarily invariant norms
$|\!|\!|\cdot|\!|\!|$ on $M_n(\cc)$.}
\bigskip

 A unitarily invariant norm $|\!|\!|\cdot|\!|\!|$ on a type ${\rm II}\sb 1$
factor $\M$ is called \emph{singular} if $\lim_{\tau(E)\rightarrow
0+}|\!|\!|E|\!|\!|$ $>0$ and \emph{continuous} if
$\lim_{\tau(E)\rightarrow 0+}|\!|\!|E|\!|\!|=0$. The following
theorem is proved in section
11.\\

\noindent{\bf Theorem G.}\,\, \emph{ Let $|\!|\!|\cdot|\!|\!|$ be
a unitarily invariant norm on $\M$ and let $\mathcal{T}$ be the
topology induced by $|\!|\!|\cdot|\!|\!|$ on
$\M_{1,\,\|\cdot\|}=\{T\in\M:\, \|T\|\leq 1\}$. If
$|\!|\!|\cdot|\!|\!|$ is singular, then $\mathcal{T}$ is the
operator norm topology on $\M_{1,\,\|\cdot\|}$. If
$|\!|\!|\cdot|\!|\!|$ is continuous, then $\mathcal{T}$ is the
measure topology $($in the sense of Nelson~}\cite{Ne}) \emph{on
$\M_{1,\,\|\cdot\|}$.}
\bigskip

Let $\M$ be a type ${\rm II}\sb 1$ factor and $|\!|\!|\cdot|\!|\!|$
be a unitarily invariant norm on $\M$. We denote by
$\overline{\M_{|\!|\!|\cdot|\!|\!|}}$ the completion of $\M$ with
respect to $|\!|\!|\cdot|\!|\!|$. Let $\widetilde{\M}$ be the
completion of $\M$ with respect to the measure topology in the sense
of Nelson~\cite{Ne}. In section 12, we prove that there is an
injective map from $\overline{\M_{|\!|\!|\cdot|\!|\!|}}$ into
$\widetilde{\M}$ that extends the identity map from $\M$ onto $\M$.
An element in $\widetilde{\M}$ can be identified with a closed,
densely defined operator affiliated with $\M$ (see~\cite{Ne}). So
generally speaking, an element in
$\overline{\M_{|\!|\!|\cdot|\!|\!|}}$ should be treated as an
unbounded operator. We will consider the
following two questions in section 13:\\

\noindent Question 1: Under what conditions is
$\overline{{\M}_{|\!|\!|\cdot|\!|\!|^\#}}$  the dual space of
$\overline{{\M}_{|\!|\!|\cdot|\!|\!|}}$ in the following sense: for
every $\phi\in \overline{{\M}_{|\!|\!|\cdot|\!|\!|}}^\#$, there is a
unique $X\in \overline{{\M}_{|\!|\!|\cdot|\!|\!|^\#}}$ such that
\[\phi(T)=\tau(TX),\quad\forall T\in \overline{{\M}_{|\!|\!|\cdot|\!|\!|}}
\] and $\|\phi\|=|\!|\!|T|\!|\!|$?\\

\noindent Question 2: Under what conditions is
$\overline{{\M}_{|\!|\!|\cdot|\!|\!|}}$  a
reflexive Banach space?\\

Let $|\!|\!|\cdot|\!|\!|_1$ be the symmetric gauge norm on
$(L^\infty[0,1],\int_0^1dx)$ corresponding to $|\!|\!|\cdot|\!|\!|$
on $\M$ as in Corollary~{\bf 2}. Then the same questions can be
asked about $\overline{{L^\infty[0,1]}_{|\!|\!|\cdot|\!|\!|_1}}$,
the completion  of $L^\infty[0,1]$ with respect to
$|\!|\!|\cdot|\!|\!|_1$.\\

 As further consequences of
Theorem~{\bf A}, we prove the following theorems that answer the
question 1 and question 2, respectively. \\

\noindent{\bf Theorem H.}\,\, \emph{Let $\M$ be a type ${\rm II}\sb
1$ factor, $|\!|\!|\cdot|\!|\!|$ be
 a unitarily invariant norm on $\M$ and let
$|\!|\!|\cdot|\!|\!|^\#$ be the dual unitarily invariant norm on
$\M$. Let $|\!|\!|\cdot|\!|\!|_1$ be the symmetric gauge norm on
$(L^\infty[0,1],\int_0^1dx)$ corresponding to $|\!|\!|\cdot|\!|\!|$
on $\M$ as in Corollary~{\bf 2}. Then the following conditions are
equivalent:
\begin{enumerate}
 \item $\overline{{\M}_{|\!|\!|\cdot|\!|\!|^\#}}$ is the
dual space of $\overline{{\M}_{|\!|\!|\cdot|\!|\!|}}$ in the sense
of question 1;
\item $\overline{{L^\infty[0,1]}_{|\!|\!|\cdot|\!|\!|_1^\#}}$ is the
dual space of $\overline{{L^\infty[0,1]}_{|\!|\!|\cdot|\!|\!|_1}}$
in the sense of question 1;
\item $|\!|\!|\cdot|\!|\!|$ is a continuous norm on $\M$;
\item $|\!|\!|\cdot|\!|\!|_1$ is a continuous norm on $L^\infty[0,1]$.
\end{enumerate}}
\bigskip

\noindent {\bf Theorem I.}\,\,\emph{Let $\M$ be a type ${\rm
II}\sb 1$ factor, $|\!|\!|\cdot|\!|\!|$ be
 a unitarily invariant norm on $\M$ and let $|\!|\!|\cdot|\!|\!|^\#$ be the dual unitarily invariant norm on
$\M$. Let $|\!|\!|\cdot|\!|\!|_1$ be the symmetric gauge norm on
$(L^\infty[0,1],\int_0^1dx)$ corresponding to $|\!|\!|\cdot|\!|\!|$
on $\M$ as in Corollary~{\bf 2}. Then the following conditions are
equivalent:
\begin{enumerate}
\item $\overline{{\M}_{|\!|\!|\cdot|\!|\!|}}$ is a reflexive space;
\item $\overline{{L^\infty[0,1]}_{|\!|\!|\cdot|\!|\!|_1}}$ is a
reflexive space;
\item both $|\!|\!|\cdot|\!|\!|$ and $|\!|\!|\cdot|\!|\!|^\#$ are
continuous norms on $\M$;
\item both $|\!|\!|\cdot|\!|\!|_1$ and $|\!|\!|\cdot|\!|\!|_1^\#$
are continuous norms on $L^\infty[0,1]$.
\end{enumerate}}
\bigskip

A key step to proving Theorem~{\bf H} is based on the following
fact: if $|\!|\!|\cdot|\!|\!|$ is a continuous unitarily invariant
norm on $\M$ and $\phi\in \overline{\M_{|\!|\!|\cdot|\!|\!|}}^\#$,
then the restriction of $\phi$ to $\M$ is an ultraweakly continuous
linear functional, i.e., $\phi$ is in the predual space of $\M$.  A
significant advantage of our approach  is that we develop a
 relatively complete theory of unitarily invariant norms on type ${\rm II}\sb 1$ factors before
 handling unbounded operators. Indeed,
 unbounded operators are slightly involved
 only in the last two sections (section 12 and section 13). Compared with the classical methods (e.g.,~\cite{Se}), which
 have to do a lot of subtle analysis on unbounded operators, our
 methods are much simpler.\\

Let $\M$ be a finite factor. Recall that a  norm
$|\!|\!|\cdot|\!|\!|$ on $\M$ is called a normalized norm if
$|\!|\!|1|\!|\!|=1$. Let $\mathfrak{N}(\M)$ be the set of normalized
unitarily invariant norms on $\M$. Then $\NN(\M)$ is a convex
compact set in the pointwise weak topology. Let $\NN_e(\M)$ be the
set of extreme points of $\NN(\M)$. By the Krein-Milman theorem,
$\NN(\M)$ is the closure of the convex hull of $\NN_e(\M)$ in the
pointwise weak topology. So it is an interesting question of
characterizing the set
$\NN_e(\M)$. In section~10, we prove the following theorems.\\

\noindent{\bf Theorem J.}\,\,\emph{$\NN_e(M_2(\cc))=\{\max\{t\|T\|,
\|T\|_1\}: 1/2\leq t\leq 1\}$, where $\|T\|_1=\tau_2(|T|)$.}
\bigskip

\noindent {\bf Theorem K.}\,\,\emph{If $\M$ is a type ${\rm II}\sb
1$ factor and $t$ is a rational number such that $0\leq t\leq 1$,
then the Ky Fan $t$-th norm is an extreme point of $\NN(\M)$.}
\bigskip

 This paper is almost self-contained and we do not assume any
 backgrounds on non-commutative $L^p$-theory.

\section{Preliminaries}
\subsection{Nonincreasing rearrangements of functions}
Throughout this paper, we denote by $m$ the Lebesgue measure on
$[0,1]$. In the following, a measurable function and a measurable
set mean a Lesbesgue measurable function and a Lebesgue measurable
set. For two measurable sets $A$ and $B$, $A=B$ means $m((A\setminus
B)\cup (B\setminus A))=0$.

Let  $f(x)$ be a real  measurable function on $[0,1]$. The
\emph{nonincreasing rearrangement function}, $f^*(x)$, of $f(x)$ is
defined by
\begin{equation}\label{E:nonincreasing}
f^*(x)=\left\{
         \begin{array}{ll}
           \sup\{y:\, m(\{f>y\})>x\}, & \hbox{$0\leq x<1$;} \\
           \essinf\,\, f, & \hbox{$x=1$.}
         \end{array}
       \right.
\end{equation}
We summarize some useful properties of $f^*(x)$ in the following
proposition.

\begin{Proposition}\label{P:non-increasing} Let $f(x), g(x)$ be
real measurable functions on $[0,1]$. Then we have the following:
 \begin{enumerate}
\item  $f^*(x)$ is a nonincreasing, right-continuous function on
$[0,1]$ such that $f^*(0)=\esssup f$;
\item if $f(x)$ and $g(x)$ are bounded functions and
$\int_0^1f^n(x)dx=\int_0^1g^n(x)dx$ for all $n=0,1,2,\cdots,$ then
$f^*(x)=g^*(x)$;
\item $f(x)$ and $f^*(x)$ are equi-measurable and $\int_0^1 f(x)dx=\int_0^1 f^*(x)dx$ when either  integral is
well-defined.
\end{enumerate}
\end{Proposition}

\subsection{Invertible measure-preserving transformations on $[0,1]$}

Let $\mathfrak{G}=\{\phi:\, \text{$\phi(x)$ is an invertible
measure-preserving transformation on $[0,1]$}\}$. It is well known
that $\mathfrak{G}$ acts on $[0,1]$ ergodically (see~\cite{Fr} page
3-4, for instance), i.e., for a measurable subset $A$ of $[0,1]$,
$\phi(A)=A$ for all $\phi\in \mathfrak{G}$ implies that $m(A)=0$ or
$m(A)=1$.

\begin{Lemma} \label{L:two measurable sets} Let $A, B$ be two measurable subsets of $[0,1]$ such
that $m(A)=m(B)$. Then there is a $\phi\in \mathfrak{G}$ such that
$\phi(A)=B$.
\end{Lemma}
\begin{proof} We can assume that $m(A)=m(B)>0$. Since $\mathfrak{G}$
acts ergodically on $[0,1]$. There is a $\phi\in \mathfrak{G}$ such
that $m(\phi(A)\cap B)>0$. Let $B_1=\phi(A)\cap B$ and
$A_1=\phi^{-1}(B_1)$. Then $m(A_1)=m(B_1)$ and $\phi(A_1)=B_1$.  By
Zorn's lemma and maximality arguments, we prove the lemma.
\end{proof}

\begin{Corollary}\label{C:two sequences of measurable sets} Let $A_1, \cdots, A_n$  and $B_1, \cdots, B_n$ be
disjoint measurable subsets of $[0,1]$ such that
 $m(A_k)=m(B_k)$ for $1\leq k\leq n$.
Then there is a $\phi\in \mathfrak{G}$ such that $\phi(A_k)=B_k$ for
$1\leq k\leq n$.
\end{Corollary}
\begin{proof} We can assume that $A_1\cup\cdots\cup A_n=B_1\cup\cdots\cup
B_n=[0,1]$.
 By Lemma~\ref{L:two measurable sets}, there is a
$\phi_k\in \GG$ such that $\phi_k(A_k)=B_k$, $1\leq k\leq n$. Define
$\phi(x)=\phi_k(x)$ for $x\in A_k$. Then $\phi\in \GG$ and
$\phi(A_k)=B_k$ for $1\leq k\leq n$.
\end{proof}

For $f(x)\in L^\infty[0,1]$, define $\tau(f)=\int_0^1f(x)dx$.  The
following theorem is a version of the Dixmier's averaging theorem
(see~\cite{Di} or~\cite{K-R}) and it has a similar proof.
\begin{Theorem}\label{T:averaging theorem} Let $f(x)\in L^\infty[0,1]$
 be a real function. Then $\tau(f)$ is in the
$L^\infty$-norm closure of the convex hull of $\{f\circ\phi(x):\,
\phi\in \mathfrak{G}\}$.
\end{Theorem}

We end this subsection with the following proposition.
\begin{Proposition}\label{P:tracial=automorphism} If $\phi(x)$ is an invertible measure-preserving
transformation on $[0,1]$, then
\[\theta(f)=f\circ\phi
\] is a $\ast$-automorphism of $L^\infty[0,1]$ preserving $\tau$.
Conversely, if $\theta$ is a $\ast$-automorphism of $L^\infty[0,1]$
preserving $\tau$, then there is an invertible measure-preserving
transformation on $[0,1]$ such that
\[\theta(f)=f\circ\phi
\] for all $f(x)\in L^\infty[0,1]$.
\end{Proposition}
\begin{proof} The first part of the proposition is easy to see. Suppose $\theta$ is a
$\ast$-automorphism of $L^\infty[0,1]$. Let $\phi(x)=\theta(f)(x)$,
where $f(x)\equiv x$. Then it is easy to see the second part of the
proposition.
\end{proof}

\subsection{$s$-numbers of operators in type ${\rm II}\sb 1$ factors}

In~\cite{F-K}, Fack and Kosaki give a rather complete exposition of
generalized $s$-numbers  of $\tau-$measurable operators affiliated
with  semi-finite von Neumann algebras. For the sake of reader's
convenience and our purpose, we provide sufficient details on
$s$-numbers of bounded operators in finite von Neumann algebras  in
the following. We will define $s$-numbers of  bounded operators in
finite von Neumann algebras from the point of view of
non-increasing rearrangement of functions.\\

 The following lemma is
well-known. The proof is an easy exercise.

\begin{Lemma}\label{L:isomorphism}
 Let $(\A,\tau)$ be a separable \emph{(}i.e., with separable
 predual\emph{)}
 diffuse abelian von Neumann algebra  with a faithful normal trace $\tau$ on $\A$. Then there is
a $*$-isomorphism $\alpha$ from $(\A,\tau)$ onto
$(L^\infty[0,1],\int_0^1 dx)$ such that $\tau=\int_0^1 dx\circ
\alpha$.
\end{Lemma}

Let $\M$ be a type ${\rm II}\sb 1$ factor and $\tau$ be the unique
trace on $\M$. For $T\in \M$, there is a separable diffuse abelian
von Neumann subalgebra $\A$ of $\M$ containing $|T|$. By
Lemma~\ref{L:isomorphism}, there is a $\ast$-isomorphism $\alpha$
from $(\A,\tau)$ onto $(L^\infty([0,1], \int_0^1 dx)$ such that
$\tau=\int_0^1dx\circ\alpha$. Let $f(x)=\alpha(|T|)$ and $f^*(x)$
be the non-increasing rearrangement of $f(x)$
(see~$(\ref{E:nonincreasing})$). Then \emph{the $s$-numbers of
$T$}, $\mu_s(T)$, are defined as
\[\mu_s(T)=f^*(s),\,\, 0\leq s\leq 1.
\]

\begin{Lemma}\label{L:s-number}$\mu_s(T)$ does not depend on $\A$ and $\alpha$.
\end{Lemma}
\begin{proof} Let $\A_1$ be another separable diffuse abelian von
Neumann subalgebra of $\M$ containing $|T|$ and  let $\beta$ be a
$\ast$-isomorphism from $(\A_1,\tau)$ onto
$(L^\infty[0,1],\int_0^1dx)$ such that
$\tau=\int_0^1dx\circ\beta$. Let $g(x)=\beta(|T|)$. For every
number $n=0,1,2,\cdots$,
$\int_0^1f^n(x)dx=\tau(|T|^n)=\int_0^1g^n(x)dx$. Since both $f(x)$
and $g(x)$ are bounded positive functions, by 2 of
Proposition~\ref{P:non-increasing}, $f^*(x)=g^*(x)$ for all $x\in
[0,1]$.
\end{proof}
\begin{Corollary}\label{C:s-number} For $T\in \M$ and $p\geq 0$, $\tau(|T|^p)=\int_0^1 \mu_s(T)^pds$.
\end{Corollary}

The following lemma says that the above definition of $s$-numbers
coincides with the definition of $s$-numbers given by Fack and
Kosaki. Recall that $\P(\M)$ is the set of projections in $\M$.

\begin{Lemma}\label{L:s-number2} For $0\leq s\leq 1$,
\[\mu_s(T)=\inf\{\|TE\|: \, E\in \P(\M),\,  \tau(E^\perp)=s\}.
\]
\end{Lemma}
\begin{proof} By the polar decomposition and the definition of
$\mu_s(T)$, we may assume that $T$ is positive. Let $\A$ be a
separable diffuse abelian von Neumann subalgebra of $\M$
containing $T$ and let $\alpha$ be a $\ast$-isomorphism from
$(\A,\tau)$ onto $(L^\infty[0,1], \int_0^1\,dx)$ such that
$\tau=\int_0^1dx\circ \alpha$. Let $f(x)=\alpha(T)$ and $f^*(x)$
be the non-increasing rearrangement of $f(x)$. Then
$\mu_s(T)=f^*(s)$. By the definition of $f^*$,
\[m(\{f^*>\mu_s(T)\})=
\lim_{n\rightarrow\infty}m\left(\left\{f^*>\mu_s(T)+\frac{1}{n}\right\}\right)\leq
s\] and
\[m(\{f^*\geq \mu_s(T)\})\geq \lim_{n\rightarrow\infty} m\left(\left\{f^*>
\mu_s(T)-\frac{1}{n}\right\}\right)\geq s.\] Since $f^*$ and $f$ are
equi-measurable, $m(\{f>\mu_s(T)\})\leq s$ and $m(\{f\geq
\mu_s(T)\})\geq s$. Therefore, there is a measurable subset $A$ of
$[0,1]$, $\{f>\mu_s(T)\}\subset [0,1]\setminus A\subset \{f\geq
\mu_s(T)\}$, such that $m([0,1]\setminus A)=s$ and
$\|f(x)\chi_A(x)\|_\infty=\mu_s(T)$ and
$\|f(x)\chi_B(x)\|_\infty\geq \mu_s(T)$ for all $B\subset
[0,1]\setminus A$ such that $m(B)>0$. Let $F=\alpha^{-1}(\chi_A)$.
Then $\tau(F^\perp)=s$,
$\|TF\|=\|\alpha^{-1}(f\chi_A)\|_\infty=\mu_s(T)$ and $\|TF'\|\geq
\mu_s(T)$ for all nonzero subprojections $F'$ of $F^\perp$. This
proves that $\mu_s(T)\geq \inf\{\|TE\|: \, E\in \P(\M),\,
\tau(E^\perp)=s\}$. Similarly, for every $\epsilon>0$, there is a
projection $F_\epsilon\in \M$ such that
$\tau(F_\epsilon^\perp)=s+\epsilon$,
$\|TF_\epsilon\|=\mu_{s+\epsilon}(T)$ and $\|TF'\|\geq
\mu_{s+\epsilon}(T)$ for all nonzero subprojections $F'$ of
$F_\epsilon^\perp$.  Suppose $E\in \M$ is a projection such that
$\tau(E^\perp)=s$. Then $\tau(E\wedge
F_\epsilon^\perp)=\tau(E)+\tau(F_\epsilon^\perp)-\tau(E\vee
F_\epsilon^\perp)=1+\epsilon-\tau(E\vee F^\perp)\geq \epsilon>0$.
Hence, $\|TE\|\geq \|T(E\wedge F_\epsilon^\perp)\|\geq
\mu_{s+\epsilon}(T)$. This proves that $\inf\{\|TE\|: \, E\in
\P(\M),\,  \tau(E^\perp)=s\}\geq \mu_{s+\epsilon}(T)$. Since
$\mu_s(T)$ is right-continuous, $\mu_s(T)\leq \inf\{\|TE\|: \, E\in
\P(\M),\,  \tau(E^\perp)=s\}.$
\end{proof}

\begin{Corollary}\label{C:property of s-number} Let $S,T\in \M$. Then $\mu_s(ST)\leq
\|S\|\mu_s(T)$ for $s\in [0,1]$.
\end{Corollary}

We refer to~\cite{Fa,F-K} for other interesting properties of
$s$-numbers of operators in type ${\rm II}\sb 1$ factors.
\subsection{$s$-numbers of operators in  finite von Neumann
algebras}

Throughout this paper, a finite von Neumann algebra $(\M,\tau)$
means a finite von Neumann algebra $\M$ with a faithful normal
tracial state $\tau$. An \emph{embedding} of a finite von Neumann
algebra $(\M,\tau)$ into another finite von Neumann algebra
$(\M_1,\tau_1)$ means a $\ast$-isomorphism $\alpha$ from $\M$ to
$\M_1$ such that $\tau=\tau_1\cdot\alpha$. Let
$(\mathcal{L}(\mathcal{F}_2),\tau')$ be the free group factor with
the faithful normal trace $\tau'$. Then the reduced free product von
Neumann algebra $\M_1=(\M,\tau)*(\mathcal{L}(\mathcal{F}_2),\tau')$
is a type ${\rm II}\sb 1$ factor with a (unique) faithful normal
trace $\tau_1$ such that the restriction of $\tau_1$  to $\M$ is
$\tau$. So every finite von Neumann algebra can be embedded into a
type ${\rm II}\sb 1$ factor.

\begin{Definition}\label{D:s-number}\emph{ Let $(\M,\tau)$ be a finite von
Neumann algebra and $T\in \M$. If $\alpha$ is an embedding of
$(\M,\tau)$ into a type ${\rm II}\sb 1$ factor $(\M_1,\tau_1)$,
then \emph{the $s$-numbers of $T$} are defined as
\[\mu_s(T)=\mu_s(\alpha(T)).
\]}
\end{Definition}

Similar to the proof of Lemma~\ref{L:s-number}, we can see that
$\mu_s(T)$ is well defined, i.e., does not depend on the choice of
$\alpha$ and $\M_1$.\\

Let $T\in (M_n(\cc),\tau_n)$, where $\tau_n$ is the normalized trace
on $M_n(\cc)$. Then $|T|$ is unitarily equivalent to a diagonal
matrix  with diagonal elements
 $s_1(T)\geq \cdots\geq
s_n(T)\geq 0$.  In the classical matrices theory~\cite{Bh, G-K},
$s_1(T),\cdots,s_n(T)$ are also called $s$-numbers of $T$. It is
easy to see that the relation between $\mu_s(T)$ and
$s_1(T),\cdots,s_n(T)$  is the following
\begin{equation}\label{E:s-number}
\mu_s(T)=s_1(T)\chi_{[0,1/n)}(s)+s_2(T)\chi_{[1/n,2/n)}(s)+\cdots+s_n(T)\chi_{[n-1/n,1]}(s).
\end{equation}
Since no confusions will arise, we will use both $s$-numbers for
matrices in $M_n(\cc)$. We refer to~\cite{Bh, G-K} for other
interesting properties of $s$-numbers of matrices.

We end this section by the following definition.
\begin{Definition}\label{D:equimeasurable}\emph{
 Positive operators $S$ and $T$ in a finite von Neumann
algebra $(\M,\tau)$ are \emph{equi-measurable} if
$\mu_s(S)=\mu_s(T)$ for $0\leq s\leq 1$.}
\end{Definition}

By 2 of Proposition~\ref{P:non-increasing} and
Corollary~\ref{C:s-number}, positive operators $S$ and $T$ in a
finite von Neumann algebra $(\M,\tau)$ are equi-measurable if and
only if $\tau(S^n)=\tau(T^n)$ for all $n=0,1,2,\cdots$.
\section{Tracial gauge semi-norms on finite von Neumann algebras
satisfying the weak Dixmier property}
\subsection{Gauge semi-norms}
\begin{Definition}\label{D:gauge seminorm}\emph{ Let $(\M,\tau)$ be a finite von Neumann
algebra. A semi-norm $|\!|\!|\cdot|\!|\!|$  on $\M$ is called
\emph{gauge invariant} if for every $T\in\M$,
\[|\!|\!|T|\!|\!|=|\!|\!|\,\, |T|\,\,|\!|\!|.\]}
\end{Definition}

\begin{Lemma}\label{L:submultiplicative}Let $(\M,\tau)$ be a finite von Neumann algebra
and let $|\!|\!|\cdot|\!|\!|$ be a semi-norm on $\M$. Then the
following conditions are equivalent:
\begin{enumerate}
\item $|\!|\!|\cdot|\!|\!|$ is gauge invariant;
\item $|\!|\!|\cdot|\!|\!|$ is left unitarily invariant, i.e., for
every unitary operator $U\in\M$ and operator $T\in \M$,
$|\!|\!|UT|\!|\!|=|\!|\!|T|\!|\!|$;
\item for operators $A,T\in \M$, $|\!|\!|AT|\!|\!|\leq \|A\|\cdot|\!|\!|T|\!|\!|$.
\end{enumerate}
\end{Lemma}
\begin{proof} $``3\Rightarrow 2"$ and $``2\Rightarrow 1"$ are easy to see.
 We only prove $``1\Rightarrow 3"$. We need to prove that if
$\|A\|<1$, then $|\!|\!|AT|\!|\!|\leq |\!|\!|T|\!|\!|$. Since
$\|A\|<1$, there are unitary operators $U_1,\cdots,U_k$ such that
$A=\frac{U_1+\cdots+U_k}{k}$ (see~\cite{K-P, R-D}). Since
$|U_1T|=\cdots=|U_kT|=|T|$, $|\!|\!|AT|\!|\!|=
|\!|\!|\frac{U_1T+\cdots+U_kT}{k}|\!|\!|\leq
\frac{|\!|\!|U_1T|\!|\!|+\cdots+|\!|\!|U_kT|\!|\!|}{k}\leq
|\!|\!|T|\!|\!|$.
\end{proof}

\begin{Definition} \label{D:normalized norm} \emph{A \emph{normalized}
 semi-norm  on a finite
von Neumann algebra $(\M,\tau)$ is a semi-norm $|\!|\!|\cdot|\!|\!|$
such that $|\!|\!|1|\!|\!|=1$.}
\end{Definition}

By Lemma~\ref{L:submultiplicative}, we have the following corollary.
\begin{Corollary}\label{C:compare gauge norm with operator norm}
Let $(\M,\tau)$ be a finite von Neumann algebra and let
$|\!|\!|\cdot|\!|\!|$ be a normalized gauge semi-norm on $\M$.
Then for every $T\in\M$,
\[|\!|\!|T|\!|\!|\leq \|T\|.\]
\end{Corollary}

A \emph{simple} operator in a finite von Neumann algebra $(\M,\tau)$
is an operator $T=a_1E_1+\cdots+a_nE_n$, where $E_1,\cdots,E_n$ are
projections in $\M$ such that $E_1+\cdots+E_n=1$.
\begin{Corollary}\label{C:simple operators are dense} Let $(\M,\tau)$ be a finite von Neumann algebra,
$|\!|\!|\cdot|\!|\!|_1$ and $|\!|\!|\cdot|\!|\!|_2$ be two gauge
invariant semi-norms on $\M$. Then
$|\!|\!|\cdot|\!|\!|_1=|\!|\!|\cdot|\!|\!|_2$ on $\M$ if
$|\!|\!|T|\!|\!|_1=|\!|\!|T|\!|\!|_2$ for all positive simple
operators $T\in \M$.
\end{Corollary}
\begin{proof} Without loss of generality, assume $|\!|\!|1|\!|\!|_1=|\!|\!|1|\!|\!|_2=1$.
 Let $T\in \M$ be a positive operator.  By
the spectral decomposition theorem, there is a sequence of positive
simple operators $T_n\in\M$ such that
$\lim_{n\rightarrow\infty}\|T-T_n\|=0$. By Corollary~\ref{C:compare
gauge norm with operator norm},
$\lim_{n\rightarrow\infty}|\!|\!|T-T_n|\!|\!|_1=\lim_{n\rightarrow\infty}|\!|\!|T-T_n|\!|\!|_2=0$.
By the assumption of the corollary,
$|\!|\!|T_n|\!|\!|_1=|\!|\!|T_n|\!|\!|_2$. Hence,
$|\!|\!|T|\!|\!|_1=|\!|\!|T|\!|\!|_2$. Since both
$|\!|\!|\cdot|\!|\!|_1$ and $|\!|\!|\cdot|\!|\!|_2$ are gauge
invariant, $|\!|\!|\cdot|\!|\!|_1=|\!|\!|\cdot|\!|\!|_2$.
\end{proof}

\subsection{Tracial gauge semi-norms}

\begin{Definition}\label{D: tracial gauge}\emph{ Let $(\M,\tau)$ be a finite von Neumann
algebra. A semi-norm $|\!|\!|\cdot|\!|\!|$ on $\M$  is called
\emph{tracial} if $|\!|\!|S|\!|\!|=|\!|\!|T|\!|\!|$ for every two
equi-measurable positive operators $S,T$ in $\M$. A semi-norm
$|\!|\!|\cdot|\!|\!|$ on $\M$ is called a \emph{tracial gauge
semi-norm} if it is both tracial and gauge invariant.}
\end{Definition}

Since for a positive operator $T$ in a finite von Neumann algebra
$(\M,\tau)$,
$\|T\|=\lim_{n\rightarrow\infty}\left(\tau(T^n)\right)^{\frac{1}{n}}$,
the operator norm $\|\cdot\|$ is a tracial gauge norm on
$(\M,\tau)$. Another less obvious example of a tracial gauge norm on
$(\M,\tau)$ is the non-commutative $L^1$-norm:
$||T||_1=\tau(|T|)=\int_0^1\mu_s(T)ds$. The less obvious part is to
show that $||\cdot||_1$ satisfies the triangle inequality. The
following lemma overcomes this difficulty.
\begin{Lemma} \label{L:1-norm} $||A||_1=\sup\{|\tau(UA)|:\, U\in
\mathcal{U}(\M)\}$, where $\U(\M)$ is the set of unitary operators
in $\M$.
\end{Lemma}
\begin{proof} By the polar decomposition theorem, there is a unitary
operator $V\in \M$ such that $A=V|A|$. By the Schwartz inequality,
$|\tau(UA)|=|\tau(UV|A|)|=|\tau(UV|A|^{1/2}|A|^{1/2})|\leq
\tau(|A|)^{1/2}\cdot \tau(|A|)^{1/2}=\tau(|A|)$. Hence
$||A||_1\geq \sup\{|\tau(UA)|:\, U\in \mathcal{U}(\M)\}$. Let
$U=V^*$, we obtain $||A||_1\leq \sup\{|\tau(UA)|:\, U\in
\mathcal{U}(\M)\}$.
\end{proof}

\begin{Corollary}
$\|A+B\|_1\leq \|A\|_1 + \|B\|_1$.
\end{Corollary}

\begin{Lemma}\label{L:sufficient condition for
tracial} Let $(\M,\tau)$ be a finite von Neumann algebra and let
$|\!|\!|\cdot|\!|\!|$ be a gauge invariant semi-norm on $\M$. Then
$|\!|\!|\cdot|\!|\!|$ is tracial if
$|\!|\!|S|\!|\!|=|\!|\!|T|\!|\!|$ for every two equi-measurable
positive simple operators $S,T$ in $\M$.
\end{Lemma}
\begin{proof} We can assume that $|\!|\!|1|\!|\!|=1$. Let $A,B$ be two equi-measurable positive operators in
$\M$. By the spectral decomposition theorem, there are two sequence
of positive simple operators $A_n, B_n$ in $\M$ such that $A_n$ and
$B_n$ are equi-measurable and
$\lim_{n\rightarrow\infty}\|A-A_n\|=\lim_{n\rightarrow\infty}\|B-B_n\|=0$.
By Corollary~\ref{C:compare gauge norm with operator norm},
$\lim_{n\rightarrow\infty}|\!|\!|A-A_n|\!|\!|=\lim_{n\rightarrow\infty}|\!|\!|B-B_n|\!|\!|=0$.
By the assumption of the lemma,
$|\!|\!|A_n|\!|\!|=|\!|\!|B_n|\!|\!|$. Hence,
$|\!|\!|A|\!|\!|=|\!|\!|B|\!|\!|$.
\end{proof}
\subsection{Symmetric gauge semi-norms}
\begin{Definition}\label{D: Symmetric gauge}
\emph{ Let $(\M,\tau)$ be a finite von Neumann algebra  and
$\Aut(\M, \tau)$ be  the set of $\ast$-automorphisms  on $\M$
preserving $\tau$.  A semi-norm
 $|\!|\!|\cdot|\!|\!|$ on
$\M$ is called \emph{symmetric} if
\[|\!|\!|\theta(T)|\!|\!|=|\!|\!|T|\!|\!|,\,\,\,\, \forall T\in \M,\, \theta\in
\Aut(\M,\tau).
\] A semi-norm $|\!|\!|\cdot|\!|\!|$ on $\M$ is called a \emph{symmetric gauge} semi-norm if it
is both symmetric and gauge invariant.}
\end{Definition}

\begin{Example}\emph{ Let $\M=\cc^n$ and
$\tau(T)=\frac{x_1+\cdots+x_n}{n}$, where $T=(x_1,\cdots,x_n)\in
\cc^n$. Then $\Aut(\M,\tau)$ is the set of permutations on
$\{1,\cdots,n\}$. So a semi-norm $|\!|\!|\cdot|\!|\!|$ on $\M$ is a
symmetric gauge semi-norm if and only if for every
$(x_1,\cdots,x_n)\in \cc^n$ and a permutation $\pi$ on
$\{1,\cdots,n\}$,
\[|\!|\!|(x_1,\cdots,x_n)|\!|\!|=|\!|\!|(|x_1|,\cdots,|x_n|)|\!|\!|,
\]
and
\[|\!|\!|(x_1,\cdots,x_n)|\!|\!|=|\!|\!|(x_{\pi(1)},\cdots,x_{\pi(n)})|\!|\!|.
\]}
\end{Example}

\begin{Lemma}\label{L:tracial implies symmetric}Let $(\M,\tau)$ be a finite von Neumann algebra and let $|\!|\!|\cdot|\!|\!|$ be a semi-norm on $\M$. If
$|\!|\!|\cdot|\!|\!|$ is tracial gauge invariant, then
$|\!|\!|\cdot|\!|\!|$ is symmetric gauge invariant.
\end{Lemma}
\begin{proof} Let $\theta\in \Aut(\M,\tau)$ and $T\in \M$. We need to prove
that $|\!|\!|\theta(T)|\!|\!|=|\!|\!|T|\!|\!|$.  Since
$|\theta(T)|=\theta(|T|)$ and $|\!|\!|\cdot|\!|\!|$ is gauge
invariant, we can assume that $T$ is positive. Since $\theta\in$
$\Aut(\M,\tau)$, $T$ and $\theta(T)$ are equi-measurable. Hence,
$|\!|\!|T|\!|\!|=|\!|\!|\theta(T)|\!|\!|$.
\end{proof}
\begin{Example}\emph{ Let $\M=\cc\oplus M_2(\cc)$ and $\tau(a\oplus
B)=\frac{a}{2}+\frac{\tau_2(B)}{2}$, where $\tau_2$ is the
normalized trace on $M_2(\cc)$. Define $|\!|\!|a\oplus B|\!|\!|$ $=$
$|a|/2+\tau_2(|B|)$. Then $|\!|\!|\cdot|\!|\!|$ is a symmetric gauge
norm but not a tracial gauge norm. Note that $1\oplus 0$ and
$0\oplus 1$ are equi-measurable, but $1/2=|\!|\!|1\oplus
0|\!|\!|\neq |\!|\!|0\oplus 1|\!|\!|=1$.}
\end{Example}

$\Aut(\M,\tau)$ acts on $\M$ \emph{ergodically} if $\theta(T)=T$ for
all $\theta\in \Aut(\M,\tau)$ implies $T=\lambda 1$.
\begin{Lemma}\label{L:symmetric implies tracial} Let $(\M,\tau)$ be a finite von Neumann algebra and let $|\!|\!|\cdot|\!|\!|$ be a semi-norm on $\M$. If $\Aut(\M,\tau)$
acts on $\M$ ergodically, then the following are equivalent:
\begin{enumerate}
\item $|\!|\!|\cdot|\!|\!|$ is a tracial gauge semi-norm;
\item $|\!|\!|\cdot|\!|\!|$ is a symmetric gauge semi-norm.
\end{enumerate}
\end{Lemma}
\begin{proof} $``1\Rightarrow 2"$ by Lemma~\ref{L:tracial implies symmetric}. We need to prove
$``2\Rightarrow 1"$. By Corollary~\ref{C:simple operators are
dense}, we need to prove $|\!|\!|S|\!|\!|=|\!|\!|T|\!|\!|$ for two
equi-measurable simple operators $S,T$ in $\M$. Similar to the proof
of Corollary~\ref{C:two sequences of measurable sets},
 there is a $\theta\in \Aut(\M,\tau)$ such that
$S=\theta(T)$. Hence $|\!|\!|S|\!|\!|=|\!|\!|T|\!|\!|$.
\end{proof}

\begin{Corollary}\label{C:symmetric norm on L infty} A semi-norm
 on $(L^\infty[0,1],\int_0^1dx)$ or $(\cc^n,\tau)$
  is a
tracial gauge norm if and only if it is a symmetric gauge norm,
where $\tau((x_1,\cdots,x_n))=\frac{x_1+\cdots+x_n}{n}$.
\end{Corollary}

\subsection{Unitarily invariant semi-norms}
\begin{Definition}\label{D:unitarily invariant norm}
\emph{ Let $(\M,\tau)$ be a  von Neumann algebra.
 A semi-norm $|\!|\!|\cdot|\!|\!|$ on $\M$ is \emph{unitarily invariant} if $|\!|\!|UTV|\!|\!|=|\!|\!|T|\!|\!|$ for
all $T\in \M$ and unitary operators $U,V\in \M$.}
\end{Definition}

\begin{Proposition}\label{P:unitarily invariant norms} Let $|\!|\!|\cdot|\!|\!|$ be a semi-norm on $\M$. Then the
following statements are equivalent:
\begin{enumerate}
\item $|\!|\!|\cdot|\!|\!|$ is unitarily invariant;
\item $|\!|\!|\cdot|\!|\!|$ is gauge invariant and unitarily conjugate
invariant, i.e.,  $|\!|\!|UTU^*|\!|\!|=|\!|\!|T|\!|\!|$ for all
$T\in\M$ and unitary operators $U\in \M$;
\item $|\!|\!|\cdot|\!|\!|$ is
left-unitarily invariant and $|\!|\!|T|\!|\!|=|\!|\!|T^*|\!|\!|$ for
every $T\in \M$;
\item for all operators
$T,A,B\in \M$, $|\!|\!|ATB|\!|\!|\leq \|A\|\cdot
|\!|\!|T|\!|\!|\cdot \|B\|$.
\end{enumerate}
\end{Proposition}
\begin{proof} $``1\Rightarrow 4"$ is similar to the proof of Lemma~\ref{L:submultiplicative}.
$``4\Rightarrow 3"$, $``3\Rightarrow 2"$, and $``2\Rightarrow 1"$
are routine.
\end{proof}

\begin{Corollary}\label{C:S<T}Let $(\M,\tau)$ be a finite von Neumann algebra and let $|\!|\!|\cdot|\!|\!|$ be
a  unitarily invariant semi-norm on $\M$.  If $0\leq S\leq T$,
then $|\!|\!|S|\!|\!|\leq |\!|\!|T|\!|\!|$.
\end{Corollary}
\begin{proof} Since $0\leq S\leq T$, there is an operator $A\in \M$
such that $S=ATA^*$ and $\|A\|\leq 1$. By
Proposition~\ref{P:unitarily invariant norms},
\[|\!|\!|S|\!|\!|=|\!|\!|ATA^*|\!|\!|\leq
\|A\|\cdot|\!|\!|T|\!|\!|\cdot\|A^*\|\leq |\!|\!|T|\!|\!|.\]
\end{proof}

For a unitary operator $U\in \M$, let $\theta(T)=UTU^*$. Then
$\theta\in \Aut(\M,\tau)$. By Proposition~\ref{P:unitarily invariant
norms}, we have the following
\begin{Corollary}\label{C:symmetric implies unitarily invariant} Let $(\M,\tau)$ be a finite von Neumann algebra and let $|\!|\!|\cdot|\!|\!|$ be a symmetric, gauge invariant semi-norm on
$\M$. Then $|\!|\!|\cdot|\!|\!|$ is a unitarily invariant semi-norm
on $\M$.
\end{Corollary}

\begin{Example}\emph{ Let $\M=\cc^n$, $n\geq 2$ and
$\tau((x_1,\cdots,x_n))=\frac{x_1+\cdots+x_n}{n}$. Define
$|\!|\!|(x_1,\cdots,x_n)|\!|\!|$ $=$ $|x_1|$. Then
$|\!|\!|\cdot|\!|\!|$ is a unitarily invariant semi-norm but not a
symmetric gauge semi-norm on $\M$.}
\end{Example}

\begin{Lemma}\label{L:unitarily invariant implies tracial} Let $(\M,\tau)$ be a finite factor and let $|\!|\!|\cdot|\!|\!|$ be a semi-norm on $\M$. Then the following
conditions are equivalent:
\begin{enumerate}
\item $|\!|\!|\cdot|\!|\!|$ is a tracial gauge semi-norm;
\item $|\!|\!|\cdot|\!|\!|$ is a symmetric gauge semi-norm;
\item $|\!|\!|\cdot|\!|\!|$ is a unitarily invariant semi-norm.
\end{enumerate}
\end{Lemma}
\begin{proof} $``1\Rightarrow 2"$ by Lemma~\ref{L:tracial implies symmetric} and $``2\Rightarrow
3"$ by Corollary~\ref{C:symmetric implies unitarily invariant}. We
need to prove $``3\Rightarrow 1"$. By Corollary~\ref{C:simple
operators are dense}, we need to prove
$|\!|\!|S|\!|\!|=|\!|\!|T|\!|\!|$ for two equi-measurable positive
simple operators $S,T\in \M$. Suppose $S=a_1E_1+\cdots+a_nE_n$ and
$T=a_1F_1+\cdots+a_nF_n$, where $E_1+\cdots+E_n=1$ and
$F_1+\cdots+F_n=1$ and $\tau(E_k)=\tau(F_k)$ for $1\leq k\leq n$.
Since $\M$ is a factor, there is a unitary operator $U\in\M$ such
that $E_k=UF_kU^*$ for $1\leq k\leq n$. Therefore, $S=UTU^*$ and
$|\!|\!|S|\!|\!|=|\!|\!|T|\!|\!|$.
\end{proof}

\subsection{Weak Dixmier property}

\begin{Definition}\label{D: weak Dixmier property}
\emph{A finite von Neumann algebra $(\M,\tau)$ satisfies \emph{the
weak Dixmier property} if for every positive operator $T\in \M$,
$\tau(T)$ is in the operator norm closure of the convex hull of
$\{S\in \M:\,\, \text{$S$ and $T$ are equi-measurable}\}$.}
\end{Definition}

A finite factor $(\M,\tau)$ satisfies the \emph{Dixmier property}
(see \cite{Di1,K-R}): for
 every operator $T\in \M$, $\tau(T)$ is in the operator norm closure of
the convex hull of $\{UTU^*:\,\, U\in \U(\M)\}$. Hence finite
factors satisfy the weak Dixmier property. In the following, we will
characterize finite von Neumann algebras satisfying the weak
Dixmier property.\\

There is a central projection $P$ in a finite von Neumann algebra
$(\M,\tau)$ such that $P\M$ is type ${\rm I}$ and $(1-P)\M$ is
type ${\rm II}$. A type ${\rm II}$ von Neumann algebra is diffuse,
i.e, there are no minimal projections in the von Neumann algebra.
Furthermore, there are central projections $P_1,\cdots,
P_n,\cdots$ in $\M$, such that $P_1+\cdots+P_n+\cdots=P$ and
$P_n\M= \A_n\otimes M_n(\cc)$, $\A_n$ is abelian.  We can
decompose $\A_n$ into an atomic part $\A_n^a$ and a diffuse part
$\A_n^c$, i.e., there is a projection $Q_n$ in $\A_n$,
$\A_n^a=Q_n\A_n$, such that $Q_n=E_{n1}+E_{n2}+\cdots$, where
$E_{nk}$ is a minimal projection in $\A_n^a$ and $\tau(E_{nk})>0$,
and $\A_n^c=(1-Q_n)\A_n$ is diffuse. Let $\M_a=\sum_\oplus
\A_n^a\otimes M_n(\cc)$ and $\M_c=\sum_\oplus \A_n^c\otimes
M_n(\cc)\oplus (1-P)\M$. Then $\M=\M_a\oplus \M_c$. We call $\M_a$
the \emph{atomic part} of $\M$ and $\M_c$ the \emph{diffuse part}
of $\M$. A finite von Neumann algebra $(\M,\tau)$ is \emph{atomic}
if $\M=\M_a$ and is \emph{diffuse} if $\M=\M_c$.

\begin{Lemma}\label{L:finite dimensional vNA satisfies} Let  $(\M,\tau)$ be a finite dimensional von Neumann algebra
such that for every two non-zero minimal projections $E,F\in \M$,
$\tau(E)=\tau(F)$. Then $(\M,\tau)$ satisfies the weak Dixmier
property.
\end{Lemma}
\begin{proof} Since $\M$ is finite dimensional, $\M\cong M_{k_1}(\cc)\oplus \cdots \oplus M_{k_r}(\cc)$.
Since $\tau(E)=\tau(F)$ for every two non-zero minimal projections
$E,F\in \M$, $(\M,\tau)$ can be embedded into $(M_n(\cc),\tau_n)$,
where $n=k_1+\cdots+k_r$. So we can assume that $(\M,\tau)$ is a von
Neumann subalgebra of $(M_n(\cc),\tau_n)$ such that $\M$ contains
all diagonal matrices $a_1E_1+\cdots+a_nE_n$. Now for every positive
operator $T\in \M$, there is a unitary operator $U\in \M$ such that
$UTU^*=a_1E_1+\cdots+a_nE_n$, $a_1,\cdots,a_n\geq 0$ and
$\tau(T)=\frac{a_1+\cdots+a_n}{n}$. Then
$\tau(T)=\frac{\sum_{\pi}(a_{\pi(1)}E_1+\cdots+a_{\pi(n)}E_n)}{n!}$.
\end{proof}

\begin{Lemma}\label{L:diffuse finite von Neumanna algebra satisfy} Let $(\M,\tau)$ be a diffuse finite von Neumann
algebra. Then $(\M,\tau)$ satisfies the weak Dixmier property.
\end{Lemma}
\begin{proof} Let $\A$ be a separable diffuse abelian von Neumann
subalgebra of $\M$. By Lemma~\ref{L:isomorphism}, there is a
$*-$isomorphism $\alpha$ from $(\A,\tau)$ onto
$(L^\infty[0,1],\int_0^1dx)$ such that $\int_0^1dx\cdot
\alpha=\tau$. For a positive operator $T\in \M$, there is an
operator $S\in \A$ such that $\alpha(S)=\mu_s(T)$. Hence
$\tau(T)=\tau(S)=\int_0^1\mu_s(T)ds$. By Theorem~\ref{T:averaging
theorem}, for any $\epsilon>0$, there are $S_1,\cdots, S_n$ in $\A$
such that $S,S_1,\cdots,S_n$ are equi-measurable and
$\|\tau(S)-\frac{S_1+\cdots+S_n}{n}\|<\epsilon$. Hence $(\M,\tau)$
satisfies the weak Dixmier property.
\end{proof}

\begin{Lemma}\label{L:no Dixmier property 1} Let $(\M,\tau)$ be an
atomic finite von Neumann algebra with two minimal projections $E$
and $F$ in $\M$ such that $\tau(E)\neq \tau(F)$. Then $(\M,\tau)$
does not satisfy the weak Dixmier property.
\end{Lemma}
\begin{proof} Since $(\M,\tau)$ is an atomic finite von Neumann
algebra, $\M\cong M_{k_1}(\cc)\oplus M_{k_2}(\cc)\oplus\cdots$. Let
$E_{ij}$ be minimal projections in $M_{k_i}$ such that $\sum
E_{ij}=1$. Without loss of generality, assume that
$\tau(E_{11})>\tau(E_{21})\geq \tau(E_{31})\geq \cdots$. Let
$T=\left(\begin{array}{ccc}1&&\\
&\ddots&\\
&&1\end{array}\right)_{k_1}\oplus A,$ where
\[A=\left(\begin{array}{ccc}\frac{1}{2}&&\\
&\ddots&\\
&&\left(\frac{1}{2}\right)^{k_2}\end{array}\right)\oplus \left(\begin{array}{ccc}\left(\frac{1}{2}\right)^{k_2+1}&&\\
&\ddots&\\
&&\left(\frac{1}{2}\right)^{k_2+k_3}\end{array}\right)\oplus \cdots.
\]
If $T_1\in \M$ and $T$ are equi-measurable, then $T_1=\left(\begin{array}{ccc}1&&\\
&\ddots&\\
&&1\end{array}\right)_{k_1}\oplus A_1$, where $A$ and $A_1$ are
equi-measurable. Hence, if $\tau(T)$ is in the operator norm closure
of the convex hull of $\{S\in \M:\,\, \text{$S$ and $T$ are
equi-measurable}\}$, then $\tau(T)=1$. It is a contradiction.
\end{proof}

Let $(\M,\tau)$ be a finite von Neumann algebra and $E\in \M$ be a
non-zero projection. The \emph{induced finite von Neumann algebra}
$(\M_E,\tau_E)$ is the von Neumann algebra $\M_E=E\M E$ with a
faithful normal trace $\tau_E(ETE)=\frac{\tau(ETE)}{\tau(E)}$. The
proof of the following lemma is similar to the proof of
Lemma~\ref{L:no Dixmier property 1}.

\begin{Lemma}\label{L:non-dixmier property 2} Let $(\M,\tau)$ be a finite von Neumann
algebra such that $\M_a\neq 0$ and $\M_c\neq 0$. Then $\M$ does not
satisfies the weak Dixmier property.
\end{Lemma}
\begin{proof} Let $P$ be the central projection such that $\M_a=P\M$
and $\M_c=(1-P)\M$.
 Let $\A$ be a separable diffuse abelian von Neumann
subalgebra of $(\M_c,\tau_{1-P})$. By Lemma~\ref{L:isomorphism},
there is a positive operator $A$ in $\M_c$ such that
$\mu_s(A)=\frac{1-s}{2}$ with respect to $(\M_c,\tau_{1-P})$.
Consider $T=P+A(1-P)$. Then
\[\mu_s(T)=\left\{
             \begin{array}{ll}
               1, & \hbox{$0\leq s<\tau(P)$;} \\
               \frac{1-s}{2\tau(1-P)}\leq \frac{1}{2}, & \hbox{$\tau(P)\leq s\leq 1$}
             \end{array}
           \right.
\] with respect to $(\M,\tau)$.  If $T_1\in \M$ and $T$ are
equi-measurable, then $T_1=P+A_1$ such that  $A_1$ and $A$ are
equi-measurable. Hence, if $\tau(T)$ is in the operator norm closure
of the convex hull of $\{S\in \M:\,\, \text{$S$ and $T$ are
equi-measurable}\}$, then $\tau(T)=1$. It is a contradiction.
\end{proof}

Summarizing Lemma~\ref{L:finite dimensional vNA satisfies},
\ref{L:diffuse finite von Neumanna algebra satisfy}, \ref{L:no
Dixmier property 1}, \ref{L:non-dixmier property 2}, we can
characterize finite von Neumann algebras satisfying the weak Dixmier
property as the following theorem.
\begin{Theorem}\label{T:the weak Dixmier property} Let $(\M,\tau)$ be a finite  von Neumann algebra.
Then $\M$ satisfies the weak Dixmier property if and only if $\M$
satisfies one of the following conditions:
\begin{enumerate}
\item $\M$ is finite dimensional $($hence atomic$)$ and for every two non-zero minimal projections $E,F\in
\M$, $\tau(E)=\tau(F)$, or equivalently, $(\M,\tau)$ can be
identified as a von Neumann subalgebra of $(M_n(\cc),\tau_n)$ that
contains all diagonal matrices;
\item $\M$ is diffuse.
\end{enumerate}
\end{Theorem}

\begin{Corollary}\label{C:local Dixmier  property} Let $(\M,\tau)$
be a finite von Neumann algebra satisfying the weak Dixmier property
and $E\in\M$ be a non-zero projection. Then $(\M_E,\tau_E)$ also
satisfies the weak Dixmier property.
\end{Corollary}

The following example shows that we can not replace the weak Dixmier
property by the following condition: $\tau(T)$ is in the operator
norm closure of the convex hull of $\{\theta(T):\,\, \theta\in
\Aut(\M,\tau)\}$.
\begin{Example}\label{E:the weak Dixmier property}\emph{ $(\cc\oplus M_2(\cc), \tau)$, $\tau(a\oplus
B)=\frac{1}{3}a+\frac{2}{3}\tau_2(B)$, satisfies the weak Dixmier
property.  On the other hand, let $T=1\oplus 2\in \cc\oplus
M_2(\cc)$. Then for every $\theta\in \Aut(\M,\tau)$, $\theta(T)=T$.
Hence, $\tau(T)$ is not in the operator norm closure of the convex
hull of $\{\theta(T):\,\, \theta\in \Aut(\M,\tau)\}$.}
\end{Example}

\subsection{A comparison theorem}

The following theorem is the main result of this section.
\begin{Theorem}\label{T:tracial} Let $(\M,\tau)$ be a finite von Neumann algebra satisfying
the weak Dixmier property. If  $|\!|\!|\cdot|\!|\!|$ is  a
normalized tracial gauge semi-norm on $\M$, then for all $T\in \M$,
\[\|T\|_1\leq |\!|\!|T|\!|\!|\leq \|T\|.
\] In particular, every tracial gauge
semi-norm on $\M$ is a norm.
\end{Theorem}
\begin{proof} By Corollary~\ref{C:compare gauge norm with operator norm},
$|\!|\!|T|\!|\!|\leq \|T\|$ for every $T\in \M$.  To prove
$\|T\|_1\leq |\!|\!|T|\!|\!|$,  we can assume $T\geq 0$.  Let
$\epsilon>0$. Since $(\M,\tau)$ satisfies the weak Dixmier property,
there are $S_1,\cdots,S_k$ in $\M$ such that $T,S_1,\cdots,S_k$ are
equi-measurable and $\|\tau(T)-\frac{S_1+\cdots+S_k}{k}\|<\epsilon$.
By Corollary~\ref{C:compare gauge norm with operator norm},
$|\!|\!|\tau(T)-\frac{S_1+\cdots+S_k}{k}|\!|\!|\leq
\|\tau(T)-\frac{S_1+\cdots+S_k}{k}\|<\epsilon$. Hence
$\|T\|_1=|\tau(T)|\leq
|\!|\!|\frac{S_1+\cdots+S_k}{k}|\!|\!|+\epsilon\leq
\frac{|\!|\!|S_1|\!|\!|+\cdots+|\!|\!|S_k|\!|\!|}{k}+\epsilon=|\!|\!|T|\!|\!|+\epsilon.$
\end{proof}

By Theorem~\ref{T:tracial} and Lemma~\ref{L:unitarily invariant
implies tracial}, we have the following corollary.
\begin{Corollary}\label{C:unitarily seminorm is norm} Let $(\M,\tau)$ be a finite factor and let $|\!|\!|\cdot|\!|\!|$ be a normalized unitarily
invariant norm on $\M$. Then
\[\|T\|_1\leq |\!|\!|T|\!|\!|\leq \|T\|,\,\,\,\, \forall T\in\M.
\] In particular, every unitarily invariant semi-norm on a finite
factor is a norm.
\end{Corollary}

By Theorem~\ref{T:tracial} and Lemma~\ref{L:symmetric implies
tracial}, we have the following corollary.

\begin{Corollary}\label{C:symmetric gauge semi-norm on L infty is norm}
 Let $|\!|\!|\cdot|\!|\!|$ be a normalized
symmetric gauge  semi-norm on $(L^\infty[0,1], \int_0^1dx)$ $($or
$(\cc^n,\tau)$, where
$\tau((x_1,\cdots,x_n))=\frac{x_1+\cdots+x_n}{n}$$)$.  Then
\[\|T\|_1\leq |\!|\!|T|\!|\!|\leq \|T\|,\,\,\,\, \forall T\in L^\infty[0,1]\,\,(\text{or}\,\,\cc^n).
\] In particular, every symmetric gauge semi-norm on $(L^\infty[0,1], \int_0^1dx)$ $($or
$(\cc^n,\tau)$ respectively$)$ is a norm.
\end{Corollary}

\section{Proof of Theorem B}

To prove Theorem~{\bf B}, we need the following lemmas.
\begin{Lemma}\label{L:diagonal operators and convex hull} Let $E_1,\cdots, E_n$ be projections in $\M$ such that
$E_1+\cdots+E_n=1$ and $T\in \M$. Then $S=E_1TE_1+\cdots+E_nTE_n$ is
in the convex hull of $\{UTU^*:\, U\in \U(\M)\}$.
\end{Lemma}
\begin{proof} Let $T=(T_{ij})$ be the matrix with respect to the
decomposition  $1=E_1+\cdots+E_n$. Let $U=-E_1+E_2+\cdots+E_n$. Then
simple computation shows that
\[\frac{1}{2}(UTU^*+T)=\left(\begin{array}{cccc}
T_{11}&0&\cdots&0\\
0&T_{22}&\cdots&T_{2n}\\
\vdots&\vdots&\ddots&\vdots\\
0&T_{n2}&\cdots&T_{nn}
\end{array}\right)=E_1TE_1+(1-E_1)T(1-E_1).
\] By induction, $S=E_1TE_1+\cdots+E_nTE_n$ is
in the convex hull of $\{UTU^*:\, U\in \U(\M)\}$.
\end{proof}

\begin{Corollary}\label{C:diagonal operators} Let $(\M,\tau)$ be a finite von Neumann algebra and let $|\!|\!|\cdot|\!|\!|$
be a unitarily invariant norm on $\M$. Let $E_1,\cdots, E_n$ be
projections in $\M$ such that $E_1+\cdots+E_n=1$ and $T\in \M$ and
$S=E_1TE_1+\cdots+E_nTE_n$. Then $|\!|\!|S|\!|\!|\leq
|\!|\!|T|\!|\!|$.
\end{Corollary}

Recall that for a (non-zero) finite projection $E$ in $\M$,
$\tau_E(ETE)=\frac{\tau(ETE)}{\tau(E)}$ is the induced trace on
$\M_E=E\M E$.

\begin{Lemma}\label{L:tracial diagonal operators} Let $(\M,\tau)$ be a finite von Neumann algebra
satisfying the weak Dixmier property and let $|\!|\!|\cdot|\!|\!|$
be a tracial gauge  norm on $\M$.
 Suppose  $T,E_1,\cdots,E_n\in \M$, $T\geq 0$,
$E_1+\cdots+E_n=1$. Then $|\!|\!|T|\!|\!|\geq
|\!|\!|\tau_{E_1}(E_1TE_1)E_1+\cdots+\tau_{E_n}(E_nTE_n)E_n|\!|\!|$.
\end{Lemma}
\begin{proof} We may assume that $|\!|\!|1|\!|\!|=1$.
 Since $\M$ satisfies the weak Dixmier property, by
Corollary~\ref{C:local Dixmier  property}, $(\M_{E_i},\tau_{E_i})$
also satisfies the weak Dixmier property, $1\leq i\leq n$. Let
$\epsilon>0$. There are operators $S_{i1},\cdots,S_{ik}$ in
$\M_{E_i}$ such that $E_iTE_i, S_{i1},\cdots, S_{ik}$ are
equi-measurable and
\[\left\|\frac{S_{i1}+\cdots+S_{ik}}{k}-\tau_{E_i}(E_iTE_i)E_i\right\|<\epsilon.\]
Let $S_j=S_{1j}E_1+\cdots+S_{nj}E_n$, $1\leq j\leq k$. Then $T, S_1,
\cdots, S_n$ are equi-measurable and
\[\left\|\frac{S_1+\cdots+S_k}{k}-(\tau_{E_1}(E_1TE_1)E_1+\cdots+
\tau_{E_n}(E_nTE_n)E_n)\right\|<\epsilon.
\]
By Corollary~\ref{C:compare gauge norm with operator norm},
\[{|\!|\!|}\frac{S_1+\cdots+S_k}{k}-(\tau_{E_1}(E_1TE_1)E_1+\cdots+
\tau_{E_n}(E_nTE_n)E_n){|\!|\!|}<\epsilon.
\]
Hence,
$|\!|\!|\tau_{E_1}(E_1TE_1)E_1+\cdots+\tau_{E_n}(E_nTE_n)E_n|\!|\!|\leq
|\!|\!|T|\!|\!|+\epsilon$. Since $\epsilon>0$ is arbitrary, we
obtain the lemma.
\end{proof}

\begin{Corollary}\label{C:Conditinal expectation to finite dimensional algebra}
Let $(\M,\tau)$ be a finite von Neumann algebra satisfying  the weak
Dixmier property and $|\!|\!|\cdot|\!|\!|$ be a tracial gauge  norm
on $\M$. If $\A$ is a finite-dimensional abelian von Neumann
subalgebra of $\M$ and $\mathbf{E}_\A$ is the normal conditional
expectation from $\M$ onto $\A$ preserving $\tau$, then for every
$T\in \M$, $|\!|\!|\mathbf{E}_{\A}(T)|\!|\!|\leq |\!|\!|T|\!|\!|$.
\end{Corollary}
\begin{proof} Let $\A=\{E_1,\cdots, E_n\}''$ such that
$E_1+\cdots+E_n=1$. Then for every $T\in \M$,
$\mathbf{E}_\A(T)=\tau_{E_1}(E_1TE_1)E_1+\cdots+\tau_{E_n}(E_nTE_n)E_n$.
By Corollary~\ref{C:diagonal operators} and Lemma~\ref{L:tracial
diagonal operators}, $|\!|\!|\mathbf{E}_{\A}(T)|\!|\!|\leq
|\!|\!|T|\!|\!|$.
\end{proof}

\begin{proof}[Proof of Theorem~{\bf B}]
 By Lemma~\ref{L:tracial implies symmetric} and
Corollary~\ref{C:symmetric implies unitarily invariant},
$|\!|\!|\cdot|\!|\!|$ is unitarily invariant. Suppose $T_\alpha$ is
a net in $\M_{1, |\!|\!|\cdot|\!|\!|}$ such that $\lim_\alpha
T_\alpha =T$ in the weak operator topology.  Let $T=U|T|$ be the
polar decomposition of $T$. Then $\lim_\alpha U^*T_\alpha=|T|$ in
the weak operator topology. Since $|\!|\!|\cdot|\!|\!|$ is unitarily
invariant, $|\!|\!|UT_\alpha|\!|\!|\leq 1$ and $|\!|\!|\,
|T|\,|\!|\!|=|\!|\!|T|\!|\!|$. So we may assume that $T\geq 0$ and
$T_\alpha=T_\alpha^*$. By the spectral decomposition theorem and
Corollary~\ref{C:compare gauge norm with operator norm}, to prove
$|\!|\!|T|\!|\!|\leq 1$, we need to prove $|\!|\!|S|\!|\!|\leq 1$
for every positive simple operator $S$ such that $S\leq T$. Let
$S=a_1E_1+\cdots+a_nE_n$ and $\epsilon>0$. Since $\lim_\alpha
T_\alpha =T\geq S$, $\lim_\alpha E_iT_\alpha E_i=E_iT E_i\geq a_i
E_i$ for $1\leq i\leq n$. Hence,
$\lim_\alpha\tau_{E_i}(E_i(T_\alpha+\epsilon) E_i)\geq
a_i+\epsilon>a_i$. So there is a $\beta$ such that
$\tau_{E_1}(E_1(T_\beta+\epsilon)
E_1)E_1+\cdots+\tau_{E_n}(E_n(T_\beta+\epsilon) E_n)E_n\geq S$. By
Lemma~\ref{L:tracial diagonal operators} and Corollary~\ref{C:S<T},
$1+\epsilon\geq |\!|\!|T_\beta+\epsilon|\!|\!|\geq
|\!|\!|\tau(E_1(T_\beta+\epsilon)E_1)E_1+\cdots+\tau(E_n(T_\beta+\epsilon)E_n)E_n|\!|\!|\geq
|\!|\!|S|\!|\!|$. Since $\epsilon>0$ is arbitrary,
$|\!|\!|S|\!|\!|\leq 1$.
\end{proof}

\begin{proof}[Proof of Corollary~{\bf 1}]
 Since $\A$ is a
separable abelian von Neumann algebra, there is a sequence of finite
dimensional abelian von Neumann subalgebras $\A_n$ such that
$\A_1\subset \A_2\subset\cdots\subset \A$ and $\A$ is the closure of
$\cup_n\A_n$ in the strong operator topology. Let
$\mathbf{E}_{\A_n}$ be the normal conditional expectation from $\M$
onto $\A_n$ preserving $\tau$. Then for every $T\in \M$,
$\mathbf{E}_\A(T)=\lim_{n\rightarrow\infty}\mathbf{E}_{\A_n}(T)$ in
the strong operator topology. By Theorem~{\bf B} and
Corollary~\ref{C:Conditinal expectation to finite dimensional
algebra}, $|\!|\!|\mathbf{E}_{\A}(T)|\!|\!|\leq |\!|\!|T|\!|\!|$.
\end{proof}

In the following we give some other useful corollaries of
Theorem~{\bf B}.

\begin{Corollary}\label{C:Monotone convergence} Let $(\M,\tau)$ be a finite von Neumann algebra satisfying
the weak Dixmier property and let $|\!|\!|\cdot|\!|\!|$ be a
tracial gauge norm on $\M$. Suppose $0\leq T_1\leq
T_2\leq\cdots\leq T$ in $\M$ such that $\lim_{n\rightarrow \infty}
T_n=T$ in the weak operator topology. Then $\lim_{n\rightarrow
\infty} |\!|\!|T_n|\!|\!|=|\!|\!|T|\!|\!|$.
\end{Corollary}
\begin{proof} By Corollary~\ref{C:S<T},
 $|\!|\!|T_1|\!|\!|\leq |\!|\!|T_2|\!|\!|\leq \cdots
\leq |\!|\!|T|\!|\!|$. Hence, $\lim_{n\rightarrow \infty}
|\!|\!|T_n|\!|\!|\leq |\!|\!|T|\!|\!|$. By Theorem~{\bf B},
$\lim_{n\rightarrow \infty} |\!|\!|T_n|\!|\!|\geq |\!|\!|T|\!|\!|$.
\end{proof}

\begin{Corollary}\label{C: type I_n operators are dense} Let $(\M,\tau)$ be a
finite von Neumann algebra satisfying the weak Dixmier property and
$|\!|\!|\cdot|\!|\!|_1$ and $|\!|\!|\cdot|\!|\!|_2$ be two tracial
gauge norms on $\M$. Then
$|\!|\!|\cdot|\!|\!|_1=|\!|\!|\cdot|\!|\!|_2$ on $\M$ if
$|\!|\!|T|\!|\!|_1=|\!|\!|T|\!|\!|_2$ for every operator
$T=a_1E_1+\cdots+a_nE_n$ in $\M$ such that  $a_1,\cdots,a_n\geq 0$
and $\tau(E_1)=\cdots=\tau(E_n)=\frac{1}{n}$, $n=1,2,\cdots$.
\end{Corollary}
\begin{proof} We need only to prove
$|\!|\!|T|\!|\!|_1=|\!|\!|T|\!|\!|_2$ for every positive operator
$T$ in $\M$. By Theorem~\ref{T:the weak Dixmier property}, $\M$ is
either a finite dimensional von Neumann algebra such that
$\tau(E)=\tau(F)$ for arbitrary two nonzero minimal projections in
$\M$ or $\M$ is diffuse. If $\M$ is a finite dimensional von
Neumann algebra such that $\tau(E)=\tau(F)$ for every pair of
nonzero minimal projections in $\M$, then the corollary is
obvious.  If $\M$ is diffuse, by the spectral decomposition
theorem, there is a sequence of operators $T_n\in\M$ satisfying
the following conditions:
\begin{enumerate}
\item $0\leq T_1\leq T_2\leq\cdots\leq T$,
\item $T_n=a_{n1}E_{n1}+\cdots+a_{nn}E_{nn}$,
$a_{n1},\cdots,a_{nn}\geq 0$ and
$\tau(E_{n1})=\cdots=\tau(E_{nn})=\frac{1}{n}$,
\item $\lim_{n\rightarrow \infty} T_n=T$ in the weak operator
topology. \end{enumerate} By the assumption of the corollary,
$|\!|\!|T_n|\!|\!|_1=|\!|\!|T_n|\!|\!|_2$. By
Corollary~\ref{C:Monotone convergence},
$|\!|\!|T|\!|\!|_1=|\!|\!|T|\!|\!|_2$.
\end{proof}

\begin{Corollary}\label{C: type I_n subfactor are dense} Let $\M$ be a type ${\rm II}\sb 1$ factor and
$|\!|\!|\cdot|\!|\!|_1$ and $|\!|\!|\cdot|\!|\!|_2$  be two
unitarily invariant norms on $\M$. Then
$|\!|\!|\cdot|\!|\!|_1=|\!|\!|\cdot|\!|\!|_2$ on $\M$ if
$|\!|\!|\cdot|\!|\!|_1=|\!|\!|\cdot|\!|\!|_2$ on all type ${\rm
I}_n$ subfactors of $\M$, $n=1,2,\cdots$.
\end{Corollary}

\section{Ky Fan norms on finite von Neumann algebras}
Let $(\M,\tau)$ be a finite von Neumann algebra and $0\leq t\leq 1$.
For $T\in \M$, define the \emph{Ky Fan $t$-th norm} by
\[|\!|\!|T|\!|\!|_{(t)}=\left\{\begin{array}{ll}
                                    \|T\|, & \hbox{$t=0$;} \\
                                     \frac{1}{t}\int_0^t \mu_s(T)ds, & \hbox{$0<t\leq 1$.}
                                   \end{array}
                                 \right.
\]

Let $\M_1=(\M,\tau)*(\mathcal{L}_{\mathcal{F}_2}, \tau')$ be the
reduced free product von Neumann algebra of $\M$ and the free group
factor $\mathcal{L}_{\mathcal{F}_2}$. Then $\M_1$ is a type ${\rm
II}\sb 1$ factor with a faithful normal trace $\tau_1$ such that the
restriction of $\tau_1$ to $\M$ is $\tau$. Recall that $\U(\M_1)$ is
the set of unitary operators in $\M_1$ and $\P(\M_1)$ is the set of
projections in $\M_1$.

\begin{Lemma}\label{L:characterization of Ky fan norms}
 For $0<t\leq 1$, $t|\!|\!|T|\!|\!|_{(t)}=\sup\{|\tau_1(UTE)|:\,
U\in\U(\M_1),\, E\in\P(\M_1),\,  \tau_1(E)=t\}$.
\end{Lemma}
\begin{proof} We may assume that $T$ is a positive operator.
Let $\A$ be a separable diffuse abelian von Neumann subalgebra of
$\M_1$ containing $T$ and $\alpha$ be a $\ast$-isomorphism from
$(\A,\tau_1)$ onto $(L^\infty[0,1], \int_0^1\,dx)$ such that
$\tau_1=\int_0^1dx\cdot \alpha$. Let $f(x)=\alpha(T)$ and $f^*(x)$
be the non-increasing rearrangement of $f(x)$. Then
$\mu_s(T)=f^*(s)$. By the definition of
$f^*$(see~(\ref{E:nonincreasing})),
\[m(\{f^*>f^*(t)\})=
\lim_{n\rightarrow\infty}m\left(\left\{f^*>f^*(t)+\frac{1}{n}\right\}\right)\leq
t\] and
\[m(\{f^*\geq f^*(t)\})\geq \lim_{n\rightarrow\infty} m\left(\left\{f^*>
f^*(t)-\frac{1}{n}\right\}\right)\geq t.\] Since $f^*$ and $f$ are
equi-measurable, $m(\{f>f^*(t)\})\leq t$ and $m(\{f\geq
f^*(t)\})\geq t$. Therefore, there is a measurable subset $A$ of
$[0,1]$, $\{f>f^*(t)\}\subset A\subset \{f\geq f^*(t)\}$, such that
$m(A)=t$. Since $f(x)$ and $f^*(x)$ are equimeasurable, $\int_A
f(s)ds=\int_0^t f^*(s)ds$.  Let $E'=\alpha^{-1}(\chi_A)$. Then
$\tau_1(E')=t$ and $\tau_1(TE')=\int_A f(s)ds=\int_0^t
f^*(s)ds=t|\!|\!|T|\!|\!|_{(t)}$. Hence, $t|\!|\!|T|\!|\!|_{(t)}\leq
\sup\{|\tau_1(UTE)|:\, U\in\U(\M_1), E\in\P(\M_1),\,  \tau_1(E)=t\}$. \\

We need to prove that if $E$ is a projection in $\M_1$,
$\tau_1(E)=t$, and $U\in \U(\M_1)$, then $t|\!|\!|T|\!|\!|_{(t)}\geq
|\tau_1(UTE)|$.
 By the Schwartz
inequality, $|\tau_1(UTE)|=\tau_1(EUT^{1/2}T^{1/2}E)\leq
\tau_1(U^*EUT)^{1/2}\tau_1(ET)^{1/2}$. By
Corollary~\ref{C:s-number}, $\tau_1(ET)=\int_0^1 \mu_s(ET)ds$. By
Corollary~\ref{C:property of s-number}, $\mu_s(ET)\leq
\min\{\mu_s(T), \mu_s(E)\|T\|\}$. Note that $\mu_s(E)=0$ for $s\geq
\tau_1(E)=t$. Hence, $\tau_1(ET)\leq
\int_0^t\mu_s(T)ds=t|\!|\!|T|\!|\!|_t$. Similarly,
$\tau_1(U^*EUT)\leq t|\!|\!|T|\!|\!|_t$.  So $|\tau_1(UTE)|\leq
t|\!|\!|T|\!|\!|_t$. This proves that $t|\!|\!|T|\!|\!|_{(t)}\geq
\sup\{|\tau_1(UTE)|:\, U\in\U(\M), E\in\P(\M_1),\,  \tau_1(E)=t\}$.
\end{proof}
\begin{Theorem}\label{T:Ky fan norms} For $0\leq t\leq 1$, $|\!|\!|\cdot|\!|\!|_{(t)}$ is a
normalized tracial gauge norm on $(\M,\tau)$.
\end{Theorem}
\begin{proof} We  only prove  the triangle inequality, since the other parts
are obvious. We may assume that $0<t\leq 1$.  Let $S, T\in \M$. By
Lemma~\ref{L:characterization of Ky fan norms},
$t|\!|\!|S+T|\!|\!|_{(t)}=\sup\{|\tau_1(U(S+T)E)|:\, U\in\U(\M_1),\,
E\in \P(\M_1),\, \tau_1(E)=t\}\leq \sup\{|\tau_1(USE)|:\,
U\in\U(\M_1),\, E\in \P(\M_1), \tau_1(E)=t\}+ \sup\{|\tau_1(UTE)|:\,
U\in\U(\M_1),\, E\in \P(\M_1),\,
\tau_1(E)=t\}=t|\!|\!|S|\!|\!|_{(t)}+t|\!|\!|T|\!|\!|_{(t)}$.
\end{proof}

\begin{Proposition}\label{P: Monotone property of Ky Fan norms} Let
$(\M,\tau)$ be a finite von Neumann algebra and $T\in (\M,\tau)$.
Then $|\!|\!|T|\!|\!|_{(t)}$ is a non-increasing continuous function
on $[0,1]$.
\end{Proposition}
\begin{proof} Let $0< t_1<t_2\leq 1$.
$|\!|\!|T|\!|\!|_{(t_1)}-|\!|\!|T|\!|\!|_{(t_2)}=\frac{1}{t_1}\int_0^{t_1}
\mu_s(T)ds- \frac{1}{t_2}\int_0^{t_2}
\mu_s(T)ds=\frac{\frac{1}{t_1}\int_0^{t_1}
\mu_s(T)ds-\frac{1}{t_2-t_1}\int_{t_1}^{t_2}\mu_s(T)ds}{t_2(t_2-t_1)}\geq
0.$ Since $\mu_s(T)$ is right-continuous, $|\!|\!|T|\!|\!|_{(t)}$
is a non-increasing continuous function on $[0,1]$.
\end{proof}

\begin{Example}\label{E:Ky fan norms for MnC} \emph{The Ky Fan $\frac{k}{n}$-th  norm
of a matrix $T\in(\M_n(\cc),\tau_n)$ is
\[|\!|\!|T|\!|\!|_{(\frac{k}{n})}=\frac{s_1(T)+\cdots+s_k(T)}{k},\quad
1\leq k\leq n.
\]}
\end{Example}

\section{Dual norms of tracial gauge norms
on finite von Neumann algebras satisfying the weak Dixmier property}

\subsection{Dual norms}
Let $|\!|\!|\cdot|\!|\!|$ be a norm on a finite von Neumann algebra
$(\M,\tau)$. For $T\in \M$, define
\[|\!|\!|T|\!|\!|^\#_\M=\sup\{|\tau(TX)|:\, X\in\M,\, |\!|\!|X|\!|\!|\leq 1\}.
\]
When no confusion arises, we simply write $|\!|\!|\cdot|\!|\!|^\#$
instead of $|\!|\!|\cdot|\!|\!|^\#_\M$.

\begin{Lemma} $|\!|\!|\cdot|\!|\!|^\#$ is a norm on $\M$.
\end{Lemma}
\begin{proof} If $T\neq 0$, $|\!|\!|T|\!|\!|^\#\geq \tau(TT^*)/|\!|\!|T^*|\!|\!|>0$.
It is easy to see that $|\!|\!|\lambda T|\!|\!|^\#=|\lambda|\cdot
|\!|\!|T|\!|\!|^\#$ and $|\!|\!|T_1+T_2|\!|\!|^\#\leq
|\!|\!|T_1|\!|\!|^\#+|\!|\!|T_2|\!|\!|^\#$.
\end{proof}

\begin{Definition}\label{D:dual norm}
\emph{ $|\!|\!|\cdot|\!|\!|^\#$ is called the \emph{dual norm} of
$|\!|\!|\cdot|\!|\!|$ on $\M$ with respect to $\tau$.}
\end{Definition}

The next lemma follows directly from the definition of dual norm.
\begin{Lemma}\label{L:Holder inequality for bounded operators}
Let $|\!|\!|\cdot|\!|\!|$ be a norm on a finite von Neumann
algebra $(\M,\tau)$ and let $|\!|\!|\cdot|\!|\!|^\#$ be the dual
norm on $\M$. Then for $S,T\in \M$, $|\tau(ST)|\leq
|\!|\!|S|\!|\!|\cdot |\!|\!|T|\!|\!|^\#.$
\end{Lemma}

The following corollary is a generalization of
H$\ddot{\text{o}}$lder's inequality for bounded operators in finite
von Neumann algebras.
\begin{Corollary}\label{C:Holder inequality for bounded operators}
 Let $(\M,\tau)$ be a finite von Neumann algebra
 and let $|\!|\!|\cdot|\!|\!|$ be a gauge norm on
$\M$. Then for $S,T\in \M$, $\|ST\|_1\leq |\!|\!|S|\!|\!|\cdot
|\!|\!|T|\!|\!|^\#$.
\end{Corollary}
\begin{proof} By Lemma~\ref{L:1-norm}, $\|ST\|_1=\sup\{|\tau(UST)|:\,
U\in \U(\M)\}$. By Lemma~\ref{L:Holder inequality for bounded
operators} and Lemma~\ref{L:submultiplicative}, $|\tau(UST)|\leq
|\!|\!|US|\!|\!|\cdot |\!|\!|T|\!|\!|^\#=|\!|\!|S|\!|\!|\cdot
|\!|\!|T|\!|\!|^\#.$
\end{proof}

\begin{Proposition}\label{P:dual unitarily invaraint norm}
 If $|\!|\!|\cdot|\!|\!|$ is a unitarily invariant norm on a finite von Neumann algebra  $(\M,\tau)$, then
 $|\!|\!|\cdot|\!|\!|^\#$ is also a unitarily invariant norm on
$\M$.
\end{Proposition}
\begin{proof}   Let $U$ be a unitary operator. Then
$|\!|\!|UT|\!|\!|^\#=\sup\{|\tau(UTX)|:\, X\in\M,\,
|\!|\!|X|\!|\!|\leq 1\}= \sup\{|\tau(TXU)|:\, X\in\M,\,
|\!|\!|X|\!|\!|\leq 1\}=\sup\{|\tau(TX)|:\, X\in\M,\,
|\!|\!|X|\!|\!|\leq 1\}=|\!|\!|T|\!|\!|$ and
$|\!|\!|TU|\!|\!|^\#=\sup\{|\tau(TUX)|:\, X\in\M,\,
|\!|\!|X|\!|\!|\leq 1\}= \sup\{|\tau(TX)|:\, X\in\M,\,
|\!|\!|X|\!|\!|\leq 1\}=|\!|\!|T|\!|\!|$.
\end{proof}

\begin{Proposition}\label{P:dual symmetric gauge norm}
 If $|\!|\!|\cdot|\!|\!|$ is a symmetric gauge
norm on a finite von Neumann algebra $(\M,\tau)$, then
$|\!|\!|\cdot|\!|\!|^\#$ is also a symmetric gauge norm on
$(\M,\tau)$.
\end{Proposition}
\begin{proof} Let $\theta\in \Aut(\M,\tau)$. Then
$|\!|\!|\theta(T)|\!|\!|^\#=\sup\{|\tau(\theta(T)X)|:\, X\in\M,\,
|\!|\!|X|\!|\!|\leq 1\}= \sup\{|\tau(\theta(T\theta^{-1}(X)))|:\,
X\in\M,\, |\!|\!|X|\!|\!|\leq 1\}=\sup\{|\tau(T\theta^{-1}(X))|:\,
X\in\M,\, |\!|\!|X|\!|\!|\leq 1\}=\sup\{|\tau(TX)|:\, X\in\M,\,
|\!|\!|X|\!|\!|\leq 1\}=|\!|\!|T|\!|\!|.$
\end{proof}

\begin{Lemma}\label{L:dual norm of positive operators}
Let $(\M,\tau)$ be a finite von Neumann algebra satisfying  the
weak Dixmier property and let $|\!|\!|\cdot|\!|\!|$ be a tracial
gauge norm on $\M$. If  $T\in \M$ is a positive operator,  then
\[|\!|\!|T|\!|\!|^\#=\sup\{\tau(TX):\,\,X\in \M,\, X\geq 0,\, XT=TX,\, |\!|\!|X|\!|\!|\leq 1\}.
\]
\end{Lemma}
\begin{proof}Let $\A$ be a separable abelian von Neumann subalgebra of
$\M$ containing $T$ and $\mathbf{E}_\A$ be the normal conditional
expectation from $\M$ onto $\A$ preserving $\tau$. For every $Y\in
\M$ such that $|\!|\!|Y|\!|\!|\leq 1$, let $X=\mathbf{E}_\A(Y)$. By
Corollary~{\bf 1}, $|\!|\!|\,\,|X|\,\,|\!|\!|=|\!|\!|X|\!|\!|\leq
|\!|\!|Y|\!|\!|\leq 1$. Furthermore,
$|\tau(TY)|=|\tau(\mathbf{E}_\A(TY))|=|\tau(T\mathbf{E}_\A(Y))|=
|\tau(TX)|\leq \tau(T|X|)$. Hence,
\[|\!|\!|T|\!|\!|^\#=\sup\{\tau(TX):\,\,X\in \M,\, X\geq 0,\, XT=TX,\, |\!|\!|X|\!|\!|\leq 1\}.\]
\end{proof}

\begin{Lemma}\label{L:Dual norm of simple operators} Let $(\M,\tau)$ be a finite von Neumann algebra
satisfying the weak Dixmier property and let $|\!|\!|\cdot|\!|\!|$
be a tracial gauge norm on $\M$. Suppose $T=a_1 E_1+\cdots+a_n
E_n$ is a positive simple operator in $\M$. Then
\begin{equation*}
\begin{split}
|\!|\!|T|\!|\!| &=\sup\left\{\tau(TX):\,\, X=b_1E_1+\cdots+b_nE_n
\geq 0 \text{ and } |\!|\!|X|\!|\!|^\#\leq 1\right\}
\\
&=\sup\left\{\sum_{k=1}^n a_kb_k\tau(E_k):\,\,
X=b_1E_1+\cdots+b_nE_n \geq 0 \text{ and } |\!|\!|X|\!|\!|^\#\leq
1\right\}.
\end{split}
\end{equation*}
\end{Lemma}
\begin{proof} By Lemma~\ref{L:dual norm of positive operators},
 $|\!|\!|T|\!|\!|^\#=\sup\{|\tau(TX)|:\, X\in\M,\,  X\geq 0,\, XT=TX,\,
|\!|\!|X|\!|\!|\leq 1\}$. Let $\A=\{E_1,\cdots,E_n\}''$ and
$\mathbf{E}_\A$ be the normal conditional expectation from $\M$ onto
$\A$ preserving $\tau$. Then
$S=\mathbf{E}_\A(X)=\tau_{E_1}(E_1XE_1)E_1+\cdots+\tau_{E_n}(E_nXE_n)E_n$
is a positive operator,
$\tau(TX)=\tau(\mathbf{E}_\A(TX))=\tau(T\mathbf{E}_\A(X))=\tau(TS)$,
and $|\!|\!|S|\!|\!|\leq |\!|\!|X|\!|\!|$ by
Corollary~\ref{C:Conditinal expectation to finite dimensional
algebra}. Combining the definition of dual norm, this proves the
lemma.
\end{proof}

\begin{Corollary}\label{C:Dual norms of equimeasurable simple operators}Let $(\M,\tau)$ be a finite von Neumann algebra
satisfying the weak Dixmier property and let $|\!|\!|\cdot|\!|\!|$
be a tracial gauge norm on $\M$. Suppose $S, T$ are
equi-measurable, positive simple operators in $\M$. Then
$|\!|\!|S|\!|\!|^\#=|\!|\!|T|\!|\!|^\#$.
\end{Corollary}

\begin{Theorem}\label{T:dual norms of tracial}
 Let $(\M,\tau)$ be a finite von Neumann algebra
satisfying the weak Dixmier property and let $|\!|\!|\cdot|\!|\!|$
be a tracial gauge norm on $\M$. Then $|\!|\!|\cdot|\!|\!|^\#$ is
also a tracial gauge norm on $\M$. Furthermore, if
$|\!|\!|1|\!|\!|=1$, then $|\!|\!|1|\!|\!|^\#=1$.
\end{Theorem}
\begin{proof}By Lemma~\ref{L:tracial implies symmetric},
$|\!|\!|\cdot|\!|\!|$ is a symmetric gauge norm on $\M$.
 By Proposition~\ref{P:dual symmetric gauge norm},
Corollary~\ref{C:Dual norms of equimeasurable simple operators} and
 Lemma~\ref{L:sufficient condition for tracial},
 $|\!|\!|\cdot|\!|\!|^\#$ is a tracial gauge norm on $\M$.
Note that $|\!|\!|1|\!|\!|=1$, hence, $|\!|\!|1|\!|\!|^\#\geq
\tau(1\cdot 1)=1$. On the other hand, by Theorem~\ref{T:tracial},
$|\!|\!|1|\!|\!|^\#=\sup\{|\tau(X)|:\, X\in\M,\, |\!|\!|X|\!|\!|\leq
1\}\leq \sup\{|\!|\!|X|\!|\!|:\, X\in\M,\, |\!|\!|X|\!|\!|\leq
1\}\leq 1$.
\end{proof}

\begin{Corollary}\label{C:Dual norms of von Neumann subalgebra}
 Let $(\M,\tau)$ be a finite von Neumann algebra
satisfying the weak Dixmier property and let $|\!|\!|\cdot|\!|\!|$
be a tracial gauge norm on $\M$. If $\N$ is a von Neumann
subalgebra of $\M$ satisfying the weak Dixmier property, then
$|\!|\!|\cdot|\!|\!|_\N^\#$ is the restriction of
$|\!|\!|\cdot|\!|\!|_\M^\#$ to $\N$.
\end{Corollary}
\begin{proof} Let $|\!|\!|\cdot|\!|\!|_1=|\!|\!|\cdot|\!|\!|_\N^\#$ and
$|\!|\!|\cdot|\!|\!|_2$ be the restriction of
$|\!|\!|\cdot|\!|\!|_\M^\#$ to $\N$. By Theorem~\ref{T:dual norms of
tracial}, both $|\!|\!|\cdot|\!|\!|_1$ and $|\!|\!|\cdot|\!|\!|_2$
are tracial gauge norms on $\N$. By Lemma~\ref{C:simple operators
are dense}, to prove $|\!|\!|\cdot|\!|\!|_1=|\!|\!|\cdot|\!|\!|_2$,
we need to prove $|\!|\!|T|\!|\!|_1=|\!|\!|T|\!|\!|_2$ for every
positive simple operator $T\in \N$. Let $\A$ be a finite dimensional
abelian von Neumann subalgebra of $\N$ containing $T$.  By
Lemma~\ref{L:Dual norm of simple operators},
$|\!|\!|T|\!|\!|_\M^\#=|\!|\!|T|\!|\!|_\N^\#=|\!|\!|T|\!|\!|_\A^\#$.
So $|\!|\!|T|\!|\!|_1=|\!|\!|T|\!|\!|_2$.
\end{proof}

\subsection{Dual norms of Ky Fan norms}
For $(x_1,\cdots,x_n)\in \cc^n$, $\tau(x)=\frac{x_1+\cdots+x_n}{n}$
defines a trace on $\cc^n$. For $1\leq k\leq n$, the Ky Fan
$\frac{k}{n}$-th  norm on $(\cc^n,\tau)$ is
$|\!|\!|(x_1,\cdots,x_n)|\!|\!|_{(\frac{k}{n})}=
\frac{x_1^*+\cdots+x_k^*}{k}$, where $(x_1^*,\cdots, x_n^*)$ is the
decreasing rearrangement of $(|x_1|,\cdots,|x_n|)$. Let
$\mathbf{\Gamma}=\{(x_1,\cdots,x_n)\in \cc^n:\, x_1\geq x_2\geq
x_k=x_{k+1}= \cdots=x_n\geq 0, \frac{x_1+\cdots+x_k}{k}\leq 1\}$ and
$\mathcal{E}$ be the set of extreme points of $\mathbf{\Gamma}$.\\

The proof of the following lemma is an easy exercise.
\begin{Lemma}\label{L:Extreme points} $\mathcal{E}$ consists of  $k+1$ points:
$\left(k,0,\cdots\right)$, $\left(\frac{k}{2}, \frac{k}{2},
0,\cdots\right)$, $\cdots, (\frac{k}{k-1},\\
 \cdots,
\frac{k}{k-1}, 0,\cdots), (1,1,\cdots,1)$ and $(0,0,\cdots,0)$.
\end{Lemma}

The following lemma is well-known. For a proof we refer to 10.2
of~\cite{HLP}.
\begin{Lemma}\label{L:HLP} Let $s_1\geq s_2\geq\cdots\geq s_n\geq 0$ and
$t_1,\cdots,t_n\geq 0$. If $t_1^*\geq t_2^*\geq \cdots \geq t_n^*$
is the decreasing rearrangement of $t_1,\cdots,t_n$, then
$s_1t_1^*+\cdots+s_nt_n^*\geq s_1t_1+\cdots+s_nt_n$.
\end{Lemma}

\begin{Lemma}\label{L:dual of ky fan norms for matrix} For $T\in (M_n(\cc),\tau_n)$,
\[|\!|\!|T|\!|\!|_{(\frac{k}{n})}^\#=\max\left\{\frac{k}{n}\|T\|,
\|T\|_1\right\}.
\]
\end{Lemma}
\begin{proof} Let
$|\!|\!|T|\!|\!|_1=|\!|\!|T|\!|\!|_{(\frac{k}{n})}^\#$ and
$|\!|\!|T|\!|\!|_2=\max\{\frac{k}{n}\|T\|,\|T\|_1\}$. Then both
$|\!|\!|\cdot|\!|\!|_1$ and $|\!|\!|\cdot|\!|\!|_2$ are unitarily
invariant norms on $M_n(\cc)$. To prove
$|\!|\!|\cdot|\!|\!|_1=|\!|\!|\cdot|\!|\!|_2$, we need only to prove
$|\!|\!|T|\!|\!|_1=|\!|\!|T|\!|\!|_2$ for every positive matrix $T$
in $M_n(\cc)$. We can assume that $T=\left(\begin{array}{ccc}
s_1&&\\
&\ddots&\\
&&s_n
\end{array}\right)$, where $s_1,\cdots, s_n$ are $s$-numbers
of $T$ such that $s_1\geq s_2\geq \cdots \geq s_n$. By
Lemma~\ref{L:Dual norm of simple operators} and~\ref{L:HLP},
\[|\!|\!|T|\!|\!|_1=\sup\left\{\frac{\sum_{i=1}^n s_it_i}{n}:\,\,
(t_1,\cdots,t_n)\in
\mathbf{\Gamma}\right\}=\sup\left\{\frac{\sum_{i=1}^n
s_it_i}{n}:\,\, (t_1,\cdots,t_n)\in \mathcal{E}\right\}.
\]
 Note that $\|T\|=s_1\geq
s_2\geq \cdots\geq s_n\geq 0$.  By Lemma~\ref{L:Extreme points} and
simple computations,
$|\!|\!|T|\!|\!|_1=\max\left\{\frac{k}{n}\|T\|,\|T\|_1\right\}=|\!|\!|T|\!|\!|_2.$
\end{proof}

The next lemma simply follows from the definition of dual norms.
\begin{Lemma}\label{L:comparison of dual norms} Let $(\M,\tau)$ be a finite von Neumann algebra and
$|\!|\!|\cdot|\!|\!|,\, |\!|\!|\cdot|\!|\!|_1,
\,|\!|\!|\cdot|\!|\!|_2$ be norms on $\M$ such that
\[|\!|\!|T|\!|\!|_1\leq|\!|\!|T|\!|\!|\leq  |\!|\!|T|\!|\!|_2,\quad \forall T\in \M.
\] Then
\[|\!|\!|T|\!|\!|_2^\#\leq|\!|\!|T|\!|\!|^\#\leq  |\!|\!|T|\!|\!|_1^\#,\quad \forall T\in \M.
\]
\end{Lemma}

\begin{Corollary}\label{C:equivalent dual norms} Let $(\M,\tau)$ be a finite von Neumann algebra and
$|\!|\!|\cdot|\!|\!|_1, \,|\!|\!|\cdot|\!|\!|_2$ be  equivalent
norms on $\M$. Then $|\!|\!|\cdot|\!|\!|_1^\#$ and
$|\!|\!|\cdot|\!|\!|_2^\#$ are equivalent norms on $\M$.
\end{Corollary}

\begin{Theorem}\label{T:dual norms of Ky Fan norms} Let $\M$ be a
type ${\rm II}\sb 1$ factor and $0\leq t\leq 1$. Then
\[|\!|\!|T|\!|\!|_{(t)}^\#=\max\{t\|T\|,\|T\|_1\},\quad\forall
T\in\M.
\]
\end{Theorem}
\begin{proof}Firstly, we assume $t=\frac{k}{n}$ is a rational
number. Let $\N_r$ be a type $I_{rn}$ subfactor of $\M$. Then the
restriction of $|\!|\!|\cdot|\!|\!|_{(t)}$ to $\N_r$ is
$|\!|\!|\cdot|\!|\!|_{(\frac{rk}{rn})}$. By Lemma~\ref{L:dual of ky
fan norms for matrix} and Corollary~\ref{C:Dual norms of von Neumann
subalgebra}, $|\!|\!|T|\!|\!|_{(t)}^\#=\max\{t\|T\|,\|T\|_1\}$ for
$T\in \N_r$. By Corollary~\ref{C: type I_n subfactor are dense},
$|\!|\!|T|\!|\!|_{(t)}^\#=\max\{t\|T\|,\|T\|_1\}$ for all $T\in \M$.
Now assume $t$ is an irrational number. Let  $t_1,t_2$ be two
rational numbers such that $t_1<t<t_2$. By Lemma~\ref{L:comparison
of dual norms}, for every $T\in \M$,
\[\max\{t_2\|T\|,\|T\|_1\}\leq |\!|\!|T|\!|\!|_{(t)}^\#\leq
\max\{t_1\|T\|,\|T\|_1\}.
\] Since  $t_1\leq t\leq t_2$ are arbitrary,
$|\!|\!|T|\!|\!|_{(t)}^\#=\max\{t\|T\|,\|T\|_1\}$.
\end{proof}
\subsection{Proof of Theorem C}
\begin{Lemma}\label{L:Double dual} Let $n\in \mathbb{N}$ and $\tau$
be an arbitrary faithful state on $\cc^n$.  If $|\!|\!|\cdot|\!|\!|$
is a norm on $(\cc^n, \tau)$ and $|\!|\!|\cdot|\!|\!|^\#$ is the
dual norm with respect to $\tau$, then
$|\!|\!|\cdot|\!|\!|^{\#\#}=|\!|\!|\cdot|\!|\!|$.
\end{Lemma}
\begin{proof}By Lemma~\ref{L:Holder
inequality for bounded operators},
$|\!|\!|T|\!|\!|^{\#\#}=\sup\{|\tau(TX)|:\, X\in \cc^n,\,
|\!|\!|X|\!|\!|^\#\leq 1\}\leq |\!|\!|T|\!|\!|$. We need  to prove
$|\!|\!|T|\!|\!|\leq |\!|\!|T|\!|\!|^{\#\#}$. By the Hahn-Banach
Theorem, there is a continuous linear functional $\phi$ on $\cc^n$
with respect to the topology induced by $|\!|\!|\cdot|\!|\!|$ on
$\cc^n$ such that $|\!|\!|T|\!|\!|=\phi(T)$ and $\|\phi\|=1$. Since
all norms on $\cc^n$ induce the same topology, there is an  element
$Y\in \cc^n$ such that $\phi(S)=\tau(SY)$ for all $S\in \cc^n$. By
the definition of dual norm, $|\!|\!|Y|\!|\!|^\#=\|\phi\|=1$. By
Lemma~\ref{L:Holder inequality for bounded operators},
$|\!|\!|T|\!|\!|=\phi(T)=\tau(TY)\leq |\!|\!|T|\!|\!|^{\#\#}$.
\end{proof}

\begin{proof}[Proof of Thereom~{\bf C}] By Theorem~\ref{T:dual norms of
tracial},
  both $|\!|\!|\cdot|\!|\!|^{\#\#}$ and
$|\!|\!|\cdot|\!|\!|$ are tracial gauge norm on $\M$. By
Corollary~\ref{C:simple operators are dense}, to prove
$|\!|\!|\cdot|\!|\!|^{\#\#}=|\!|\!|\cdot|\!|\!|$, we need to prove
that $|\!|\!|T|\!|\!|=|\!|\!|T|\!|\!|^{\#\#}$ for every positive
simple operator $T\in\M$. Let $\A$ be the abelian von Neumann
subalgebra generated by $T$. By Corollary~\ref{C:Dual norms of von
Neumann subalgebra} and Lemma~\ref{L:Double dual},
$|\!|\!|T|\!|\!|_\M^{\#\#}=|\!|\!|T|\!|\!|_\A^{\#\#}=|\!|\!|T|\!|\!|$.
\end{proof}

\section{Proof of Theorem A}

 Let $(\M,\tau)$ be a finite von Neumann
algebra.
\begin{Lemma}\label{L:T_f} Let $n\in \nn$, $a_1\geq a_2\geq \cdots\geq
a_n\geq a_{n+1}=0$ and
$f(x)=a_1\chi_{[0,\frac{1}{n})}(x)+a_2\chi_{[\frac{1}{n},\frac{2}{n})}(x)+\cdots+
 a_n\chi_{[\frac{n-1}{n},1]}(x)$. For $T\in \M$, define
\begin{equation}
|\!|\!|T|\!|\!|_{f}=\int_0^1 f(s)\mu_s(T)ds.
\end{equation}
  Then
\begin{equation}
|\!|\!|T|\!|\!|_{f}=\sum_{k=1}^n \frac{k(a_k - a_{k+1})}{n}
|\!|\!|T|\!|\!|_{\left(\frac{k}{n}\right)}.
\end{equation}
\end{Lemma}
\begin{proof} Since $t|\!|\!|T|\!|\!|_{(t)}=\int_{0}^t \mu_s(T)ds$,
summation by parts shows that
\begin{equation*}
\begin{split}
|\!|\!|T|\!|\!|_{f} &= \int_0^1 f(s)\mu_s(T)dt
=a_1\int_0^{\frac{1}{n}}\mu_s(T)ds+a_2\int_{\frac{1}{n}}^{\frac{2}{n}}\mu_s(T)ds+\cdots+
a_n\int_{\frac{n-1}{n}}^1\mu_s(T)ds \\
&= \sum_{k=1}^n \frac{k(a_k - a_{k+1})}{n}
|\!|\!|T|\!|\!|_{\left(\frac{k}{n}\right)}. \qedhere
\end{split}
\end{equation*}
\end{proof}
\begin{Corollary}\label{L:T_f is a norm} The norm $|\!|\!|\cdot|\!|\!|_{f}$ defined as above is a
tracial  gauge norm on $\M$ and $|\!|\!|1|\!|\!|_{f}=\int_0^1
f(x)dx=\frac{a_1+\cdots+a_n}{n}$.
\end{Corollary}

\begin{Lemma}\label{L:maximal unitarily invariant norms} Let
$(\M,\tau)$ be a finite von Neumann algebra satisfying the weak
Dixmier property and
 $\{|\!|\!|\cdot|\!|\!|_\alpha\}$ be a set of
tracial gauge norms on $(\M,\tau)$ such that
$|\!|\!|1|\!|\!|_\alpha\leq 1$ for all $\alpha$. For every $T\in
\M$, define
\[|\!|\!|T|\!|\!|=\sup_\alpha|\!|\!|T|\!|\!|_\alpha.
\] Then $|\!|\!|\cdot|\!|\!|\triangleq \bigvee_\alpha|\!|\!|\cdot|\!|\!|_\alpha$ is also
a tracial gauge norm on $(\M,\tau)$.
\end{Lemma}
\begin{proof} By Corollary~\ref{C:compare gauge norm with operator
norm}, $|\!|\!|T|\!|\!|\leq \|T\|$ is well defined. It is easy to
check that $|\!|\!|\cdot|\!|\!|$ is a tracial gauge norm on
$(\M,\tau)$.
\end{proof}

\begin{proof}[Proof of Theorem~{\bf A}] Let
\begin{multline*}
\F'=\{\mu_s(X):\, X\in \M,\,
   |\!|\!|X|\!|\!|^\#\leq 1,\, X=b_1F_1+\cdots+b_kF_k\geq 0, \\
\text{where}\, F_1+\cdots+F_k=1\,\text{and}\,
\tau(F_1)=\cdots=\tau(F_k)=\frac{1}{k}, k=1,2,\cdots\}.
\end{multline*}
 For every positive operator $X\in \M$ such that $|\!|\!|X|\!|\!|^\#\leq 1$,  $\int_0^1\mu_s(X)ds =\tau(X)=\|X\|_1\leq
|\!|\!|X|\!|\!|^\#\leq 1$ by Theorem~\ref{T:tracial}. Hence
$\F'\subset \F$ and $\mu_s(1)=\chi_{[0,1]}(s)\in \F'$ by Theorem
6.10. For $T\in \M$, define
\[|\!|\!|T|\!|\!|'=\sup\{|\!|\!|T|\!|\!|_{f}:\, f\in \F'\}.
\] By Corollary~\ref{L:T_f is a norm}, $|\!|\!|\cdot|\!|\!|'$ is a tracial gauge norm on $\M$.
 To prove that
$|\!|\!|\cdot|\!|\!|'=|\!|\!|\cdot|\!|\!|$, by Corollary~\ref{C:
type I_n operators are dense}, we need to prove that
$|\!|\!|T|\!|\!|'=|\!|\!|T|\!|\!|$ for every positive operator $T\in
\M$ such that $T=a_1E_1+\cdots+a_nE_n$ and
$\tau(E_1)=\cdots=\tau(E_n)=\frac{1}{n}$.\\

  By Lemma~\ref{L:Dual
norm of simple operators} and Theorem~{\bf C},
\begin{equation*}
|\!|\!|T|\!|\!|=\sup\left\{\frac{1}{n}\sum_{k=1}^n a_kb_k:\,\,
X=b_1E_1+\cdots+b_nE_n \geq 0 \text{ and } |\!|\!|X|\!|\!|^\#\leq
1\right\}.
\end{equation*}
Note that if $X=b_1E_1+\cdots+b_nE_n$ is a positive simple operator
in $\M$ and $|\!|\!|X|\!|\!|^\#\leq 1$, then $\mu_s(X)\in \F'$ and $
|\!|\!|T|\!|\!|_{\mu_s(X)}=\int_0^1\mu_s(X)\mu_s(T)ds=\frac{1}{n}\sum_{k=1}^n
a_k^*b_k^*$, where $\{a_k^*\}$ and $\{b_k^*\}$ are non-increasing
rearrangements of $\{a_k\}$ and $\{b_k\}$, respectively. By
Lemma~\ref{L:HLP}, $ |\!|\!|T|\!|\!|\leq
\sup\{|\!|\!|T|\!|\!|_{f}:\, f\in\F'\} = |\!|\!|T|\!|\!|'.$
\\

We need to  prove $|\!|\!|T|\!|\!|\geq |\!|\!|T|\!|\!|'$. Let
$X=b_1F_1+\cdots+b_kF_k$ be a positive operator in $\M$ such that
$F_1+\cdots+F_k=1$, $\tau(F_1)=\cdots=\tau(F_k)=\frac{1}{k}$ and
$|\!|\!|X|\!|\!|^\#\leq 1$. We need only to prove that
$|\!|\!|T|\!|\!|\geq |\!|\!|T|\!|\!|_{\mu_s(X)}$.
 Since
$(\M,\tau)$ satisfies the weak Dixmier property, by
Theorem~\ref{T:the weak Dixmier property}, $(\M,\tau)$ is either a
von Neumann subalgebra of $(M_n(\cc),\tau_n)$ that contains all
diagonal matrices or $\M$ is a diffuse von Neumann algebra. In
either case, we may assume that
$T=\tilde{a}_1\tilde{E}_1+\cdots+\tilde{a}_r\tilde{E}_r$ and
$X=\tilde{b}_1\tilde{F}_1+\cdots+\tilde{b}_r\tilde{F}_r$, where
$\tilde{E}_1+\cdots+\tilde{E}_r$ $=$
$\tilde{F}_1+\cdots+\tilde{F}_r=1$ and
$\tau(\tilde{E}_i)=\tau(\tilde{F}_i)=\frac{1}{r}$ for $1\leq i\leq
r$, $\tilde{a}_1\geq \cdots\geq \tilde{a}_r\geq 0$ and
$\tilde{b}_1\geq \cdots\geq \tilde{b}_r\geq 0$. Let
$Y=\tilde{b}_1\tilde{E}_1+\cdots+\tilde{b}_r\tilde{E}_r$. Then $X$
and $Y$ are two equi-measurable operators in $\M$ and
$\mu_s(X)=\mu_s(Y)$.   By Theorem~\ref{T:dual norms of tracial},
$|\!|\!|Y|\!|\!|^\#\leq 1$. By Lemma~\ref{L:Holder inequality for
bounded operators},
\[|\!|\!|T|\!|\!|\geq\tau(TY)=\frac{1}{r}\sum_{i=1}^r
\tilde{a}_i\tilde{b}_i=\int_0^1\mu_s(Y)\mu_s(T)ds=\int_0^1\mu_s(X)\mu_s(T)ds=
|\!|\!|T|\!|\!|_{\mu_s(X)}.\]
\end{proof}

Combining Theorem~{\bf A} and Lemma~\ref{L:unitarily invariant
implies tracial}, we obtain the following corollary.
\begin{Corollary}\label{C:representation theorem for UIN on 21 factor}
 Let $(\M,\tau)$ be a finite factor and let $|\!|\!|\cdot|\!|\!|$ be a normalized unitarily invariant norm on
$\M$. Then there is a subset $\F'$ of $\F$ containing the constant 1
function on $[0,1]$ such that for all $T\in \M$,
$|\!|\!|T|\!|\!|=\sup\{|\!|\!|T|\!|\!|_{f}:\, f\in \F'\}.$
\end{Corollary}

Combining Theorem~{\bf A} and Lemma~\ref{L:symmetric implies
tracial} we obtain the following corollary.
\begin{Corollary}\label{C:representation theorem for UIN on  L(0,1)} Let $|\!|\!|\cdot|\!|\!|$ be a normalized symmetric gauge
norm on $(L^\infty[0,1], \int_0^1dx)$. Then there is a subset $\F'$
of $\F$ containing the constant 1 function on $[0,1]$ such that for
all $T\in L^\infty[0,1]$,
$|\!|\!|T|\!|\!|=\sup\{|\!|\!|T|\!|\!|_{f}:\, f\in \F'\}.$

\end{Corollary}

\section{Proof of Theorem D and Theorem E}

\begin{Lemma}\label{L: two embedding are equivalent}
 Let $\theta_1, \theta_2$ be two embeddings from
$(L^\infty[0,1], \int_0^1dx)$ into a finite von Neumann algebra
$(\M,\tau)$.  If $|\!|\!|\cdot|\!|\!|$ is a tracial gauge norm on
$\M$, then $|\!|\!|\theta_1(f)|\!|\!|=|\!|\!|\theta_2(f)|\!|\!|$ for
every $f\in L^\infty[0,1]$.
\end{Lemma}
\begin{proof}  If $f\in L^\infty[0,1]$ is a positive
function, then $\theta_1(f)$ and $\theta_2(f)$ are equi-measurable
operators in $\M$. Hence
$|\!|\!|\theta_1(f)|\!|\!|=|\!|\!|\theta_2(f)|\!|\!|$.
\end{proof}

\begin{proof}[Proof of Theorem~{\bf D}] We  prove Theorem {\bf D} for
diffuse finite von Neumann algebras. The proof of the atomic case is
similar.  We may assume that the norms on $\M$ or $L^\infty[0,1]$
are normalized. By the definition of Ky Fan norms, there is a
one-to-one correspondence  between  Ky Fan $t$-th norms on $(\M,
\tau)$ and Ky Fan $t$-th norms on $(L^\infty[0,1], \int_0^1dx)$ as
in Theorem~{\bf D}. By Lemma~\ref{L:T_f}, Theorem~\ref{T:the weak
Dixmier property} and Theorem~{\bf A}, there is a one-to-one
correspondence between normalized tracial norms on $(\M,\tau)$ and
normalized symmetric gauge norms on $(L^\infty[0,1], \int_0^1dx)$ as
in Theorem~{\bf D}.
\end{proof}

\begin{Example}\label{E:Lp space} \emph{For $1\leq p\leq \infty$, the $L^p$-norm  on
$L^\infty[0,1]$ defined by
\[\|f(x)\|_p=\left\{
               \begin{array}{ll}
                 \left(\int_0^1|f(x)|^pdx\right)^{1/p}, & \hbox{$1\leq p<\infty$;} \\
               \esssup\, |f|, & \hbox{$p=\infty$}
               \end{array}
             \right.
\] is a normalized symmetric gauge norm on $(L^\infty[0,1],
\int_0^1dx)$. By Corollary~{\bf 2} and Corollary~\ref{C:s-number},
the induced norm
\[\|T\|_p=\left\{
               \begin{array}{ll}
                \left(\tau(|T|^p)\right)^{1/p}= \left(\int_0^1|\mu_s(T)|^pds\right)^{1/p}, & \hbox{$1\leq p<\infty$;} \\
                 \|T\|, & \hbox{$p=\infty$}
               \end{array}
             \right.
\] is a normalized unitarily invariant norm on a type ${\rm II}\sb 1$ factor $\M$.
The norms $\{\|\cdot\|_p:\, 1\leq p\leq \infty\}$ are called
$L^p$-norms on $\M$.}
\end{Example}

\begin{Corollary}\label{C:extension of norms} Let $(\M,\tau)$ be a
finite von Neumann algebra satisfying the weak Dixmier property
and let $|\!|\!|\cdot|\!|\!|$ be a tracial gauge norm on
$(\M,\tau)$. If $(\M,\tau)$ can be embedded into a finite factor
$(\M_1,\tau_1)$, then there is a unitarily invariant norm
$|\!|\!|\cdot|\!|\!|_1$  on $(\M_1,\tau_1)$ such that
$|\!|\!|\cdot|\!|\!|$ is the restriction of
$|\!|\!|\cdot|\!|\!|_1$ to $(\M,\tau)$.
\end{Corollary}

The following example shows that without the weak Dixmier property,
Corollary~\ref{C:extension of norms} may fail.

\begin{Example}\emph{On $(\cc^2,\tau)$,
$\tau((x,y))=\frac{1}{3}x+\frac{2}{3}y$, let
$|\!|\!|(x,y)|\!|\!|=\frac{2}{3}|x|+\frac{1}{3}|y|$. It is easy to
see that $|\!|\!|\cdot|\!|\!|$ is a tracial gauge norm on
$(\cc^2,\tau)$. Let $\M_1$ be the reduced free product of
$(\cc^2,\tau)$ with the free group factor $\L(\F_2)$. Then $\M_1$ is
a type ${\rm II}\sb 1$ factor with a faithful normal trace $\tau_1$
such that the restriction of $\tau_1$ to $\cc^2$ is $\tau$. Suppose
$|\!|\!|\cdot|\!|\!|_1$ is a unitarily invariant norm on $\M_1$ such
that the restriction of $|\!|\!|\cdot|\!|\!|_1$ to $\cc^2$ is
$|\!|\!|\cdot|\!|\!|$. Let $E=(1,0)$ and $F=(0,1)$ in $\cc^2$. Then
$\tau_1(E)=\tau(E)<\tau(F)=\tau_1(F)$. So there is a unitary
operator $U$ in $\M_1$ such that $UEU^*\leq F$.  By
Corollary~\ref{C:S<T},
$\frac{2}{3}=|\!|\!|E|\!|\!|=|\!|\!|E|\!|\!|_1=|\!|\!|UEU^*|\!|\!|_1\leq
|\!|\!|F|\!|\!|_1=|\!|\!|F|\!|\!|=\frac{1}{3}$. It is a
contradiction.}
\end{Example}

\begin{proof}[Proof of Theorem~{\bf E}] Let
 $|\!|\!|\cdot|\!|\!|_2$ be the tracial gauge norm on
 $\M$ corresponding to the symmetric gauge norm
 $|\!|\!|\cdot|\!|\!|_1^\#$
on $(L^\infty[0,1],\int_0^1dx)$  as in Theorem~{\bf D}.  By
Lemma~\ref{C: type I_n operators are dense}, to prove
$|\!|\!|\cdot|\!|\!|_2=|\!|\!|\cdot|\!|\!|^\#$ on $\M$, we need to
prove $|\!|\!|T|\!|\!|_2=|\!|\!|T|\!|\!|^\#$ for every positive
simple operator $T=a_1E_1+\cdots+a_nE_n$ in $\M$ such that
$\tau(E_1)=\cdots=\tau(E_n)=\frac{1}{n}$. We may assume that
$a_1\geq \cdots\geq a_n\geq 0$. Then
$\mu_s(T)=a_1\chi_{[0,\frac{1}{n})}(s)+\cdots+a_n\chi_{[\frac{n-1}{n},1]}(s)$.
By Lemma~\ref{L:Dual norm of simple operators},
\[|\!|\!|T|\!|\!|^\#=\sup\left\{\frac{1}{n}\sum_{k=1}^n a_kb_k:\,\,
X=b_1E_1+\cdots+b_nE_n \geq 0 \text{ and } |\!|\!|X|\!|\!|^\#\leq
1\right\}.\] By Lemma~\ref{L:HLP},
\[|\!|\!|T|\!|\!|^\#=\sup\left\{\frac{1}{n}\sum_{k=1}^na_kb_k:\,
X=b_1E_1+\cdots+b_nE_n\geq 0,\, b_1\geq \cdots\geq b_n\geq 0,\,
|\!|\!|X|\!|\!|\leq 1\right\}.\] By Theorem~{\bf D} and
Lemma~\ref{L:Dual norm of simple operators},
\begin{multline*}
|\!|\!|T|\!|\!|_2
=|\!|\!|\mu_s(T)|\!|\!|^\#=\sup\{\frac{1}{n}\sum_{k=1}^na_kb_k:\,
g(s)=b_1\chi_{[0,\frac{1}{n})}(s)+\cdots+b_n\chi_{[\frac{n-1}{n},1]}(s)\geq
0,\\ |\!|\!|g(s)|\!|\!|\leq 1\}.
\end{multline*}
 By Lemma~\ref{L:HLP},
 \begin{multline*}
|\!|\!|T|\!|\!|_2=|\!|\!|\mu_s(T)|\!|\!|^\#=\sup{
\{}\frac{1}{n}\sum_{k=1}^na_kb_k:\,
g(s)=b_1\chi_{[0,\frac{1}{n})}(s)+\cdots+b_n\chi_{[\frac{n-1}{n},1]}(s)\geq
0,\\
 b_1\geq \cdots b_n\geq 0,\, |\!|\!|g(s)|\!|\!|\leq 1{ \}}.
\end{multline*}
 Note
that if $b_1\geq \cdots\geq b_n\geq 0$, then
$\mu_s(b_1E_1+\cdots+b_nE_n)=b_1\chi_{[0,\frac{1}{n})}(s)+\cdots+b_n\chi_{[\frac{n-1}{n},1]}(s)$.
Since $|\!|\!|\cdot|\!|\!|$ is the tracial gauge norm on $(\M,\tau)$
corresponding to the symmetric gauge norm $|\!|\!|\cdot|\!|\!|_1$ on
$(L^\infty[0,1],\int_0^1dx)$ as in Theorem~{\bf D},
$|\!|\!|b_1E_1+\cdots+b_nE_n|\!|\!|\leq 1$ if and only if
$|\!|\!|b_1\chi_{[0,\frac{1}{n})}(s)+\cdots+b_n\chi_{[\frac{n-1}{n},1]}(s)|\!|\!|_1\leq
1$. Therefore, $|\!|\!|T|\!|\!|_2=|\!|\!|T|\!|\!|^\#$.
\end{proof}

\begin{Example}\label{E:dual of Lp space} \emph{If $p=1$, let $q=\infty$. If $1<p<\infty$, let $q=\frac{p}{p-1}$.
Then the $L^q$-norm  on $L^\infty[0,1]$  is the dual norm of the
$L^p$-norm on $L^\infty[0,1]$. By Theorem~{\bf E}, the $L^q$-norm on
a type ${\rm II}\sb 1$ factor $\M$ is the dual norm of the
$L^p$-norm on $\M$.}
\end{Example}

\section{Proof of Theorem F}

\begin{proof}[Proof of Theorem~{\bf F}]  Let $|\!|\!|\cdot|\!|\!|$ be a
tracial gauge norm on $\M$. By Lemma~\ref{L:T_f},
$|\!|\!|S|\!|\!|_f\leq |\!|\!|T|\!|\!|_f$ for every $f\in \F$. By
Theorem~{\bf A}, $|\!|\!|S|\!|\!|\leq |\!|\!|T|\!|\!|$.
\end{proof}

\begin{Corollary}\label{C:Ky Fan Theorem for 21factor}
 Let $\M$ be a type ${\rm II}\sb 1$ factor and $S, T\in \M$. If
$|\!|\!|S|\!|\!|_{(t)}\leq |\!|\!|T|\!|\!|_{(t)}$ for all Ky Fan
$t$-th norms, $0\leq t\leq 1$, then $|\!|\!|S|\!|\!|\leq
|\!|\!|T|\!|\!|$ for all unitarily invariant norms
$|\!|\!|\cdot|\!|\!|$ on $\M$.
\end{Corollary}

By Corollary~\ref{C:representation theorem for UIN on matrices}, we
obtain Ky Fan's Dominance Theorem~\cite{Fan}.\\

\noindent{\bf Ky Fan's Dominance Theorem.}\,\, \emph{If $S,T\in
M_n(\cc)$ and $|\!|\!|S|\!|\!|_{(k/n)}\leq |\!|\!|T|\!|\!|_{(k/n)}$,
i.e., $\sum_{i=1}^k s_i(S)\leq \sum_{i=1}^k s_i(T)$ for $1\leq k\leq
n$, then $|\!|\!|S|\!|\!|\leq |\!|\!|T|\!|\!|$ for all unitarily
invariant norms $|\!|\!|\cdot|\!|\!|$ on $M_n(\cc)$.}

\section{Extreme points of normalized unitarily invariant norms on
finite factors}

In this section, we assume that $\M$ is a finite factor with the
unique tracial state $\tau$.

\subsection{$\NN(\M)$}

 Let $\mathfrak{N}(\M)$ be the set of normalized unitarily invariant
norms on $\M$. It is easy to see that $\mathfrak{N}(\M)$ is a convex
set.  Let $\mathfrak{F}(\M)$ be the set of complex functions defined
on $\M$. Then $\mathfrak{F}(\M)$ is a locally convex space such that
a neighborhood of $f\in \mathfrak{F}(\M)$ is
\[N(f,T_1,\cdots,T_n,\epsilon)=\{g\in \mathfrak{F}(\M):\, |g(T_i)-
f(T_i)|<\epsilon\}.\] In this  topology, $f_\alpha\rightarrow f$
means $\lim_\alpha f_\alpha(T)=f(T)$ for every $T\in \M$. We call
this topology the \emph{pointwise weak} topology.

\begin{Lemma}\label{L:the set of normalized norms is compact} $\mathfrak{N}(\M)\subseteq \FF(\M)$ is a compact convex
subset in the pointwise weak topology.
\end{Lemma}
\begin{proof} It is clear that $\NN(\M)$ is a convex subset of
$\mathfrak{F}(\M)$. Suppose $|\!|\!|\cdot|\!|\!|_\alpha\in \FF(\M)$
and $f(T)=\lim_\alpha |\!|\!|T|\!|\!|_\alpha$ for every $T\in \M$.
It is easy to check that $f(T)$ defines a unitarily invariant
semi-norm on $\M$ such that $f(1)=1$. By Corollary~\ref{C:unitarily
seminorm is norm}, $f(T)$ is a norm and $f\in \NN(\M)$.
\end{proof}

 Let $\NN_e(\M)$ be the subset of extreme points of
$\NN(\M)$. By the Krein-Milman Theorem,  the closure of the convex
hull of $\NN_e(\M)$  is $\NN(\M)$  in the pointwise weak topology.
It is an interesting question of characterizing $\NN_e(\M)$. In the
following, we will provide some results on $\NN_e(\M)$.

\subsection{$\NN_e(M_n(\cc))$}
For $n\geq 2$, let $1\oplus s_2\oplus\cdots\oplus s_n$ be the matrix $\left(\begin{array}{cccc} 1&&&\\
&s_2&&\\
&&\ddots&\\
 &&&s_{n} \end{array}\right)\in M_n(\cc)$. Let $|\!|\!|\cdot|\!|\!|$ be
 a normalized
unitarily invariant norm on $M_n(\cc)$. For $0\leq
s_{n}\leq\cdots\leq s_2\leq 1$, define
\begin{equation}\label{Eq:norm function}
f(s_2,\cdots,s_{n})=f_{|\!|\!|\cdot|\!|\!|}(s_2,\cdots,s_{n})=|\!|\!|1\oplus
s_2\oplus\cdots\oplus s_{n}|\!|\!|.
\end{equation}
 In the following, let
$\Omega_{n-1}=\{(s_2,\cdots,s_{n}): 0\leq s_n\leq \cdots \leq
s_{2}\leq 1\}$.\\

By Lemma 3.2 of~\cite{G-K} and Corollary~\ref{C:unitarily seminorm
is norm}, we have the following lemma.
\begin{Lemma}\label{L:convex increasing multivariable
function} Let  $f(s_2,\cdots,s_{n})$ be a function defined on
$\Omega_{n-1}$. In order that $f(s_2,\cdots,
s_{n})=f_{|\!|\!|\cdot|\!|\!|}(s_2,\cdots,s_{n})$ for some
$|\!|\!|\cdot|\!|\!|\in \NN(M_n(\cc))$, it is necessary and
sufficient that the following conditions are satisfied:
\begin{enumerate}
\item $f(s_2,\cdots,s_{n})>0$ for all $(s_2,\cdots, s_{n})\in
\Omega_{n-1}$ and $f(1,\cdots, 1)=1$;
\item $f(s_2,\cdots,s_{n})$ is a convex function on $\Omega_{n-1}$;
\item  for $0\leq s_{n}\leq s_{n-1}\leq \cdots\leq s_1,\,
0\leq t_{n}\leq t_{n-1}\leq \cdots\leq t_1$, if $\sum_{i=1}^ks_i\leq
\sum_{i=1}^kt_i$ for $1\leq k\leq n$, then $s_1\cdot
f(\frac{s_2}{s_1},\cdots,\frac{s_{n}}{s_1})\leq t_1\cdot
f(\frac{t_2}{t_1},\cdots,\frac{t_{n}}{t_1})$.
\end{enumerate}
If $f(s_2,\cdots,s_{n})$ satisfies the above conditions, then $f$
satisfies
\[\frac{1+s_2+\cdots+s_{n}}{n}\leq f(s_2,\cdots,s_{n})\leq 1\]
for all $(s_2,\cdots,s_{n})\in \Omega_{n-1}$.
\end{Lemma}

Let $|\!|\!|\cdot|\!|\!|_1,\,
|\!|\!|\cdot|\!|\!|_2\in\NN(M_n(\cc))$. If
$|\!|\!|S|\!|\!|_1=|\!|\!|S|\!|\!|_2$ for all $S=1\oplus
s_2\oplus\cdots\oplus s_{n}$, $(s_2,\cdots,s_{n})\in \Omega_{n-1}$,
then $|\!|\!|T|\!|\!|_1=|\!|\!|T|\!|\!|_2$ for every matrix $T\in
M_n(\cc)$. This implies the following lemma.

\begin{Lemma}\label{L:norms are determined by the induced
functions}Let $|\!|\!|\cdot|\!|\!|_1,\,
|\!|\!|\cdot|\!|\!|_2\in\NN(M_n(\cc))$. Then
$|\!|\!|\cdot|\!|\!|_1=|\!|\!|\cdot|\!|\!|_2$ if and only if
$f_{|\!|\!|\cdot|\!|\!|_1}(s_2,\cdots,s_{n})=f_{|\!|\!|\cdot|\!|\!|_2}(s_2,\cdots,s_{n})$
for all $(s_2,\cdots,s_{n})\in \Omega_{n-1}$.
\end{Lemma}

Let $1\leq m\leq n$. Suppose $|\!|\!|\cdot|\!|\!|$ is a normalized
unitarily invariant norm on $M_m(\cc)$ and  $g(s_2,\cdots,
s_{m})=g_{|\!|\!|\cdot|\!|\!|}(s_2,\cdots,s_m)$ is the function on
$\Omega_{m-1}$ induced by $|\!|\!|\cdot|\!|\!|$ (see~(\ref{Eq:norm
function})). Define $f(s_2,\cdots,s_n)$ on $\Omega_{n-1}$ by
\[f(s_2,\cdots,s_n)=g(s_2,\cdots,s_{m}),\quad (s_2,\cdots,s_n)\in
\Omega_{n-1}.
\] It is easy to check that $f(s_2,\cdots,s_n)$ is a function on
$\Omega_{n-1}$ satisfying Lemma~\ref{L:convex increasing
multivariable function}. By Lemma~\ref{L:convex increasing
multivariable function} and Lemma~\ref{L:norms are determined by the
induced functions},  there is a unique normalized unitarily
invariant norm $|\!|\!|\cdot|\!|\!|_1\in \NN(\M_n(\cc))$ such that
$f(s_2,\cdots,s_n)=f_{|\!|\!|\cdot|\!|\!|_1}(s_2,\cdots,s_{n})=g(s_2,\cdots,s_{m})$
for all $(s_1,\cdots,s_{n})\in \Omega_{n-1}$ (This fact can also be
obtained by Corollary~\ref{C:representation theorem for UIN on
matrices} and Lemma~\ref{L:norms are determined by the induced
functions}).  $|\!|\!|\cdot|\!|\!|_1$ is called the \emph{induced
norm} of $|\!|\!|\cdot|\!|\!|$. \\

Conversely, suppose $|\!|\!|\cdot|\!|\!|_1$ is a normalized
unitarily invariant norm on $M_n(\cc)$ and
$f(s_2,\cdots,s_n)=f_{|\!|\!|\cdot|\!|\!|_1}(s_2,\cdots,s_n)$ is the
function on $\Omega_{n-1}$ induced by  $|\!|\!|\cdot|\!|\!|_1$. If
$f(s_2,\cdots,s_n)=g(s_2,\cdots,s_m)$ for all $(s_2,\cdots,s_n)\in
\Omega_{n-1}$, then $g(s_2,\cdots,s_m)$ satisfies
Lemma~\ref{L:convex increasing multivariable function}. Hence, there
is a unique normalized unitarily invariant norm
$|\!|\!|\cdot|\!|\!|$ on $M_m(\cc)$ such that $g(s_2,\cdots,s_m)$
$=$ $g_{|\!|\!|\cdot|\!|\!|}(s_2,\cdots,s_m)$ for all
$(s_2,\cdots,s_m)\in \Omega_{m-1}$. $|\!|\!|\cdot|\!|\!|$ is called
the \emph{reduced
norm} of $|\!|\!|\cdot|\!|\!|_1$. \\

\begin{Proposition}\label{P:ky fan norms are extreme points} For $1\leq k\leq n$, the Ky Fan $\frac{k}{n}$-th
norm $($see Example~\ref{E:Ky fan norms for MnC}$)$ on $M_n(\cc)$ is
an extreme point of $\NN(M_n(\cc))$.
\end{Proposition}
\begin{proof} Suppose $0<\alpha<1$ and
$|\!|\!|\cdot|\!|\!|_1,\, |\!|\!|\cdot|\!|\!|_2\in \NN(M_n(\cc))$
satisfy
$|\!|\!|\cdot|\!|\!|_{\left(\frac{k}{n}\right)}=\alpha|\!|\!|\cdot|\!|\!|_1+(1-\alpha)|\!|\!|\cdot|\!|\!|_2$.
Let
$f(s_2,\cdots,s_n)=f_{|\!|\!|\cdot|\!|\!|_{\left(\frac{k}{n}\right)}}(s_2,\cdots,s_n)$,
$f_1(s_2,\cdots, s_n)=f_{|\!|\!|\cdot|\!|\!|_1}(s_2,\cdots, s_n)$
and $f_2(s_2,\cdots,s_n)=f_{|\!|\!|\cdot|\!|\!|_2}(s_2,\cdots,s_n)$
for $(s_2,\cdots, s_n)\in\Omega_{n-1}$. Then
$f(s_2,\cdots,s_{n-1})=\alpha
f_1(s_2,\cdots,s_{n-1})+(1-\alpha)f_2(s_2,\cdots,s_{n-1})$.\\

Since $f(s_2,\cdots,s_n)=\frac{1+s_2+\cdots+s_k}{k}$,
$\frac{\partial f}{\partial s_{k+1}}=\cdots=\frac{\partial
f}{\partial s_{n}}=0$. Since $f_1(s_2,\cdots,s_n)$,
$f_2(s_2,\cdots,s_n)$ are convex functions on $\Omega_{n-1}$,
$\frac{\partial f_i}{\partial s_{j}}\geq 0$ for $i=1,2$ and $k+1\leq
j\leq n$. Since $f=\alpha f_1+(1-\alpha)f_2$, $\frac{\partial
f_i}{\partial s_{j}}=0$ for  $i=1,2$ and $k+1\leq j\leq n$. This
implies that $f_i(s_2,\cdots,s_n)=g_i(s_2,\cdots,s_k)$ for all
$(s_2,\cdots,s_n)\in \Omega_{n-1}$ and $i=1,2$.\\

By the discussions above the proposition, there are normalized
unitarily invariant norms $|\!|\!|\cdot|\!|\!|_1$,
$|\!|\!|\cdot|\!|\!|_2$ on $M_k(\cc)$ such that
$g_i(s_2,\cdots,s_k)=(g_i)_{|\!|\!|\cdot|\!|\!|_i}(s_2,\cdots,s_k)$
for all $(s_2,\cdots,s_k)\in \Omega_{k-1}$ and $i=1,2$. By
Lemma~\ref{L:convex increasing multivariable function},
\[g_i(s_2,\cdots,s_k)\geq \frac{1+s_2+\cdots+s_k}{k}
\]for all $(s_2,\cdots,s_k)\in \Omega_{k-1}$ and $i=1,2$. Since
$f=\alpha f_1+(1-\alpha)f_2$, $\frac{1+s_2+\cdots+s_k}{k}=\alpha
g_1(s_2,\cdots,s_k)+(1-\alpha)g_2(s_2,\cdots,s_k)$. This implies
that
$g_1(s_2,\cdots,s_k)=g_2(s_2,\cdots,s_k)=\frac{1+s_2+\cdots+s_k}{k}$.
So $f=f_1=f_2$.
\end{proof}

The proof of the following proposition is similar to that of
Proposition~\ref{P:ky fan norms are extreme points}.
\begin{Proposition}\label{P:induced extreme points} Let $1\leq m\leq n$ and let $|\!|\!|\cdot|\!|\!|$ be
a normalized unitarily invariant norm on $M_m(\cc)$. If
$|\!|\!|\cdot|\!|\!|$ is an extreme point of $\NN(M_m(\cc))$, then
the induced norm $|\!|\!|\cdot|\!|\!|_1$ on $M_n(\cc)$ is also an
extreme point of $\NN(M_n(\cc))$.
\end{Proposition}

 \noindent{\bf Question.} \emph{For $n\geq 3$, find all  extreme
points of  $\NN(\M_n(\cc))$.}

\subsection{$\NN_e(M_2(\cc))$}

 In this subsection, we will prove
Theorem~{\bf J}. We need the following auxiliary results. The
following lemma is a corollary of Lemma~\ref{L:convex increasing
multivariable function} in the case $n=2$.
\begin{Lemma}\label{L:increasing convex function} Let  $f(s)$ be a
function on $[0,1]$. If there is a normalized unitarily invariant
norm $|\!|\!|\cdot|\!|\!|$ on $M_2(\cc)$ such that
$f(s)=f_{|\!|\!|\cdot|\!|\!|}(s)=|\!|\!|1\oplus s|\!|\!|$, then
$f(s)$
 is an increasing convex function
on $[0,1]$ satisfying
\[\frac{1+s}{2}\leq f(s)\leq 1,\quad\forall s\in [0,1].
\]
\end{Lemma}

\begin{Corollary}\label{L:derivative of convex function}
 For $0\leq a\leq b\leq 1$, we have \[0\leq f'(a-)\leq f'(a+)
\leq f'(b-)\leq f'(b+)\leq f'(1-)\leq \frac{1}{2}.\]
\end{Corollary}
\begin{proof} Since $f(s)$ is an increasing convex function, $0\leq f'(a-)\leq f'(a+)
\leq f'(b-)\leq f'(b+)\leq f'(1-)$. By Lemma~\ref{L:increasing
convex function},
\[f'(1-)=\lim_{h\rightarrow 0+}\frac{f(1)-f(1-h)}{h}\leq
\lim_{h\rightarrow 0+}\frac{1-(2-h)/2}{h}=\frac{1}{2}.
\]
\end{proof}

For $\frac{1}{2}\leq t\leq 1$, define
$|\!|\!|\cdot|\!|\!|_{<t>}=\max\{t\|T\|, \|T\|_1\}$.
\begin{Lemma}\label{L:Extreme points in NN(M2cc)} For $1/2\leq t\leq 1$, $|\!|\!|\cdot|\!|\!|_{<t>}$ is an
extreme point of $\NN(M_2(\cc))$.
\end{Lemma}
\begin{proof} Suppose $0<\alpha<1$ and $|\!|\!|\cdot|\!|\!|_1,\, |\!|\!|\cdot|\!|\!|_2\in
\NN(M_2(\cc))$ such that
$|\!|\!|\cdot|\!|\!|_{<t>}=\alpha|\!|\!|\cdot|\!|\!|_1+(1-\alpha)|\!|\!|\cdot|\!|\!|_2$.
Let $f(s)=f_{|\!|\!|\cdot|\!|\!|_{<t>}}(s)$,
$f_1(s)=f_{|\!|\!|\cdot|\!|\!|_1}(s)$ and
$f_2(s)=f_{|\!|\!|\cdot|\!|\!|_2}(s)$. Then $f(s)=\alpha
f_1(s)+(1-\alpha)f_2(s)$. Note that
\[f(s)=\begin{cases} t&0\leq
s\leq
\frac{2t-1}{2};\\
\frac{s+1}{2}&\frac{2t-1}{2}\leq s\leq 1.
\end{cases}
\]
Hence, $f'(s)=0$ if $0\leq s<\frac{2t-1}{2}$ and $f'(s)=\frac{1}{2}$
if $\frac{2t-1}{2}< s\leq 1$. By Corollary~\ref{L:derivative of
convex function}, $f'_1(s)=f'_2(s)=0$ if $0\leq s<\frac{2t-1}{2}$
and $f'_1(s)=f'_2(s)=\frac{1}{2}$ if $\frac{2t-1}{2}< s\leq 1$.
Since $f(s), f_1(s), f_2(s)$ are convex functions and hence
continuous and $f(1)=f_1(1)=f_2(1)=1$, $f(s)=f_1(s)=f_2(s)$ for all
$0\leq s\leq 1$. This implies that
$|\!|\!|\cdot|\!|\!|_{<t>}=|\!|\!|\cdot|\!|\!|_1=|\!|\!|\cdot|\!|\!|_2$.
\end{proof}

\begin{Lemma}\label{L:map from [0,1] to extreme points} The mapping: $t\rightarrow |\!|\!|\cdot|\!|\!|_{<t>}$ is
continuous with respect to the usual topology on $[1/2,1]$ and the
pointwise weak topology on $\NN(M_2(\cc))$. In particular,
$\{|\!|\!|\cdot|\!|\!|_{<t>}:\, 1/2\leq t\leq 1\}$ is compact in the
pointwise weak topology.
\end{Lemma}
\begin{proof} For every $0\leq s\leq 1$, $|\!|\!|1\oplus s|\!|\!|_{<t>}=\max\{t, \frac{1+s}{2}\}$ is a
continuous function on $[0,1]$. Hence, The mapping: $t\rightarrow
|\!|\!|\cdot|\!|\!|_{<t>}$ is continuous with respect to the usual
topology on $[1/2,1]$ and the pointwise weak topology on
$\NN(M_2(\cc))$.
\end{proof}

\begin{Lemma}\label{L:Choque theorem} The set
\[\S=\{|\!|\!|\cdot|\!|\!|:\, |\!|\!|\cdot|\!|\!|=\int_{1/2}^1 |\!|\!|\cdot|\!|\!|_{<t>}d\mu(t),\,
\mu \, \text{ is a regular Borel probability measure on $[1/2,1]$}\}
\] is a convex compact subset of $\NN(M_2(\cc))$ in the pointwise weak topology.
\end{Lemma}
\begin{proof} Suppose $\{|\!|\!|\cdot|\!|\!|_\alpha\}$ is a net in $\S$ such
that $|\!|\!|\cdot|\!|\!|_\alpha\rightarrow |\!|\!|\cdot|\!|\!|\in
\NN(M_2(\cc))$ in the pointwise weak topology. Let $\mu_\alpha$ be
the regular Borel probability measure on $[1/2,1]$ corresponding to
$|\!|\!|\cdot|\!|\!|_\alpha$. Then there is a subnet of $\mu_\alpha$
that converges weakly to a regular Borel probability measure $\mu$
on $[1/2,1]$, i.e., for every continuous function $\phi(t)$ on
$[1/2,1]$,
\[\lim_\alpha \int_{1/2}^1 \phi(t) d\mu_{\alpha_\beta}(t)=\int_{1/2}^1
\phi(t)d\mu(t).
\]
In particular, for every $T\in M_2(\cc)$, we have
\[|\!|\!|T|\!|\!|=\lim_{\alpha_\beta}|\!|\!|T|\!|\!|_{\alpha_\beta}=\lim_{\alpha_\beta}
\int_{1/2}^1|\!|\!|T|\!|\!|_{<t>}d\mu_{\alpha_\beta}(t)=\int_{1/2}^1|\!|\!|T|\!|\!|_{<t>}d\mu(t).\]
Hence $|\!|\!|\cdot|\!|\!|\in \S$.
\end{proof}

\begin{Lemma}\label{L:one-to-one correspondence between functions and norms}
Let $f(s)$ be a convex, increasing function on $[0,1]$ such that
\[\frac{1+s}{2}\leq f(s)\leq 1,\quad,\forall s\in [0,1].
\] Then there is an element  $|\!|\!|\cdot|\!|\!|\in \S$
such that $f(s)=|\!|\!|1\oplus s|\!|\!|.$
\end{Lemma}
\begin{proof} We can approximate $f$ uniformly by  piecewise linear functions
satisfying the conditions of the lemma. By Lemma~\ref{L:Choque
theorem}, we may assume that $f(s)$ is a piecewise linear function.
Furthermore, we may assume that $0=a_0<a_1<a_2<\cdots<a_n=1$ and
$f(s)$ is linear on $[a_i,a_{i+1}]$ for $0\leq i\leq n-1$. Let
$f'(s)=\alpha_i/2$ on $[a_i,a_{i+1}]$. By
Corollary~\ref{L:derivative of convex function}, $0=\alpha_0\leq
\alpha_1\leq\cdots\leq \alpha_{n-1}\leq 1$. Let
$g(s)=(1-\alpha_{n-1})\|1\oplus
s\|+(\alpha_{n-1}-\alpha_{n-2})|\!|\!|1\oplus
s|\!|\!|_{<\alpha_{n-1}>}+\cdots +(\alpha_1-\alpha_0)|\!|\!|1\oplus
s|\!|\!|_{<\alpha_1>}+\alpha_0\|1\oplus s\|_1$. Then $g(1)=f(1)=1$
and $g'(s)=\alpha_i/2$ on $[a_i,a_{i+1}]$. So $g'(s)=f'(s)$ except
$s=\alpha_i$ for $1\leq i\leq n$. Hence $f(s)=g(s)$ for all $0\leq
s\leq 1$.
\end{proof}

\begin{proof}[Proof of Theorem~{\bf J}] By Lemma~\ref{L:Extreme points in
NN(M2cc)}, $\{|\!|\!|\cdot|\!|\!|_{<t>}: 1/2\leq t\leq 1\}$ are
extreme points of $\NN(\M)$. By Lemma~\ref{L:Choque theorem},
Lemma~\ref{L:one-to-one correspondence between functions and norms}
and Lemma~\ref{L:norms are determined by the induced functions}, the
closure of the convex hull of $\{|\!|\!|\cdot|\!|\!|_{<t>}: 1/2\leq
t\leq 1\}$ in the pointwise weak topology is $\NN(\M_2(\cc))$. By
Lemma~\ref{L:Extreme points in NN(M2cc)}, \ref{L:map from [0,1] to
extreme points} and Theorem 1.4.5 of~\cite{K-R},
$\NN_e(M_2(\cc))=\{|\!|\!|\cdot|\!|\!|_{<t>}: 1/2\leq t\leq 1\}$.
\end{proof}

\begin{Corollary}\label{C:increasing convex function 2} Let  $f(s)$ be a
function on $[0,1]$. Then the following conditions are equivalent:
\begin{enumerate}
\item $f(s)=f_{|\!|\!|\cdot|\!|\!|}(s)=|\!|\!|1\oplus s|\!|\!|$ for
some normalized unitarily invariant norm $|\!|\!|\cdot|\!|\!|$ on
$M_2(\cc)$;
\item $f(s)$  is an increasing convex function
on $[0,1]$ such that  $\frac{1+s}{2}\leq f(s)\leq 1$ for all $s\in
[0,1]$;
\item $f(s)$  is an increasing convex function
on $[0,1]$ such that $f(1)=1$ and $f'(1-)\leq \frac{1}{2}$.
\end{enumerate}
\end{Corollary}

In the following, we will show how to write the $L^p$-norms on
$M_2(\cc)$ in terms of  extreme points of $\NN(M_2(\cc))$. Recall
that for $1\leq p<\infty$, the $L^p$-norm of $1\oplus s$ is
\[\|1\oplus s\|_p=\left(\frac{1+s^p}{2}\right)^{1/p}.
\]  Let $f_p(s)=f_{\|\cdot\|_p}(s)=\left(\frac{1+s^p}{2}\right)^{1/p}$, $0\leq s\leq
1.$ Then $f_p(1)=1$ and
\[f_p'(s)=\frac{s^{p-1}}{2}\left(\frac{1+s^p}{2}\right)^{1/p-1},
\] $f_p'(0)=0$, $f_p'(1)=\frac{1}{2}$.

\begin{Lemma} For $1<p<\infty$ and $0\leq s\leq 1$,
\[f_p(s)=\int_{1/2}^1 |\!|\!|1\oplus s|\!|\!|_{<t>} 4f_p^{''}(2t-1)dt.
\]
\end{Lemma}
\begin{proof}\[\int_{1/2}^1 |\!|\!|1\oplus s|\!|\!|_{<t>} 4f_p^{''}(2t-1)dt=
\int_0^1 |\!|\!|1\oplus s|\!|\!|_{<\frac{x+1}{2}>} 2f_p^{''}(x)dx
=\int_0^s \frac{1+s}{2}2f_p^{''}(x)dx+\]\[\int_s^1
\frac{1+x}{2}2f_p^{''}(x)dx
=(1+s)f'(s)-(1+s)f'(0)+2f'(1)-(1+s)f'(s)-\int_s^1
f_p'(x)dx\]\[=1-f_p(1)+f_p(s)=f_p(s).\]
\end{proof}

\begin{Corollary} For $1<p<\infty$ and $T\in M_2(\cc)$,
\[\|T\|_p=\int_{1/2}^1 |\!|\!|T|\!|\!|_{<t>} 4f_p''(2t-1)dt.
\]
\end{Corollary}

\subsection{Proof of Theorem K}

\begin{Lemma}\label{L:Restriction of an extreme norm} Let $\M$ be a type ${\rm II}\sb 1$ factor and let $|\!|\!|\cdot|\!|\!|$ be a normalized unitarily invariant norm on
$\M$. Suppose $\N_1\subset\N_2\subset\cdots$ are a sequence of type
$I_{n_r}$ subfactors of $\M$ such that  $\N_r\cong M_{n_r}(\cc)$ and
$\lim_{r\rightarrow \infty} n_r=\infty$. If the restriction of
$|\!|\!|\cdot|\!|\!|$ to $\N_r$ is an extreme point of $\NN(\N_r)$
for all $r=1,2,\cdots,$ then $|\!|\!|\cdot|\!|\!|$ is an extreme
point of $\NN(\M)$.
\end{Lemma}
\begin{proof}Suppose $0<\alpha<1$ and $|\!|\!|\cdot|\!|\!|_1, |\!|\!|\cdot|\!|\!|_2\in
\NN(\M)$ such that
$|\!|\!|\cdot|\!|\!|=\alpha|\!|\!|\cdot|\!|\!|_1+(1-\alpha)|\!|\!|\cdot|\!|\!|_2$
on $\M$. Then for every $r=1,2,\cdots$,
$|\!|\!|\cdot|\!|\!|=\alpha|\!|\!|\cdot|\!|\!|_1+(1-\alpha)|\!|\!|\cdot|\!|\!|_2$
on $\N_r$. By the assumption of the lemma,
$|\!|\!|\cdot|\!|\!|=|\!|\!|\cdot|\!|\!|_1=|\!|\!|\cdot|\!|\!|_2$ on
$\N_r$. By Corollary~\ref{C: type I_n subfactor are dense},
$|\!|\!|\cdot|\!|\!|=|\!|\!|\cdot|\!|\!|_1=|\!|\!|\cdot|\!|\!|_2$ on
$\M$. So $|\!|\!|\cdot|\!|\!|$ is an extreme point of $\NN(\M)$.
\end{proof}

\begin{proof}[Proof of Theorem~{\bf K}] By the assumption of the
theorem, $t=\frac{k}{n}$ is a rational number. Then we can construct
a sequence of type ${\rm I}_{rn}$ subfactor $\M_{rn}$ of $\M$ such
that $\M_{n}\subseteq \M_{2n}\subseteq \cdots$. Then the restriction
of $\!|\!|\cdot|\!|\!|_{(t)}$ on $\M_{rn}$ is
$\!|\!|\cdot|\!|\!|_{\left(\frac{rk}{rn}\right)}$. By
Proposition~\ref{P:ky fan norms are extreme points}, the restriction
of $\!|\!|\cdot|\!|\!|_{(t)}$ on $\M_{rn}$ is an extreme point of
$\NN(M_{rn}(\cc))$. By Lemma~\ref{L:Restriction of an extreme norm},
$\!|\!|\cdot|\!|\!|_{(t)}$ is an extreme point of $\NN(\M)$.\\
\end{proof}

\begin{Remark} \emph{Here we point out
other interesting examples of extreme points of $\NN(\M)$. For
$0\leq t\leq 1$, recall that $\!|\!|\cdot|\!|\!|_{(t)}$ is the
$t$-th Ky-fan norm on $\M$. For any non-negative function $c(t)$ on
$[0,1]$ such that $\|c(t)\|_\infty=1$ and $T\in\M$,  define
\[\!|\!|T|\!|\!|_{[c(t)]}=\|c(t)\!|\!|T|\!|\!|_{(t)}\|_\infty.
\]
Then it is easy to see that $\!|\!|\cdot|\!|\!|_{[c(t)]}$ is a
normalized unitarily invariant norm on $\M$. It can be proved that
if $c(t)$ is a simple function or if $tc(t)$ is a simple function,
then $\!|\!|\cdot|\!|\!|_{[c(t)]}$ is an extreme point of
$\NN(\M)$.}
\end{Remark}

\section{Proof of Theorem G}

In this section, we assume that $\M$ is a type ${\rm II}\sb 1$
factor with the unique tracial state $\tau$ and
$|\!|\!|\cdot|\!|\!|$ is a unitarily invariant norm on $\M$. For two
projections $E,F$ in $\M$, $\tau(E)\leq \tau(F)$ if and only if
there is a unitary operator $U\in \M$ such that $UEU^*\leq F$. By
Corollary~\ref{C:S<T}, if $\tau(E)\leq \tau(F)$,
$|\!|\!|E|\!|\!|\leq |\!|\!|F|\!|\!|$. So we can define
\[r(|\!|\!|\cdot|\!|\!|)=\lim_{\tau(E)\rightarrow 0+}|\!|\!|E|\!|\!|.
\]
\begin{Definition}\label{D:singular and continuous norms}
\emph{ A unitarily invariant norm $|\!|\!|\cdot|\!|\!|$ on $\M$ is
\emph{singular} if $r(|\!|\!|\cdot|\!|\!|)>0$ and \emph{continuous}
if $r(|\!|\!|\cdot|\!|\!|)=0$.}
\end{Definition}

\begin{Example}\label{E: noncommutative Lp norm}
\emph{The operator norm is singular since
$r(\|\cdot\|)=\lim_{\tau(E)\rightarrow 0+}\|E\|=1$. If $0<t\leq 1$,
the Ky Fan $t$-th  norm $|\!|\!|\cdot|\!|\!|_{(t)}$ is continuous
since
$r(|\!|\!|\cdot|\!|\!|_{(1)})=r(\|\cdot\|_1)=\lim_{\tau(E)\rightarrow
0+}\tau(E)=0$ and $r(|\!|\!|\cdot|\!|\!|_{(t)})\leq \frac{1}{t}\cdot
r(|\!|\!|\cdot|\!|\!|_{(1)})=0$. If $1\leq p<\infty$, it is easy to
see that the $L^p$-norm  on $\M$ is also continuous.}
\end{Example}
\begin{Lemma}\label{L:singular norms} If $|\!|\!|\cdot|\!|\!|$ is singular, then $|\!|\!|\cdot|\!|\!|$ is
equivalent to the operator norm $\|\cdot\|$. Indeed, for every
$T\in\M$, we have
\[r(|\!|\!|\cdot|\!|\!|) \|T\|\leq |\!|\!|T|\!|\!|\leq |\!|\!|1|\!|\!|\cdot\|T\|.
\]
\end{Lemma}
\begin{proof} By Lemma~\ref{L:submultiplicative}, $|\!|\!|T|\!|\!|\leq
|\!|\!|1|\!|\!|\cdot\|T\|$.  We need to prove
$r(|\!|\!|\cdot|\!|\!|) \|T\|\leq |\!|\!|T|\!|\!|$. We may assume
that $T>0$. For any $\epsilon>0$, let
$E=\chi_{[\|T\|-\epsilon,\|T\|]}(T)>0$. Then $T\geq
(\|T\|-\epsilon)E$. By Corollary~\ref{C:S<T} and
Lemma~\ref{L:submultiplicative}, $|\!|\!|T|\!|\!|\geq
|\!|\!|(\|T\|-\epsilon)E|\!|\!|\geq (\|T\|-\epsilon)\cdot
|\!|\!|E|\!|\!|\geq (\|T\|-\epsilon)r(|\!|\!|\cdot|\!|\!|)$. Since
$\epsilon>0$ is arbitrary, $r(|\!|\!|\cdot|\!|\!|) \|T\|\leq
|\!|\!|T|\!|\!|$.
\end{proof}

Recall that a neighborhood $N(\epsilon,\delta)$ of $0\in\M$ in the
measure topology (see~\cite{Ne}) is
\[N(\epsilon,\delta)=\{T\in \M, \,\text{there is a projection $E\in\M$ such that
$\tau(E)<\delta$ and $\|TE^\perp\|<\epsilon$}\}
\]

\begin{proof}[Proof of Theorem {\bf G}] By Lemma~\ref{L:singular norms}, if
$|\!|\!|\cdot|\!|\!|$ is singular, then  $\mathcal{T}$ is the
operator topology on $\M_{1,\, \|\cdot\|}$. Suppose
$|\!|\!|\cdot|\!|\!|$ is continuous. For $\epsilon,\delta>0$ and
$T\in \M$ such that $\|T\|\leq 1$ and
$|\!|\!|T|\!|\!|<\epsilon\delta$, by Corollary~\ref{C:unitarily
seminorm is norm}, $\tau(\chi_{[\epsilon,1]}(|T|)\leq
\frac{\|T\|_1}{\epsilon}\leq \frac{|\!|\!|T|\!|\!|}{\epsilon}<
\delta$ and $\|T\cdot\chi_{[0,\epsilon)}(|T|)\|<\epsilon$. This
implies that $\{T\in\M_{1,\, \|\cdot\|}:\,
|\!|\!|T|\!|\!|<\epsilon\delta\}\subseteq N(\epsilon,\delta)$.
Conversely, let $\omega>0$. Since $r(|\!|\!|\cdot|\!|\!|)=0$, there
is an $\epsilon$, $0<\epsilon<\omega/2$,  such that if
$\tau(E)<\epsilon$ then $|\!|\!|E|\!|\!|<\omega/2$. For every $T\in
N(\epsilon, \omega/2)$ and $\|T\|\leq 1$, choose $E\in \M$ such that
$\tau(E)<\epsilon$ and $\|TE^\perp\|<\omega/2$. By
Proposition~\ref{P:unitarily invariant norms} and
Corollary~\ref{C:compare gauge norm with operator norm},
$|\!|\!|T|\!|\!|\leq
|\!|\!|TE|\!|\!|+|\!|\!|TE^\perp|\!|\!|<\|T\|\cdot
|\!|\!|E|\!|\!|+\|TE^\perp\|<\omega/2+\omega/2=\omega$. Hence
$\{T\in N(\epsilon,\omega/2):\, \|T\|\leq 1\}\subseteq \{T\in\M:\,
|\!|\!|T|\!|\!|<\omega\}$.
\end{proof}

\begin{Corollary}\label{C:topology induced by Lp norms are same}
Topologies induced by the $L^p$-norms, $1\leq p<\infty$, on the unit
ball of a type ${\rm II}\sb 1$ factor are same.
\end{Corollary}

\section{Completion of type ${\rm II}\sb 1$ factors with respect to unitarily invariant norms}
In this section, we assume that $\M$ is a type ${\rm II}\sb 1$
factor with the unique tracial state $\tau$ and
$|\!|\!|\cdot|\!|\!|$ is a unitarily invariant norm on $\M$. The
completion of $\M$ with respect to $|\!|\!|\cdot|\!|\!|$ is denoted
by $\overline{{\M}_{|\!|\!|\cdot|\!|\!|}}$.  We will use the
traditional notation $L^p(\M,\tau)$ to denote the completion of $\M$
with respect to the $L^p$-norm defined as in Example~\ref{E:
noncommutative Lp norm}. Note that $L^\infty(\M,\tau)=\M$. Let
$\widetilde{\M}$ be the completion of $\M$ in the measure-topology
in the sense of~\cite{Ne}.

\subsection{Embedding of $\overline{{\M}_{|\!|\!|\cdot|\!|\!|}}$
into $\widetilde{\M}$}
\begin{Lemma}\label{L:epsilon-delta for one operator} Let $|\!|\!|\cdot|\!|\!|$ be a continuous unitarily invariant
norm on $\M$ and $T\in \M$. For every $\epsilon>0$, there is a
$\delta>0$ such that if $\tau(E)<\delta$, then
$|\!|\!|TE|\!|\!|<\epsilon.$
\end{Lemma}
\begin{proof}  Since $|\!|\!|\cdot|\!|\!|$ is continuous,
$\lim_{\tau(E)\rightarrow 0}|\!|\!|E|\!|\!|=0$. Hence, for every
$\epsilon>0$, there is a $\delta>0$ such that if $\tau(E)<\delta$,
then $|\!|\!|E|\!|\!|<\frac{\epsilon}{1+\|T\|}$.  By
Proposition~\ref{P:unitarily invariant norms}, $|\!|\!|TE|\!|\!|\leq
\|T\|\cdot |\!|\!|E|\!|\!|<\epsilon$.
\end{proof}

\begin{Lemma}\label{L:epsilon-delta for a sequence operators}Let $|\!|\!|\cdot|\!|\!|$ be a continuous unitarily invariant
norm on $\M$ and $\{T_n\}$ in  $\M$ be a Cauchy sequence with
respect to $|\!|\!|\cdot|\!|\!|$. For every $\epsilon>0$, there is a
$\delta>0$ such that if $\tau(E)<\delta$, then
$|\!|\!|T_nE|\!|\!|<\epsilon$ for all $n$.
\end{Lemma}
\begin{proof}  Since $\{T_n\}$ is a Cauchy sequence with respect to
$|\!|\!|\cdot|\!|\!|$, there is an $N$ such that for all $n\geq N$,
$|\!|\!|T_n-T_N|\!|\!|<\epsilon/2$. By Lemma~\ref{L:epsilon-delta
for one operator}, there is a $\delta_1$ such that if
$\tau(E)<\delta_1$ then $|\!|\!|T_NE|\!|\!|<\epsilon/2$. By
Proposition~\ref{P:unitarily invariant norms}, for $n\geq N$,
$|\!|\!|T_nE|\!|\!|\leq
|\!|\!|(T_n-T_N)E|\!|\!|+|\!|\!|T_NE|\!|\!|<|\!|\!|(T_n-T_N)|\!|\!|\cdot
\|E\|+\epsilon/2<\epsilon$. A simple argument shows that we can
choose $0<\delta<\delta_1$ such that if $\tau(E)<\delta$ then
$|\!|\!|T_nE|\!|\!|<\epsilon$ for all $n$.
\end{proof}

The following proposition generalizes Theorem 5 of~\cite{Ne}.
\begin{Proposition}\label{P:embedding into measure topology} Let
$\M$ be a type ${\rm II}\sb 1$ factor and let
$|\!|\!|\cdot|\!|\!|$ be a unitarily invariant norm on $\M$. There
is an injective map from $\overline{{\M}_{|\!|\!|\cdot|\!|\!|}}$
to $\widetilde{\M}$ that extends the identity map from $\M$ to
$\M$.
\end{Proposition}
\begin{proof} If $|\!|\!|\cdot|\!|\!|$ is singular, by Lemma~\ref{L:singular
norms}, $\overline{{\M}_{|\!|\!|\cdot|\!|\!|}}=\M$. So we will
assume that  $|\!|\!|\cdot|\!|\!|$ is continuous. If $\{T_n\}$ in
$\M$ is a Cauchy sequence with respect to $|\!|\!|\cdot|\!|\!|$,
then $\{T_n\}$ is a Cauchy sequence in the $L^1$-norm by
Corollary~\ref{C:unitarily seminorm is norm}. For every $\delta>0$
and $T\in \M$, $\tau(\chi_{(\delta,\infty)}(|T|)\leq
\frac{\tau(|T|)}{\delta}$. Hence, if $\{T_n\}$ is a Cauchy sequence
in $\M$ in the $L^1$-norm, then $\{T_n\}$ is a Cauchy sequence in
the measure topology. So there is a natural map $\Phi$ from
$\overline{{\M}_{|\!|\!|\cdot|\!|\!|}}$ to $\widetilde{\M}$ that
extends the identity map from $\M$ to $\M$. To prove that $\Phi$ is
injective, we need to prove that if $\{T_n\}$ in $\M$ is a Cauchy
sequence with respect to $|\!|\!|\cdot|\!|\!|$ and $T_n\rightarrow
0$ in the measure topology, then
$\lim_{n\rightarrow\infty}|\!|\!|T_n|\!|\!|=0$. Let $\epsilon>0$. By
Lemma~\ref{L:epsilon-delta for a sequence operators}, there is a
$\delta>0$ such that if $\tau(E)<\delta$ then
$|\!|\!|T_nE|\!|\!|<\epsilon/2$ for all $n$.  Since $T_n\rightarrow
0$ in the measure topology, there are $N$ and $\delta_1$,
$0<\delta_1<\delta$, such that for all $n\geq N$, there is a
projection $E_n$ such that $\tau(E_n)<\delta_1$ and
$\|T_nE_n^\perp\|<\epsilon/2$. By Corollary~\ref{C:unitarily
seminorm is norm}, $|\!|\!|T_n|\!|\!|\leq
|\!|\!|T_nE_n^\perp|\!|\!|+|\!|\!|T_nE_n|\!|\!|<
\|T_nE_n^\perp\|+\epsilon/2<\epsilon$. This proves that
$\lim_{n\rightarrow\infty}|\!|\!|T_n|\!|\!|=0$ and hence $\Phi$ is
an injective map from $\overline{{\M}_{|\!|\!|\cdot|\!|\!|}}$ to
$\widetilde{\M}$ that extends the identity map from $\M$ to $\M$.
\end{proof}

By the proof of Proposition~\ref{P:embedding into measure topology},
we have the following
\begin{Corollary}\label{L:The completion of M is in L1} There is an
injective map from $\overline{{\M}_{|\!|\!|\cdot|\!|\!|}}$ to
$L^1(\M,\tau)$ that extends the identity map from $\M$ to $\M$.
\end{Corollary}

By Proposition~\ref{P:embedding into measure topology}, we will
consider $\overline{{\M}_{|\!|\!|\cdot|\!|\!|}}$ as a subset of
$\widetilde{\M}$. The following corollary is very useful.
\begin{Corollary}\label{C:embedding into measure topology} Let
$\M$ be a type ${\rm II}\sb 1$ factor and let
$|\!|\!|\cdot|\!|\!|$ be a unitarily invariant norm on $\M$. If
$\{T_n\}\subset\M$ is a Cauchy sequence with respect to
$|\!|\!|\cdot|\!|\!|$ and $\lim_{n\rightarrow\infty}T_n=T$ in the
measure topology, then $T\in
\overline{{\M}_{|\!|\!|\cdot|\!|\!|}}$ and
$\lim_{n\rightarrow\infty}T_n=T$ in the topology induced by
$|\!|\!|\cdot|\!|\!|$.

\end{Corollary}

\begin{Corollary}\label{C:extension of norm on the completion algebra}
 $\overline{{\M}_{|\!|\!|\cdot|\!|\!|}}$ is a
 linear subspace of $\widetilde{\M}$ satisfying the following conditions:
\begin{enumerate}
\item if $T\in \overline{{\M}_{|\!|\!|\cdot|\!|\!|}}$, then
$T^*\in \overline{{\M}_{|\!|\!|\cdot|\!|\!|}}$;
\item $T\in \overline{{\M}_{|\!|\!|\cdot|\!|\!|}}$ if and only if
$|T|\in \overline{{\M}_{|\!|\!|\cdot|\!|\!|}}$;
\item if $T\in \overline{{\M}_{|\!|\!|\cdot|\!|\!|}}$ and $A, B\in
\M$,
 then $ATB\in \overline{{\M}_{|\!|\!|\cdot|\!|\!|}}$ and
$|\!|\!|ATB|\!|\!|\leq \|A\|\cdot |\!|\!|T|\!|\!|\cdot\|B\|$.
\end{enumerate}
In particular, $|\!|\!|\cdot|\!|\!|$ can be extended to  a unitarily
invariant norm, also denoted by $|\!|\!|\cdot|\!|\!|$, on
$\overline{{\M}_{|\!|\!|\cdot|\!|\!|}}$.
\end{Corollary}

\subsection{$\widetilde{\M}$ and $L^1(\M,\tau)$}
The following theorem is due to Nelson~\cite{Ne}.
\begin{Theorem}\label{T:Nelson}\emph{(Nelson, \cite{Ne})}
$\widetilde{\M}$ is a $\ast$-algebra and $T\in \widetilde{\M}$ if
and only if $T$ is a closed, densely defined operator affiliated
with $\M$. Furthermore, if $T\in \widetilde{\M}$ is a positive
operator, then $\lim_{n\rightarrow\infty}\chi_{[0,n]}(T)=T$ in the
measure-topology.
\end{Theorem}

  In the
following, we define $s$-numbers for unbounded operators in
$\widetilde{\M}$ as in~\cite{F-K}.
\begin{Definition}\emph{ For $T\in \widetilde{\M}$ and $0\leq s\leq 1$, define
the \emph{$s$-th numbers} of $T$ by
\[\mu_s(T)=\inf\{\|TE\|:\,\,\text{$E\in \M$ is a projection such that $\tau(E^\perp)=s$
}\}.
\]}
\end{Definition}

\begin{Theorem}\emph{(Fack and  Kosaki, \cite{F-K})}
\label{T:THIERRY FACK AND HIDEKI KOSAKI} Let $T$ and $T_n$ be a
sequence of operators in $\widetilde{\M}$ such that
$\lim_{n\rightarrow \infty}T_n=T$ in the measure-topology. Then for
almost all $s\in [0,1]$,
$\lim_{n\rightarrow\infty}\mu_s(T_n)=\mu_s(T)$.
\end{Theorem}

Let $\{T_n\}$ be a sequence of operators in $\M$ such that
$T=\lim_{n\rightarrow\infty} T_n$ in the $L^1$-norm. By
Lemma~\ref{L:1-norm},  $\{\tau(T_n)\}$ is a Cauchy sequence in
$\cc$. Define $\tau(T)=\lim_{n\rightarrow\infty}\tau(T_n)$. It is
obvious that $\tau(T)$ does not depend on the sequence $\{T_n\}$. In
this way, $\tau$ is extended to a linear functional on
$L^1(\M,\tau)$.

\begin{Lemma}\label{L:TX in L1} Let $|\!|\!|\cdot|\!|\!|$ be a normalized
unitarily invariant norm on a type ${\rm II}\sb 1$ factor $\M$. If
 $T\in \overline{{\M}_{|\!|\!|\cdot|\!|\!|}}$ and $X\in \M$, then $TX\in
 L^1(\M,\tau)$.
\end{Lemma}
\begin{proof}By the proof Proposition~\ref{P:embedding into measure topology},
$\lim_{n\rightarrow\infty}T_n=T$ in the measure topology. Hence
$\lim_{n\rightarrow\infty}T_nX=TX$ in the measure topology (see
Theorem 1 of~\cite{Ne}). By Corollary~\ref{C:Holder inequality for
bounded operators}, $\|T_nX-T_mX\|_1\leq |\!|\!|T_n-T_m|\!|\!|\cdot
|\!|\!|X|\!|\!|^\#$. So $\{T_nX\}$ is a Cauchy sequence in the
$L^1$-norm. By Corollary~\ref{C:embedding into measure topology},
$TX\in L^1(\M,\tau)$ and $\lim_{n\rightarrow\infty}T_nX=TX$ in the
$L^1$-norm.
\end{proof}

\subsection{Elements in $\overline{{\M}_{|\!|\!|\cdot|\!|\!|}}$}

\begin{Lemma}\label{L:dual norms for unbounded operators} For all $T\in \overline{{\M}_{|\!|\!|\cdot|\!|\!|}}$,
\[|\!|\!|T|\!|\!|=\sup\{|\tau(TX)|:\, X\in\M,\, |\!|\!|X|\!|\!|^\#\leq 1\}.
\]
\end{Lemma}
\begin{proof}Let $\{T_n\}$ be a sequence of operators in $\M$ such
that $\lim_{n\rightarrow\infty} T_n=T$ with respect to
$|\!|\!|\cdot|\!|\!|$. By Corollary~\ref{C:Holder inequality for
bounded operators}, if $X\in \M$ and $|\!|\!|X|\!|\!|^\#\leq 1$,
then $|\tau(TX)|=\lim_{n\rightarrow \infty}|\tau(T_nX)|\leq
\lim_{n\rightarrow \infty} |\!|\!|T_n|\!|\!|=|\!|\!|T|\!|\!|$.
Therefore, $|\!|\!|T|\!|\!|\geq \sup\{|\tau(TX)|:\, X\in\M,\,
|\!|\!|X|\!|\!|^\#\leq 1\}.$\\

We need to prove that $|\!|\!|T|\!|\!|\leq \sup\{|\tau(TX)|:\,
X\in\M,\, |\!|\!|X|\!|\!|^\#\leq 1\}.$ Let $\epsilon>0$. Since
$\lim_{n\rightarrow\infty} T_n=T$ with respect to
$|\!|\!|\cdot|\!|\!|$, there is an $N$ such that
$|\!|\!|T-T_N|\!|\!|<\epsilon/3$. For $T_N$, there is an $X\in\M$,
$|\!|\!|X|\!|\!|^\#\leq 1$, such that $|\!|\!|T_N|\!|\!|\leq
|\tau(T_NX)|+\epsilon/3$. By the proof Lemma~\ref{L:TX in L1} and
Corollary~\ref{C:Holder inequality for bounded operators},
\[|\tau(TX)-\tau(T_NX)|=\lim_{n\rightarrow\infty}|\tau(T_nX)-\tau(T_NX)|
\leq \lim_{n\rightarrow\infty}|\!|\!|T_n-T_N|\!|\!|\cdot
|\!|\!|X|\!|\!|^\#\leq |\!|\!|T-T_N|\!|\!| <\epsilon/3.
\]
  So $|\tau(TX)|\geq
|\tau(T_NX)|-|\tau((T_N-T)X)|\geq
|\!|\!|T_N|\!|\!|-\epsilon/3-\epsilon/3\geq
|\!|\!|T|\!|\!|-\epsilon$. Therefore, $|\!|\!|T|\!|\!|\leq
\sup\{|\tau(TX)|:\, X\in\M,\, |\!|\!|X|\!|\!|^\#\leq 1\}.$
\end{proof}

The following theorem generalizes Theorem~{\bf A}. Its proof   is
based on Lemma~\ref{L:dual norms for unbounded operators} and is
similar to the proof of Theorem~{\bf A}. So we omit the proof.
\begin{Theorem}\label{T: a representation theorem 2}  If
 $|\!|\!|\cdot|\!|\!|$ is a unitarily invariant norm on a type
 ${\rm II}\sb 1$ factor
$\M$, then there is a subset $\F'$ of $\F$ containing the constant 1
function on $[0,1]$ such that for all $T\in
\overline{{\M}_{|\!|\!|\cdot|\!|\!|}}$,
\[|\!|\!|T|\!|\!|=\sup\{|\!|\!|T|\!|\!|_{f}:\, f\in \F'\},\]
where $|\!|\!|T|\!|\!|_{f}$ is defined in Lemma~\ref{L:T_f} by
equation $(3)$ or by equation $(4)$ and
$\F=\{f(x)=a_1\chi_{[0,\frac{1}{n})}(x)+a_2\chi_{[\frac{1}{n},\frac{2}{n})}(x)+\cdots+
 a_n\chi_{[\frac{n-1}{n},1]}(x):\, a_1\geq a_2\geq \cdots\geq a_n\geq
0, \frac{a_1+\cdots+a_n}{n}\leq 1, n=1,2,\cdots\}$.
\end{Theorem}

Combining Theorem~\ref{T: a representation theorem 2} and
Theorem~\ref{T:THIERRY FACK AND HIDEKI KOSAKI}, we have the
following corollary.
\begin{Corollary}Let
 $|\!|\!|\cdot|\!|\!|$ be a unitarily invariant norm on a type
 ${\rm II}\sb 1$ factor
$\M$ and
 $|\!|\!|\cdot|\!|\!|'$ be the corresponding symmetric gauge norm on
 $(L^\infty[0,1],\int_0^1dx)$ as in Corollary~{\bf 2}.
  If $T\in \widetilde{\M}$, then $T\in
 \overline{{\M}_{|\!|\!|\cdot|\!|\!|}}$ if and only if
 $\mu_s(T)\in \overline{{L^\infty[0,\infty)}_{|\!|\!|\cdot|\!|\!|'}}$. In this case,
 $|\!|\!|T|\!|\!|=|\!|\!|\mu_s(T)|\!|\!|'$.
\end{Corollary}

\begin{Example}\emph{Let $T\in \widetilde{\M}$ and $1\leq p\leq \infty$. Then
$T\in L^p(\M,\tau)$ if and only if $\mu_s(T)\in L^p([0,1])$. In this
case,
$\|T\|_p=\left(\int_0^1\mu_s(T)^pds\right)^{1/p}=\left(\int_0^\infty
\lambda^p\,d\mu_{|T|}\right)^{1/p}$.}
\end{Example}

\subsection{A generalization of H$\ddot{\text{o}}$lder's inequality}
\begin{Lemma}\label{L:Tn converges to T} Let
$|\!|\!|\cdot|\!|\!|$ be a  unitarily invariant norm on a finite
factor $\M$ and
 $T\in \overline{{\M}_{|\!|\!|\cdot|\!|\!|}}$ be a positive
 operator. Then $\lim_{n\rightarrow\infty}\chi_{[0,n]}(T)=T$ with
 respect to $|\!|\!|\cdot|\!|\!|$.
\end{Lemma}
\begin{proof} If $|\!|\!|\cdot|\!|\!|$ is singular, then $T\in \M$
by Lemma~\ref{L:singular norms} and the lemma is obvious. We may
assume that $|\!|\!|\cdot|\!|\!|$ is continuous. Let
$T_n=\chi_{[0,n]}(T))$ and $\epsilon>0$. By
Lemma~\ref{L:epsilon-delta for a sequence operators}, there is a
$\delta>0$ such that if $\tau(E)<\delta$ then
$|\!|\!|TE|\!|\!|<\epsilon$. There is an $N$ such that
$\mu_s([N,\infty))<\delta$. So for $m>n\geq N$,
$|\!|\!|T_m-T_n|\!|\!|=|\!|\!|T\cdot\chi_{(m,n]}(T)|\!|\!|<\epsilon$.
This implies that $\{T_n\}$ is a Cauchy sequence of $\M$ with
respect to $|\!|\!|\cdot|\!|\!|$. Since
$\lim_{n\rightarrow\infty}T_n=T$ in the measure topology, by
Corollary~\ref{C:embedding into measure topology},
$\lim_{n\rightarrow\infty}T_n=T$ in the topology induced by
$|\!|\!|\cdot|\!|\!|$.
\end{proof}

The following theorem is a generalization of
H$\ddot{\text{o}}$lder's inequality.
\begin{Theorem}\label{T:Holder inequality} Let
$|\!|\!|\cdot|\!|\!|$ be a normalized unitarily invariant norm on a
finite factor $\M$. If
 $T\in \overline{{\M}_{|\!|\!|\cdot|\!|\!|}}$ and $S\in
 \overline{{\M}_{|\!|\!|\cdot|\!|\!|^\#}}$,
  then $TS\in L^1(\M,\tau)$ and $\|TS\|_1\leq
|\!|\!|T|\!|\!| \cdot |\!|\!|S|\!|\!|^\#$.
\end{Theorem}
\begin{proof} By the polar decomposition and Corollary~\ref{C:extension of norm on the completion
algebra}, we may assume that $S$ and $T$ are positive operators. Let
$T_n=\chi_{[0,n]}(T)$ and $S_n=\chi_{[0,n]}(S)$. By Lemma~\ref{L:Tn
converges to T},
$\lim_{n\rightarrow\infty}|\!|\!|T-T_n|\!|\!|=\lim_{n\rightarrow\infty}|\!|\!|S-S_n|\!|\!|^\#=0$.
Let $K$ be a positive number such that $|\!|\!|T_n|\!|\!|\leq K$ and
$|\!|\!|S_n|\!|\!|^\#\leq K$ for all $n$ and $\epsilon>0$. Then
there is an $N$ such  that for all $m>n\geq N$,
$|\!|\!|T_m-T_n|\!|\!|<\epsilon/(2K)$ and
$|\!|\!|S_m-S_n|\!|\!|^\#<\epsilon/(2K)$. By Corollary~\ref{C:Holder
inequality for bounded operators},
 $\|T_mS_m-T_nS_n\|_1\leq
\|(T_m-T_n)S_m\|_1+\|T_n(S_m-S_n)\|_1\leq |\!|\!|T_m-T_n|\!|\!|\cdot
|\!|\!|S_m|\!|\!|^\#+|\!|\!|T_n|\!|\!|\cdot
|\!|\!|S_m-S_n|\!|\!|^\#<\epsilon.$ This implies that $\{T_nS_n\}$
is a Cauchy sequence in $\M$ with respect to $\|\cdot\|_1$. Since
$\lim_{n\rightarrow\infty}T_nS_n=TS$ in the measure topology, by
Proposition~\ref{P:embedding into measure topology},
$\lim_{n\rightarrow\infty} T_nS_n=TS$ in $\|\cdot\|_1$.  By
Corollary~\ref{C:Holder inequality for bounded operators},
$\|T_nS_n\|_1\leq |\!|\!|T_n|\!|\!| \cdot |\!|\!|S_n|\!|\!|^\#$ for
every $n$. Hence, $\|TS\|_1\leq |\!|\!|T|\!|\!| \cdot
|\!|\!|S|\!|\!|^\#$.
\end{proof}

Combining Example~\ref{E:dual of Lp space} and Theorem~\ref{T:Holder
inequality}, we obtain the noncommutative H$\ddot{\text{o}}$lder's
inequality.
\begin{Corollary}Let $\M$ be a finite factor with the faithful
normal tracial state $\tau$. If $T\in L^p(\M,\tau)$ and $S\in
L^q(\M,\tau)$, then $TS\in L^1(\M,\tau)$ and
\[\|TS\|_1\leq \|T\|_p\cdot\|S\|_q,
\] where $1\leq p,q\leq \infty$ and $\frac{1}{p}+\frac{1}{q}=1$.
\end{Corollary}
\section{Proof of Theorem H and Theorem I}

 In this section, we assume that $\M$
is a type ${\rm II}\sb 1$ factor with the unique tracial state
$\tau$, $|\!|\!|\cdot|\!|\!|$ is a unitarily invariant norm on
$\M$ and $|\!|\!|\cdot|\!|\!|^\#$ is the dual unitarily invariant
norm on $\M$ (see Definition~\ref{D:dual norm}). Let
$|\!|\!|\cdot|\!|\!|_1$ be the corresponding symmetric gauge norm
on $(L^\infty[0,1],\int_0^1dx)$ as in Theorem~{\bf D} and
$|\!|\!|\cdot|\!|\!|_1^\#$ be the dual norm on
$(L^\infty[0,1],\int_0^1dx)$.\\

\begin{Lemma}\label{L:dual spaces on 21 factors and abelian} If $\overline{{\M}_{|\!|\!|\cdot|\!|\!|^\#}}$ is the
dual space of $\overline{{\M}_{|\!|\!|\cdot|\!|\!|}}$ in the sense
of question 1, then
$\overline{{L^\infty[0,1]}_{|\!|\!|\cdot|\!|\!|_1^\#}}$ is the dual
space of $\overline{{L^\infty[0,1]}_{|\!|\!|\cdot|\!|\!|_1}}$ in the
sense of question 1.
\end{Lemma}
\begin{proof} By Corollary~{\bf 2} and Lemma~\ref{L:isomorphism}, there is a separable diffuse
abelian von Neumann subalgebra $\A$ of $\M$ and a
$\ast$-isomorphism $\alpha$ from $\A$ onto $L^\infty[0,1]$ such
that $\tau=\int_0^1dx\circ \alpha$ and
$|\!|\!|\alpha(T)|\!|\!|_1=|\!|\!|T|\!|\!|$ for each $T\in \A$. By
Theorem~{\bf E}, $|\!|\!|\alpha(T)|\!|\!|_1^\#=|\!|\!|T|\!|\!|^\#$
for each $T\in \A$. So we need only to prove that
$\overline{{\A}_{|\!|\!|\cdot|\!|\!|^\#}}$ is the dual space of
$\overline{{\A}_{|\!|\!|\cdot|\!|\!|}}$ in the sense of question
1. Let $\phi\in \overline{{\A}_{|\!|\!|\cdot|\!|\!|}}^\#$. By the
Hahn-Banach extension theorem, $\phi$ can be extended to a bounded
linear functional $\psi$ on
$\overline{{\M}_{|\!|\!|\cdot|\!|\!|}}$ such that
$\|\psi\|=\|\phi\|$. By the assumption of the lemma, there is an
operator $X\in \overline{{\M}_{|\!|\!|\cdot|\!|\!|}^\#}$ such that
$\psi(S)=\tau(SX)$ for all $S\in
\overline{{\M}_{|\!|\!|\cdot|\!|\!|}}$ and
$\|\psi\|=|\!|\!|X|\!|\!|^\#$. Let $X=U|X|$ be the polar
decomposition of $X$ and $X_n=U\cdot\chi_{[0,n]}(|X|)$. By
Lemma~\ref{L:Tn converges to T}, $\lim_{n\rightarrow\infty}X_n=X$
with respect to the norm $|\!|\!|\cdot|\!|\!|^\#$. Let
$Y_n=\mathbf{E}_{\A}(X_n)$ for $n=1,2,\cdots$. By Corollary~{\bf
1}, $\{Y_n\}$ is a Cauchy sequence in $\A$ with respect to the
norm $|\!|\!|\cdot|\!|\!|^\#$ and $|\!|\!|Y_n|\!|\!|^\#\leq
|\!|\!|X_n|\!|\!|^\#$. Let $Y=\lim_{n\rightarrow \infty}Y_n$ with
respect to the norm $|\!|\!|\cdot|\!|\!|^\#$. Then $Y\in
\overline{{\A}_{|\!|\!|\cdot|\!|\!|}^\#}$ and
$|\!|\!|Y|\!|\!|^\#\leq |\!|\!|X|\!|\!|^\#=\|\psi\|=\|\phi\|$. For
$T\in \overline{{\A}_{|\!|\!|\cdot|\!|\!|}}$,
$\phi(T)=\psi(T)=\tau(TX)=\lim_{n\rightarrow\infty}\tau(TX_n)=
\lim_{n\rightarrow\infty}\tau(\mathbf{E}_\A(TX_n))=\lim_{n\rightarrow\infty}\tau(TY_n)=\tau(TY)$.
By Lemma~\ref{L:dual norms for unbounded operators},
$\|\phi\|=|\!|\!|Y|\!|\!|^\#.$
\end{proof}

Recall that $|\!|\!|\cdot|\!|\!|$ is a singular norm on $\M$ if
$\lim_{\tau(E)\rightarrow 0+}|\!|\!|E|\!|\!|>0$ and is a continuous
norm on $\M$ if $\lim_{\tau(E)\rightarrow 0+}|\!|\!|E|\!|\!|=0$ (see
section 11).

\begin{Corollary}\label{C:Dual space of singular norms} If $|\!|\!|\cdot|\!|\!|$ is a
singular unitarily invariant norm on $\M$, then
$\overline{{\M}_{|\!|\!|\cdot|\!|\!|^\#}}$ is not the dual space of
$\overline{{\M}_{|\!|\!|\cdot|\!|\!|}}$ in the sense of question 1.
\end{Corollary}
\begin{proof}
   Since
$|\!|\!|\cdot|\!|\!|$ is a singular norm on $\M$, by
Lemma~\ref{L:singular norms}, $|\!|\!|\cdot|\!|\!|$ is equivalent to
the operator norm  on $\M$ and
$\overline{{\M}_{|\!|\!|\cdot|\!|\!|}}=\M$. By
Corollary~\ref{C:equivalent dual norms} and Theorem~\ref{T:dual
norms of Ky Fan norms}, $|\!|\!|\cdot|\!|\!|^\#$ is equivalent to
the $L^1$-norm on $\M$. So $|\!|\!|\cdot|\!|\!|_1$ is equivalent to
the $L^\infty$-norm on $L^\infty[0,1]$ and
$|\!|\!|\cdot|\!|\!|_1^\#$ is equivalent to the $L^1$-norm on
$L^\infty[0,1]$ by Theorem~{\bf E}. Note that
$\overline{{L^\infty[0,1]}_{|\!|\!|\cdot|\!|\!|_1}}=L^\infty[0,1]$
is not separable with respect to $|\!|\!|\cdot|\!|\!|_1$ but
$\overline{{L^\infty[0,1]}_{|\!|\!|\cdot|\!|\!|_1^\#}}$ is separable
with respect to $|\!|\!|\cdot|\!|\!|_1^\#$. So
$\overline{{L^\infty[0,1]}_{|\!|\!|\cdot|\!|\!|_1^\#}}$ is not the
dual space of $\overline{{L^\infty[0,1]}_{|\!|\!|\cdot|\!|\!|_1}}$
in the sense of question 1. By Lemma~\ref{L:dual spaces on 21
factors and abelian}, $\overline{{\M}_{|\!|\!|\cdot|\!|\!|^\#}}$ is
not the dual space of $\overline{{\M}_{|\!|\!|\cdot|\!|\!|}}$ in the
sense of question 1.
\end{proof}

\begin{Lemma}\label{L:Dual space for continuous norms} If $|\!|\!|\cdot|\!|\!|$ is a
continuous unitarily invariant norm on $\M$, then
$\overline{{\M}_{|\!|\!|\cdot|\!|\!|^\#}}$ is  the dual space of
$\overline{{\M}_{|\!|\!|\cdot|\!|\!|}}$ in the sense of question 1.
\end{Lemma}
\begin{proof}We may assume that $|\!|\!|1|\!|\!|=1$.
By Theorem~\ref{T:Holder inequality},
$\overline{{\M}_{|\!|\!|\cdot|\!|\!|^\#}}$ is a subspace of the dual
space of $\overline{{\M}_{|\!|\!|\cdot|\!|\!|}}$ in the sense of
question 1.
 Let $\phi$ be a linear functional in the dual space of
$\overline{{\M}_{|\!|\!|\cdot|\!|\!|}}$. Then for every $T\in
\overline{{\M}_{|\!|\!|\cdot|\!|\!|}}$, $|\phi(T)|\leq \|\phi\|\cdot
|\!|\!|T|\!|\!|$. By Corollary~\ref{C:unitarily seminorm is norm},
for every $T\in \M$, $|\phi(T)|\leq \|\phi\|\cdot\|T\|$. So $\phi$
is a bounded linear functional on $\M$. Since $|\!|\!|\cdot|\!|\!|$
is a continuous norm on $\M$, $\lim_{\tau(E)\rightarrow
0}|\!|\!|E|\!|\!|=0$. Hence, $\lim_{\tau(E)\rightarrow 0}\phi(E)=0$.
This implies that $\phi$ is an ultraweakly continuous linear
functional on $\M$ and hence in the predual space of $\M$. So there
is an operator $X\in L^1(\M,\tau)$ such that for all $T\in \M$,
$\phi(T)=\tau(TX)$. By Lemma~\ref{L:dual norms for unbounded
operators}, $|\!|\!|X|\!|\!|^\#=\|\phi\|<\infty$. This implies that
$X\in\overline{{\M}_{|\!|\!|\cdot|\!|\!|^\#}}$. So
$\phi(T)=\tau(TX)$ for all $T\in
\overline{{\M}_{|\!|\!|\cdot|\!|\!|}}$ and
$\|\phi\|=|\!|\!|X|\!|\!|^\#$. This proves the lemma.
\end{proof}

\begin{proof}[Proof of Theorem~{\bf H} and Theorem~{\bf I}]Combining Lemma~\ref{L:dual
spaces on 21 factors and abelian}, \ref{L:Dual space for continuous
norms} and Theorem~{\bf A} gives the proof of Theorem~{\bf H} and
Theorem~{\bf I}.
\end{proof}

\begin{Example}\emph{If $1\leq p<\infty$ and $\frac{1}{p}+\frac{1}{q}=1$,
then $L^q(\M,\tau)$ is the dual space of $L^p(\M,\tau)$.
$L^1(\M,\tau)$ is not the dual space of $\M$.}
\end{Example}

\begin{Example}\emph{ For $1<p<\infty$, $L^p(\M,\tau)$ is a reflexive
space. $L^1(\M,\tau)$ and $\M$ are not reflexive spaces.  By
Theorem~\ref{T:dual norms of Ky Fan norms}, for $0\leq t\leq 1$,
$\overline{{\M}_{|\!|\!|\cdot|\!|\!|_{(t)}}}$ is not a reflexive
space.}
\end{Example}

\end{document}